\documentclass[sn-mathphys,Numbered]{sn-jnl}

\usepackage{graphicx}%
\usepackage{multirow}%
\usepackage{amsmath,amssymb,amsfonts}%
\usepackage{amsthm}%
\usepackage{mathrsfs}%
\usepackage[title]{appendix}%
\usepackage{xcolor}%
\usepackage{textcomp}%
\usepackage{manyfoot}%
\usepackage{booktabs}%
\usepackage{algorithm}%
\usepackage{algorithmicx}%
\usepackage{algpseudocode}%
\usepackage{listings}%
\usepackage{caption}
\usepackage{subcaption}

\newcommand{\R}{\mathbb{R}}
\newcommand{\norm}[1]{||{#1}||}

\renewcommand{\norm}[2][]{\left\Vert#2\right\Vert_{#1}}

\newcommand{\kron}{\otimes}
\DeclareMathOperator*{\argmin}{arg\,min}                   
\renewcommand{\t} {^{\top}}                                

\renewcommand{\phi}{\mathbf{\varphi}}

\newcommand{\bfzero}{{\bf0}}

\newcommand{\bfA}{\mathbf{A}}

\newcommand{\bfU}{\mathbf{U}}

\newcommand{\bfR}{\mathbf{R}}
\newcommand{\bfI}{\mathbf{I}}
\newcommand{\bfD}{\mathbf{D}}

\newcommand{\bfb}{\mathbf{b}}

\newcommand{\bfH}{\mathbf{H}}

\newcommand{\bfx}{\mathbf{x}}

\newcommand{\bfe}{\mathbf{e}}
\newcommand{\bfu}{\mathbf{u}}
\newcommand{\bfy}{\mathbf{y}}

\newcommand{\bfM}{\mathbf{M}}

\newcommand{\bfL}{\mathbf{L}}

\newcommand{\bfW}{\mathbf{W}}
\newcommand{\bfz}{\mathbf{z}}

\newcommand{\bfs}{\mathbf{s}}
\newcommand{\bfv}{\mathbf{v}}
\newcommand{\bfV}{\mathbf{V}}

\newcommand{\bfQ}{\mathbf{Q}}

\newcommand{\bfS}{\mathbf{S}}
\newcommand{\bfZ}{\mathbf{Z}}

\newcommand{\calN}{\mathcal{N}}

\newcommand{\bfPsi}{{\boldsymbol{\Psi}}}

\newcommand{\bflambda}{{\boldsymbol{\lambda}}}

\newcommand{\bfxi}{{\boldsymbol{\xi}}}

\newcommand{\bfpsi}{{\boldsymbol{\psi}}}

\newcommand{\bbR}{\mathbb{R}}

\begin{document}

\title[Flexible Krylov Methods for Group Sparsity Regularization]{Flexible Krylov Methods for Group Sparsity Regularization}


\author[1]{\fnm{Julianne} \sur{Chung}}

\author*[1]{\fnm{Malena} \sur{Sabat\'e Landman}}\email{malena.sabate.landman@emory.edu}


\affil[1]{\orgdiv{Department of Mathematics}, \orgname{Emory University}, \orgaddress{ \city{Atlanta}, \postcode{30322}, \state{GA}, \country{USA}}}




\abstract{This paper introduces new solvers for efficiently computing solutions to large-scale inverse problems with group sparsity regularization, including both non-overlapping and overlapping groups.
Group sparsity regularization refers to a type of structured sparsity regularization, where the goal is to impose additional structure in the regularization process by assigning variables to predefined groups that may represent graph or network structures.
Special cases of group sparsity regularization include $\ell_1$ and isotropic total variation regularization. In this work, we develop hybrid projection methods based on flexible Krylov subspaces, where we first recast the group sparsity regularization term as a sequence of 2-norm penalization terms using adaptive regularization matrices in an iterative reweighted norm fashion. Then we exploit flexible preconditioning techniques to efficiently incorporate the weight updates.  The main advantages of these methods are that they are computationally efficient (leveraging the advantages of flexible methods), they are general (and therefore very easily adaptable to new regularization term choices), and they are able to select the regularization parameters automatically and adaptively (exploiting the advantages of hybrid methods). Extensions to multiple regularization terms and solution decomposition frameworks (e.g., for anomaly detection) are described, and a variety of numerical examples demonstrate both the efficiency and accuracy of the proposed approaches compared to existing solvers.}



\maketitle

\section{Introduction}
Large-scale linear ill-posed inverse problems of the form
\begin{equation}\label{eq:linear_pbm}
\bfA \bfx_{\rm true} + \bfe = \bfb,
\end{equation}
arise in the discretization of problems coming from various scientific and engineering applications, such as biomedical, atmospheric and medical imaging \cite{Hansen2010,Vogel2002}. In this form, $\bfx_{\rm true}$ is the unknown solution, $\bfb$ is the measured data affected by unknown Gaussian white noise $\bfe$, and $\bfA\in\R^{m\times n}$ models the forward model. Given $\bfb$ and $\bfA$, the goal of the inverse problem is to approximate $\bfx_{\rm true}$, but computing solutions can be challenging especially for large-scale and ill-posed problems. More specifically, we are interested in cases where $\bfA$ is ill-conditioned and has ill-determined rank, meaning that the singular values of $\bfA$ decay and cluster at zero without a clear gap to indicate numerical rank. Due to the presence of noise in the measured data, the naive solution $\bfA^{\dagger}\bfb$ (where $\bfA^{\dagger}$ is the Moore-Penrose pseudoinverse of $\bfA$) can differ significantly from the desired true solution $\bfA^{\dagger}(\bfb-\bfe)$. This is due to noise amplification, as described in \cite{Hansen2010}. To obtain a meaningful approximation of the true solution, regularization is required.  There are many forms of regularization, but the general idea of regularization is to include prior knowledge about the solution (e.g., smoothness assumptions, solution structure, or hard constraints) in the solution process.  

In this work, we are interested in group sparsity regularization, which is a type of structured sparsity regularization, that promotes sparsity among pre-defined groups.  Consider the general optimization framework where the goal is to solve the following variational regularization problem,
\begin{equation}\label{eqn:problem_gs_1}
\min_{\bfx} {\{\|\bfA\bfx-\bfb\|}^{2}_{2}+\lambda{\|\bfPsi \bfx\|}_{2,1}\},
\end{equation}
where $\lambda>0$ is a regularization parameter, $\bfPsi$ is an invertible matrix, and the $\ell{}_{2,1}$ norm for a vector $\bfz$ whose elements belong to $s$ different groups is defined as 
\begin{equation}\label{eq:gs_def}
    {\|\bfz\|}_{2,1}=\sum_{i=1}^{s}{\|\bfz_{g_i}\|}_{2},
\end{equation}
with the vector $\bfz_{g_i}$ denoting all components in $\bfz$ belonging to the $i$-th group.  The choice of groups promotes structural patterns in the unknown parameters, and the methodology presented in this paper can handle both non-overlapping and overlapping groups.

Group sparsity regularization has gained significant interest in the literature, due to recent applications to machine learning (e.g., in training of deep neural networks \cite{scardapane2017group} and spatio-temporal problems \cite{pasha2021efficient}) and new theories (e.g., for strong group sparsity \cite{huang2017some}).  The benefits of using group sparsity regularization have been studied and demonstrated on various problems; however, various computational challenges remain that have hindered the use and extension of these group sparsity regularizers for large-scale problems.  First, selecting an appropriate choice of the regularization parameter $\lambda$ can be a difficult task, and most solution approaches require time-consuming tuning of this parameter.  Second, for problems where the forward model matrix $\bfA$ lacks an exploitable structure and may not even be explicitly constructed and stored, the only viable approach to solving problem \eqref{eq:linear_pbm} is to employ iterative methods.  Iterative methods exhibit semiconvergence behavior, where solutions at later iterations become dominated by errors, so a good stopping criteria is critical \cite{jensen2007iterative}.  Third, in many scenarios, group sparsity on its own is not enough, and further refinement is needed.  For example, regularizers that combine a group sparsity prior as well as an $\ell_1$ (LASSO) regularizer \cite{chartrand2013nonconvex,scardapane2017group} result in problems of the form,
\begin{equation}\label{eq:combinedproblem}
\min_{\bfx} {\{\|\bfA \bfx - \bfb\|}^{2}_{2}+\lambda \|\bfx\|_{1}+\alpha {\| \bfx\|}_{2,1}\},
\end{equation} 
for some $\alpha>0.$
One interpretation is that in addition to the nonzero components being clustered into groups, the nonzero groups themselves may be sparse.

Anomaly detection presents another scenario that requires going beyond the standard group sparsity problem \eqref{eqn:problem_gs_1}.  In this case, the solution can be represented as $\bfx = \bfxi+ \bfs$, where $\bfxi$ captures background smoothness and $\bfs$ represents anomalous events \cite{chung2022hybrid}.  Such problems arise in  the context of atmospheric inverse modeling, and these are often spatio-temporal problems where anomalies are consistent over time but sparse in terms of spatial location.  Thus, group sparsity regularization can be used for $\bfs$ with groups defined via time, while a spatio-temporal Gaussian prior can be used for $\bfxi,$ e.g., $\bfxi \sim \calN(\bfzero, \alpha^{-1}\bfQ)$, where $\bfQ$ is a symmetric positive definite (SPD) matrix.  Assuming that $\bfR$ is an SPD matrix representing the noise covariance matrix, we are interested in anomaly detection problems of the form, 
\begin{equation}\label{eq:sdproblem}
\min_{\bfxi, \bfs} {\{\|\bfA(\bfxi+\bfs)-\bfb\|}^{2}_{\bfR^{-1}}+\alpha \|\bfxi\|_{\bfQ^{-1}}^2+\lambda{\|\bfs\|}_{2,1}\},
\end{equation}
where $\|\bfx\|_{\bfM}^2 = \bfx\t \bfM \bfx$ for SPD matrix $\bfM.$ 
For both \eqref{eq:combinedproblem} and \eqref{eq:sdproblem}, a further challenge is to estimate the additional regularization parameter $\alpha$.

\paragraph{Contributions and overview.} In this paper, we develop flexible Krylov methods, including hybrid and iteratively reweighted variants, for efficiently computing solutions to large-sale inverse problems with group sparsity regularization.  Although Krylov methods have been considered for for sparsity-promoting regularization terms \cite{gazzola2014generalized,chung2019flexible,gazzola2020krylov,6aff08faf77e4015b89d53a954ba37df}, this work is the first to use flexible methods to address group sparsity regularizers.  We address particular challenges that come with extensions to group sparsity regularizers, including handling overlapping groups, combining multiple regularization terms, and extending solution decomposition frameworks for anomaly detection.
The proposed methods have the advantage of being general (and therefore very easily adaptable to new regularization term choices),
very efficient (leveraging the advantages of flexible methods), and allowing for automatic and adaptive regularization parameter choices (exploiting the advantages of hybrid methods).
The performance of this framework is shown through a variety of numerical examples showcasing both the efficiency of the proposed framework compared to existing solvers and showing novel applications of group sparsity.

The paper is organized as follows.  In Section \ref{sec:background} we provide a some background on group sparsity. Then in Section \ref{sec:flexible} we describe how iterative reweighted norm schemes can be used to handle a group sparsity regularization term and
describe various flexible hybrid projection methods, with particular focus on extensions to multiple regularization terms and to the solution decomposed problem.  In  Section \ref{sec:numerics} we provide various examples from image deblurring, dynamic image deblurring and anomaly detection in atmospheric inverse modeling. Conclusions are provided in Section \ref{sec:conclusions}.

\section{Background on group sparsity}
\label{sec:background}
The group sparsity regularized problem \eqref{eqn:problem_gs_1} corresponds to a very general framework for regularization.  In fact, the choice of specific groups leads to standard regularization terms. For example, $\ell_1$ regularization can be obtained by considering $n$ groups, with each group containing a single element, and 2-norm regularization can be obtained by considering one group containing all parameters.  Note that this is not Tikhonov regularization, since there is no square on the regularizer. Moreover, isotropic total variation (TV) regularization in any dimension (e.g. in 2D, 3D, and considering a temporal dimension) can be obtained by taking $\bfPsi$ to be a discrete derivative operator and grouping all directional derivatives.  Notice that a similar case to the combined regularized problem \eqref{eq:combinedproblem} is the elastic net regularizer, which combines an $\ell_1$ and an $\ell_2$ regularizer \cite{zou2005regularization}; hence, the methods proposed in this paper to solve \eqref{eq:combinedproblem} can be straightforwardly adapted to solve the elastic net problem.

A major distinction of the group sparsity regularizer is the ability to select groups, where the choice of groups promotes structural patterns in the unknown parameters.  For many problems, there is a natural choice for the groups that comes from the structure of the problem.  For example, in source localization and anomaly detection where spatial-temporal images are desired, a natural grouping arises where each pixel at every time point constitutes a group. 
This has been studied, e.g. for image-based biochemical assay \cite{AguilaPla2017a}, low-dimensional nonlinear signal modeling \cite{COSTA2017142}, and classification \cite{4959720}.  In other applications, groups are defined by networks (e.g., deep neural networks) or by wavelet hierarchies \cite{6115845}.
 
The groups can be non-overlapping (e.g., \eqref{eqn:problem_gs} with $\bfPsi = \bfI$ corresponds to group LASSO), where a partition of the unknown vector is made, and the regularizer promotes sparsity for the entire group.  Or, the groups can be overlapping, where a component of the unknown vector can belong to more than one group simultaneously.  If the unknown parameters themselves do not lend themselves to natural grouping, an invertible matrix $\bfPsi$ in \eqref{eqn:problem_gs_1} can be used to arrange the groups. Notice that \eqref{eqn:problem_gs_1} can be expressed as
\begin{equation}\label{eqn:problem_gs}
\min_{\bfz} {\{\|\bfA \bfPsi^{-1} \bfz-\bfb\|}^{2}_{2}+\lambda{\|\bfz\|}_{2,1}\}, \quad \text{for}  \quad \bfx = \bfPsi^{-1} \bfz.
\end{equation}

\begin{figure}[htb]
     \centering
     \begin{subfigure}[b]{0.32\textwidth}
         \centering
         \includegraphics[trim = 6cm 4cm 5.5cm 3cm,clip,height=4cm]{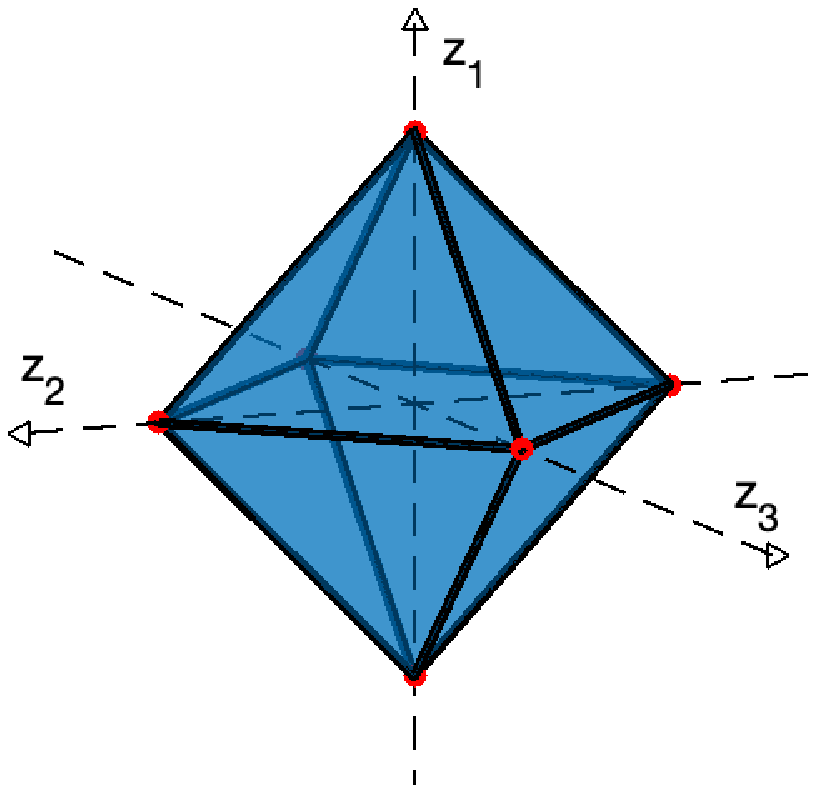}
         \caption{$\ell_{1}$ regularization}
         \label{fig:l1}
     \end{subfigure}
     \begin{subfigure}[b]{0.32\textwidth}
         \centering
         \includegraphics[trim = 6cm 4cm 5.5cm 3cm,clip,height=4cm]{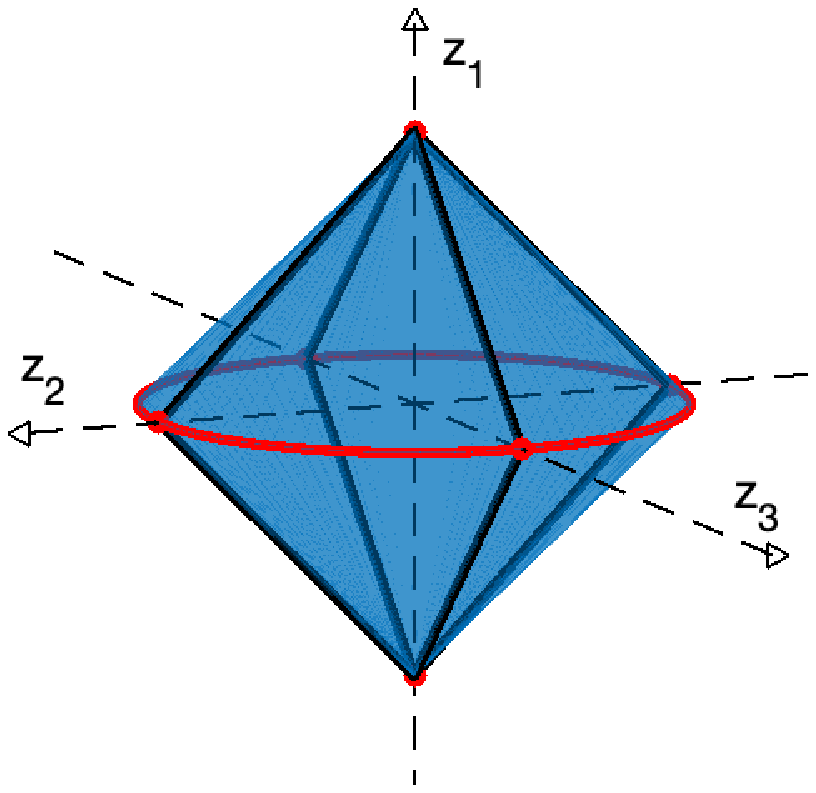}
         \caption{$\ell_{2,1}$ regularization}
         \label{fig:l21}
    \end{subfigure} 
    \begin{subfigure}[b]{0.32\textwidth}
         \centering
         \includegraphics[trim = 6cm 4cm 5.5cm 3cm,clip,height=4cm]{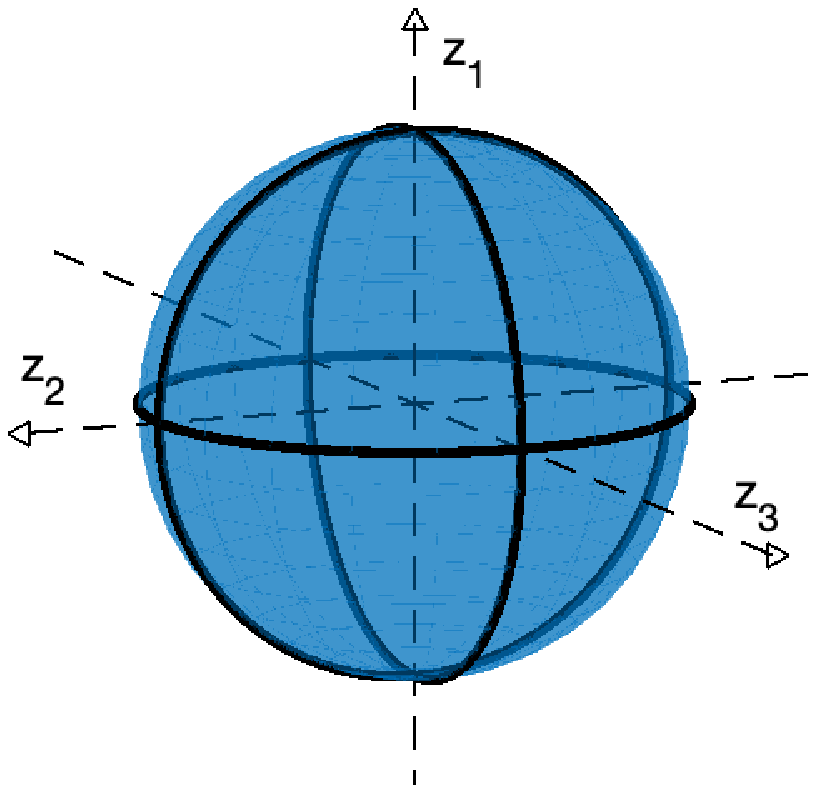}
         \caption{$\ell_{2}$ regularization}
         \label{fig:l2}
     \end{subfigure}
        \caption{Graphical illustration of the $\ell_{2,1}$ regularizer (with predefined groups $\bfz_{g_1}=[z_1]$ and $\bfz_{g_2}=[z_2, z_3]$) in (b) compared to the $\ell_{1}$ and $\ell_{2}$ norms in (a) and (c) respectively.}
        \label{fig:l112_graphic_explanation}
\end{figure}

A natural question is to how to interpret the group sparsity regularizer. In Figure \ref{fig:l112_graphic_explanation}, we provide a graphical representation of non-overlapping group sparsity for a schematic 3D example, compared to standard $\ell_1$ and $\ell_2$ norms. 
The edges of the solid shapes in Figures \ref{fig:l1}, \ref{fig:l21} and \ref{fig:l2}, represent the different vectors $\bfz=[z_1,z_2,z_3]\in \mathbb{R}^{3}$ with $\|\bfz\|_1=1$, $\|\bfz\|_{2,1}=1$ and $\|\bfz\|^2_2=1$, respectively, where the groups for the $\ell_{2,1}$ norm have been taken to be $\bfz_{g_1}=[z_1]$ and $\bfz_{g_2}=[z_2, z_3]$. The geometry of the norm ball is directly related to the solution of \eqref{eqn:problem_gs}. 
In particular, the positions of the singularities are different for each norm.
For example, the $\ell_1$ norm promotes sparsity in the individual coefficients.  To see this, imagine a plane corresponding to the set of solutions with equal discrepancy (i.e., giving the same values for the fit-to-data term $\|\bfA \bfPsi^{-1} \bfz-\bfb\|_2^2$).  The intersection of this plane with the octahedron in Figure \ref{fig:l1} will most likely occur at a corner (i.e., denoted by the red dots in Figure \ref{fig:l1}), which corresponds to a sparse solution.
On the contrary, the differentiable-everywhere $\ell_2$-norm which represents a sphere in $\bbR^3$ (see Figure \ref{fig:l2}) does not favor any particular direction (i.e., does not promote sparsity in the solution).  The $\ell_{2,1}$ norm is represented as a double-cone in Figure \ref{fig:l21}, and it provides a natural balance between the $\ell_1$ and $\ell_2$ norms. Since the two cones are attached at their bases, $\ell_{2,1}$ regularization promotes sparsity, but only at the group level, i.e. either $z_1$ or both $z_2$ and $z_3$.  Notice the absence of singularities at standard unit vectors $\bfe_2$ and $\bfe_3$, such that the likely points of intersection (e.g., with a plane) would occur at the points denoted in red.


For overlapping group sparsity, a graphical illustration is not as intuitive and depends on the choice of overlapping weights.  Nevertheless, overlapping group sparsity can be very useful in some applications (e.g. tree structured groups, contiguous groups, and directed-acyclic-graph groups).  It has been shown in \cite{Baraniuk2010} that structured sparsity assumptions can improve computational efficiency and accuracy in the solution in an analysis analogous to traditional compressive sensing for sparsity. This can be exploited, for example, when using tree-based image representations, see e.g., \cite{542157} for a characterization using the block-DCT transform associated with JPEG, or \cite{6115845} for image representations using wavelets. In this setting, sub-trees are formed by grouping hierarchically different scale representations at the same location, so for natural images it is reasonable to assume that the image would be zero at all resolutions in a given point; therefore, one can group parents and children together in different ways. For example, this is exploited in \cite{Chen2014}, where an algorithm based on preconditioned conjugate gradient with restarts is presented, and \cite{NIPS2012_65658fde}, where group-wise soft thresholding is used to deal with the $\ell_{2,1}$ regularization term.

Many of the challenges to computing solutions to \eqref{eqn:problem_gs_1} mimic those in LASSO and compressive sensing (e.g., a discontinuous first derivative at zero from the 1-norm, causing challenges for optimization algorithms).  Various methods have been proposed to solve such problems, ranging from convex relaxations (e.g., iteratively reweighted schemes) and nonconvex optimization methods (e.g., ADMM).  A main caveat of all of these methods is the need to select the regularization parameter a priori.   More recently, flexible Krylov methods have been shown to be powerful alternatives to solve $\ell_1$ regularization problems (see \cite{Gazzola21}, \cite{chung2019flexible}) compared to other classical solvers such as suitable versions of FISTA \cite{beck2009fast} or SpaRSA \cite{4518374}, or using iteratively reweighted norms (maybe) in combination with Krylov methods (see, e.g. \cite{https://doi.org/10.1002/cpa.20303}, \cite{rodriguez-2008-sparse}). Similarly, flexible Krylov methods have been used for TV regularization (see \cite{jimaging7100216},\cite{6aff08faf77e4015b89d53a954ba37df}) and proved competitive against the most popular solvers for this problem such as proximal gradient methods \cite{5173518}, primal–dual methods \cite{doi:10.1137/S1064827596299767}, split Bregman methods \cite{doi:10.1137/040605412} or methods based on iteratively reweighted norms \cite{4729670}. Another example of modified TV using group sparsity can be found in \cite{pasha2021efficient}, where group sparsity is used to promote piece-wise constant structures in space and time, sparsity in gradient images shared across time instances, but not among the group. This is solved using a majorization-minimization approach based on a generalized Krylov subspace approach 
\cite{bach2012optimization,herzog2012directional}. 
 For non-overlapping groups, accelerated proximal gradient methods have been considered with added non-negativity constraints to solve \eqref{eqn:problem_gs_1} in \cite{AguilaPla2017a} and alternating direction method of multipliers was used in \cite{chartrand2013nonconvex}. 
%

It is beyond the scope of this paper to handle cases where the groups are not know a-priori. However, the methodology presented in Section \ref{sec:flexible} could be used in combination with a strategy to determine (or modify) groups on-the-fly once a new approximation of the solution is available, since the weights are iteration dependent. Similarly, the flexible Krylov methods can be easily generalised to $\ell_{2,p}$ by using suitably defined weights. Moreover, if memory requirements are a problem, iteratively reweighted schemes with fixed preconditioning can be used in an inner-outer scheme.

\section{Flexible Krylov methods for group sparsity} 
\label{sec:flexible}
In this section, we describe flexible Krylov methods and hybrid variants for solving the group sparsity regularized problem \eqref{eqn:problem_gs_1} for both overlapping and non-overlapping groups, and we develop extensions to solve \eqref{eq:combinedproblem} and \eqref{eq:sdproblem}. We begin in Section \ref{sub:irn} with a description of the iteratively reweighted norm approach.  
  Then in Section \ref{sub:flexible}, we describe flexible Krylov methods for efficiently generating a single solution basis for approximately solving the re-weighted least-squares subproblems. A method for automatically selecting the regularization parameter within the hybrid framework is described.

\subsection{Iteratively reweighted norm schemes}
\label{sub:irn}
Traditional iteratively reweighted (IRW) schemes are majorization-minimization methods that involve constructing a sequence of least-squares problems approximating a variational regularization problem of interest, which are then solved at each iteration. These methods intrinsically rely on the interpretation of a given $\ell_p$ norm as a non-linear weighted $\ell_{2}$ norm (for $0 < p \leq 2$), such that evaluating the weights at an available approximation of the solution $\bfx_{k-1}$ can be used to obtain a quadratic tangent majorant of the original problem at $\bfx_{k-1}$. Recall that a quadratic tangent majorant of a non-quadratic function at a point $\bar{\bfz}$ is a quadratic upper bound of the original function constructed such that both the value of the function and the value of the gradient of both functions coincide at $\bar{\bfz}$. More specifically, for $\Phi(\bfz)=\frac{1}{p}\|\bfz\|^p_p=\frac{1}{p} \|W(\bfz)\bfz\|_2^2$, then $\Phi^{\bar{\bfz}}(\bfz)=\frac{1}{2}\|W(\bar{\bfz})\bfz\|_2^2+(\frac{1}{p}-\frac{1}{2})\|W(\bar{\bfz})\bar{\bfz}\|_2^2$ is a quadratic tangent majorant of $\Phi(\bfz)$ at $\bar{\bfz}$. This can be explicitly checked as $\Phi^{\bar{\bfz}}$ is quadratic, $\Phi^{\bar{\bfz}}(\bar{\bfz}) = \Phi(\bar{\bfz})$, the gradients are equal at $\bar{\bfz}$, i.e. $\nabla \Phi^{\bar{\bfz}}(\bar{\bfz}) = \nabla \Phi(\bar{\bfz})$, and $\Phi^{\bar{\bfz}}(\bfz)$ is an upper bound of $\Phi(\bfz)$ for all $\bfz$ (see, e.g. \cite{lplq}). 

IRW norm schemes can be used to solve minimization problems involving $\ell_{2,1}$ norm terms. Assume the components of $\bfz \in \mathbb{R}^n$ belong to $s$ different (possibly overlapping) groups and define $\bfz_{g_i}$ to be the vector constructed by concatenating the components in $\bfz$ belonging to the $i$-th group. Let $G_j$ be the indexes of the (possibly multiple) groups where the component $z_j$ of $\bfz$ belongs to and, with a small notation abuse,  let $z_j \in \bfz_{g_i}$ be the components of $\bfz$ that belong to the $i$-th group. Then, we can define a diagonal weighting matrix (as a function of $\bfz$) as
\begin{equation}\label{weights}
    \bfW_{j,j}(\bfz) = \sqrt{ \sum_{i \in G_j} \frac{1}{{\|\bfz_{g_i}\|}_2}},
\end{equation}
such that the $\ell_{2,1}$ norm defined in \eqref{eq:gs_def} can be written as a reweighted $\ell_2$ norm:
\begin{eqnarray}\label{weights_check}
    {\|\bfW(\bfz)\bfz\|}_{2}^{2} &=&  \sum_{j=1}^{n} \left(\sum_{i \in G_j} \frac{1}{{\|\bfz_{g_i}\|}_2}  z_{j}^2 \right)= \sum_{i=1}^{s} \left(\sum_{z_j \in \bfz_{g_i}} \frac{1}{{\|\bfz_{g_i}\|}_2}  z_{j}^2\right) \\ &=& \sum_{i=1}^{s} \frac{{\|\bfz_{g_i}\|}_2^2}{{\|\bfz_{g_i}\|}_2} = \sum_{i=1}^s {\|\bfz_{g_i}\|}_2 = {\|\bfz\|}_{2,1}. \nonumber
\end{eqnarray}
The solution-dependent non-linear weights defined in \eqref{weights_check}, can be then used to re-write the original regularization problem \eqref{eqn:problem_gs} with group sparsity regularization as
\begin{equation}\label{ls_weighted}
\min_{\bfz} \{{\|\bfA \bfPsi^{-1} \bfz-\bfb\|}^{2}_{2}+\lambda{\|\bfW(\bfz) \bfz\|}_{2}^2\} \quad \text{for} \quad \bfx = \bfPsi^{-1} \bfz.
\end{equation}
Using an IRW scheme corresponds to evaluating (and fixing) the weights \eqref{weights} at each (outer) iteration using the available approximation of the solution coming from the previous problem in the sequence. However, this can be numerically unstable due to divisions by zero when one has to evaluate the weights on a vector with groups whose components are all 0 valued. This situation is in fact expected, as we are assuming sparsity in the groups, and is caused by the lack of smoothness of the functional in \eqref{eqn:problem_gs} at vectors with 0 valued components.  Thus, we consider a smooth approximation of the original problem, and the following weights are considered instead of the ones defined in \eqref{weights}:
\begin{equation}\label{eq:smooth_W}
\widetilde{\bfW}_{j,j}(\bfz) = \sqrt{ \sum_{i \in G_j} \frac{1}{
    \sqrt{{\|\bfz_{g_i}\|}^2_2+\tau^2}
    } }  ,
\end{equation}
so that $\bfW_{k} = \widetilde{\bfW}(\bfz_{k-1}) \approx \bfW(\bfz_{k-1})$. This yields the following sequence of least-squares problems approximating \eqref{eqn:problem_gs}:
\begin{equation}\label{eq:reweighted}
        \min_\bfz \{{\|\bfA \bfPsi^{-1} \bfz-\bfb\|}_{2}^{2}+\lambda{\|\bfW_k \bfz \|}_{2}^{2}\} \quad \text{for} \quad \bfx = \bfPsi^{-1} \bfz.    
\end{equation}
Note that this minimization is equivalent to minimizing a sequence of quadratic tangent majorants of a smoothed version of \eqref{eqn:problem_gs}, where multiplicative constants have been absorbed by $\lambda$ and additive constants have been disregarded. Moreover, note that one could transform problem \eqref{eq:reweighted} to the following equivalent expression using a solution-dependent right preconditioner (or, equivalently, a suitable change of variables):
\begin{equation}\label{eq:reweighted_std}
\min_{\bfs} \{{\|\bfA \bfPsi^{-1} \bfW_k^{-1} \bfs-\bfb\|}^{2}_{2}+\lambda {\|\bfs \|}_{2}^{2}\}
\quad \text{for} \quad \bfx = \bfPsi^{-1} \bfz = \bfPsi^{-1} \bfW_k^{-1} \bfs.  \end{equation}
Expression \eqref{eq:reweighted_std} motivates the use of flexible Krylov methods which allow for iteration-dependent preconditioning, see e.g. \cite{chung2019flexible} \cite{Gazzola21} and references therein, and will be explained in detail the following section. Note that reweighting schemes can also be used to solve \eqref{eq:sdproblem} using \eqref{weights_check} and \eqref{eq:smooth_W} to derive analogous expressions to \eqref{ls_weighted}, \eqref{eq:reweighted} and \eqref{eq:reweighted_std}, see \cite{chung2022hybrid}.

Moreover, we can handle the minimization \eqref{eq:combinedproblem} with multiple regularization terms by defining two solution-dependent diagonal weight matrices $\bfW_1$ and $\bfW_2$ such that \eqref{eq:combinedproblem} can be written as the following reweighted problem:
\begin{equation}
\min_{\bfx} {\{\|\bfA \bfx - \bfb\|}^{2}_{2}+\lambda \|\bfW_1(\bfx) \bfx\|^{2}_{2}+\alpha{\|\bfW_2(\bfx)\bfx\|}^{2}_{2}\}
\end{equation}
or, equivalently, 
\begin{equation}
\min_{\bfx} \left\{ {\|\bfA \bfx - \bfb\|}^{2}_{2}+ \left\|\begin{bmatrix}\sqrt{\lambda} \bfW_1(\bfx) \\ \sqrt{\alpha} \bfW_2(\bfx) \end{bmatrix} \bfx\right\|_{2}^2 \right\}.
\end{equation}

One point of concern is the estimation of regularization parameters, which is difficult if we want to solve for both parameters.  An alternative is to define $\alpha = \tau_{\lambda}^2 \lambda$ for some constant $\tau_{\lambda} $. Then, let 
\begin{equation}\label{eq:tau_lambda}
\bfL = \begin{bmatrix}\sqrt{\lambda} \bfW_1(\bfx) \\ \sqrt{\alpha} \bfW_2(\bfx) \end{bmatrix} = \sqrt{\lambda} \begin{bmatrix}  \bfW_1(\bfx) \\ \tau_{\lambda} \, \bfW_2(\bfx) \end{bmatrix} =  \sqrt{\lambda}\, \bfQ^{\bfL} \bfD
\end{equation}
be the QR factorization of $\bfL$. Notice that since $\bfL$ is a concatenation of diagonal matrices which are invertible (by definition), $\bfL$ is guaranteed to have orthogonal columns and is full column rank. Thus, the diagonal elements of $\bfD$ can be easily obtained by taking the norms of the columns of $\bfL$. Next, notice that $\bfD$ is a diagonal matrix that depends on $\bfx$ and $\bfQ^{\bfL}$ contains orthonormal columns, so we get the IRW regularized problem,
\begin{equation}\label{flsqr-C}
\min_{\bfx} {\{\|\bfA \bfx - \bfb\|}^{2}_{2}+ \lambda \left\|\bfD(\bfx) \bfx\right\|_{2}^2 \}.
\end{equation} 
Flexible hybrid methods can be used to solve \eqref{flsqr-C}, and extensions to more than two regularization terms is straightforward. Numerical results demonstrating this approach are provided in Section \ref{sec:combined}.

\subsection{Flexible Krylov methods} 
\label{sub:flexible}
Flexible Krylov methods are a subclass of Krylov methods that allow for iteration-dependent right preconditioning. The methodology presented here can be used generally for any flexible Krylov method, and in particular we will present three different methods to showcase the use of these algorithms in different applications.  

One interpretation of flexible Krylov methods where the preconditioning is motivated by a reweighted norm in the regularization term is that some information about the solution is embedded in the solution space \cite{10.1016/j.cam.2005.10.038}. Moreover, this is achieved efficiently by building a single solution space, thereby avoiding inner-outer schemes.

Since the projected space already contains information about the (variational) regularization term, one can also set $\lambda=0$ in \eqref{eq:reweighted_std} and obtain a suitably regularized solution equipped with early stopping. Without appropriate stopping, however, the solution will approach the naive solution $\bfA^{\dagger}\bfb$. This is the original approach followed in flexible Krylov methods and will be denoted with the prefix (F). Another approach is to add regularization in the projected problem in order to avoid semi-convergence. With a suitable regularization parameter, which can be adaptively set at each iteration, this method delivers a good reconstruction. However, note that the regularization parameter for the projected problems might not be suitable for the original problem \eqref{eqn:problem_gs}. This can be understood since, in the limit, the solution obtained with this method will converge to the solution of the fit-to-data term in \eqref{eqn:problem_gs} with added Tikhonov regularization. This is usually noted with the prefix (hybrid-F) and will be the preferred method used throughout the paper. For a more detailed explanation, see \cite{Gazzola21}. Lastly, problem \eqref{eqn:problem_gs} can be projected onto the flexible Krylov subspace. Assuming no break-down has happened at $k=\min(m,n)$, one can extend the algorithm for $k>=\min(m,n)$ by updating the weights in \eqref{eqn:problem_gs}. In this case, the solution provided by a flexible method following this scheme will converge to the solution of the (smoothed version of) \eqref{eqn:problem_gs}. This comes at the cost of a QR factorization of a tall and skinny matrix at each iteration and will be denoted with the prefix (IRW-F). The convergence proof for this scheme can be found in \cite{Gazzola21}.

\subsubsection{hybrid-FGMRES} Given $\bfA \in \mathbb{R}^{n \times n}$, $\bfb\in \mathbb{R}^{n}$, iteration independent right preconditioner $\bfPsi^{-1} \in \mathbb{R}^{n \times n}$ and iteration dependent right preconditioning matrices  $\bfW_k^{-1} \in \mathbb{R}^{n \times n}$, the flexible Arnoldi method is a process that, at the $k$th iteration, constructs vectors $\bfv_{k+1}$ and $\bfz_k = {\bfW}_{k}^{-1} \bfv_k$ such that
\begin{equation}\label{Farnoldi}
    \bfA \bfPsi^{-1} \bfZ_k=\bfV_{k+1} \bfH_k ,
\end{equation}
where, if $\bfx_0=\bf0$, $\bfv_1 = \bfb/ \|\bfb\|_2$. Here $\bfH_k \in \mathbb{R}^{(k+1)\times k}$ is upper Hessenberg and $\bfV_k \in \mathbb{R}^{n \times k}$ has orthonormal columns.  A detailed algorithm for flexible Arnoldi can be found in e.g. \cite{chung2019flexible}.

Using hybrid flexible GMRES (hybrid-FGMRES) to promote group sparsity regularization involves projecting the least-squares problem $\|\bfA\bfx-\bfb\|_2^2$ onto a modified (flexible) Krylov subspace using the decomposition in \eqref{Farnoldi} and adding regularization. Then, the solution at each iteration $k$ can be computed as 
\begin{eqnarray}\label{eq:gmres_full}
       \bfx_k =\bfPsi^{-1} \bfz_k  = \bfPsi^{-1} \bfZ_k \bfy_k  \quad \text{where} \quad \bfy_k = \argmin_\bfy \{{\|\bfH_k \bfy-{\|\bfb\|}_{2} \bfe_1\|}_{2}^{2}+\lambda{\| \bfy \|}_{2}^{2}\}.
\end{eqnarray}
Note that the iteration dependent preconditioners are built using ${\bfW}_{k} = \widetilde{\bfW}(\bfZ_{k-1} \bfy_{k-1}) = \widetilde{\bfW}(\bfz_{k-1})$ as defined in \eqref{eq:smooth_W}, which depend on the solution computed at the previous iteration.

\subsubsection{hybrid-FLSQR} For rectangular matrices $\bfA \in \mathbb{R}^{m \times n}$, we can use a flexible LSQR approach, where for $\bfb\in \mathbb{R}^{m}$, an iteration independent right preconditioner $\bfPsi^{-1} \in \mathbb{R}^{n \times n}$ and iteration dependent right preconditioning matrices $\bfW_k^{-1} \in \mathbb{R}^{n \times n}$, the flexible Golub-Kahan decomposition computes at the $k$th iteration vectors $\bfu_{k+1},$ $\bfv_k$ and $\bfz_k = {\bfW}_{k}^{-1} \bfv_k$ such that
\begin{equation}\label{FGK}
    \bfA \bfPsi^{-1} \bfZ_k=\bfU_{k+1} \bfM_k \quad \text{and} \quad \bfPsi^{-\top} \bfA\t \bfU_{k+1}=\bfV_{k+1}\bfS_{k+1},
\end{equation}
where, if $\bfx_0=\bf0$, $\bfu_1 = \bfb/ \|\bfb\|_2$. 
Here $\bfM_k \in \mathbb{R}^{(k+1)\times k}$ is upper Hessenberg, $\bfS_{k+1} \in \mathbb{R}^{(k+1)\times (k+1)}$ is upper triangular and both $\bfV_{k+1} \in \mathbb{R}^{n\times (k+1)}$ and $\bfU_{k+1} \in \mathbb{R}^{m\times (k+1)}$ contain orthonormal columns. For more details see, e.g. \cite{simoncini2007recent}.

Using hybrid flexible LSQR (hybrid-FLSQR) to promote group sparsity regularization involves projecting the least-squares problem $\|\bfA\bfx-\bfb\|_2^2$ onto a modified (flexible) Krylov subspace using the decomposition in \eqref{FGK} and adding regularization to the projected problem. Then, we can compute an approximation of the solution at each iteration $k$ as
\begin{eqnarray}\label{eq:flsqr_full}
        \bfx_k =\bfPsi^{-1} \bfz_k  = \bfPsi^{-1} \bfZ_k \bfy_k  \quad \text{where} \quad \bfy_k = \argmin_\bfy \{{\|\bfM_k \bfy-{\|\bfb\|}_{2} \bfe_1\|}_{2}^{2}+\lambda{\| \bfy \|}_{2}^{2}\}.
\end{eqnarray} 
In particular, we are interested in the case where ${\bfW}_{k} = \widetilde{\bfW}(\bfZ_{k-1} \bfy_{k-1}) = \widetilde{\bfW}(\bfz_{k-1})$ and therefore it depends on the solution computed at the previous iteration. 

\subsubsection{hybrid-SD} With appropriately defined weights, both hybrid-GMRES and hybrid-FLSQR can be used to solve group sparsity regularized problems \eqref{eqn:problem_gs_1} and \eqref{eq:combinedproblem}, but a different projection approach is needed to handle the solution decomposition problem \eqref{eq:sdproblem} (e.g., for anomaly detection).  In this work, we use the Flexible Generalized Golub-Kahan (FGGK) approach described in \cite{chung2022hybrid} to generate a basis for the solution, where the weights are determined from the group sparsity regularizer.

Given  $\bfA \in \mathbb{R}^{m \times n}$, $\bfb\in \mathbb{R}^{m}$, SPD matrices $\bfQ\in \mathbb{R}^{n\times n}$ and $\bfR\in \mathbb{R}^{m\times m}$, and iteration-dependent right preconditioning matrices  $\bfW_k^{-1} \in \mathbb{R}^{n \times n}$, the flexible, generalized Golub-Kahan iterative process generates vectors $\bfv_k, \bfz_k = \bfW_k^{-1}\bfv_k, $ and $\bfu_{k+1}$ such that at iteration $k$,
\begin{equation}\label{sdHybr_decomp}
\begin{bmatrix} \bfA \bfQ & \bfA\end{bmatrix} \widehat{\bfZ}_k = \bfU_{k+1} \bfM_k \quad \text{and} \quad \bfA\t \bfR^{-1} \bfU_{k+1}=\bfV_{k+1}\bfS_{k+1},
\end{equation}
where 
\begin{equation}
\widehat{\bfZ}_k = \begin{bmatrix} \bfv_1 &...& \bfv_k \\  \bfz_1 &...&   \bfz_k\end{bmatrix} = \begin{bmatrix} \bfV_k \\  \bfZ_k \end{bmatrix},
\end{equation}
and $\bfu_1 = \bfb / \|\bfb\|_{\bfR^{-1}}$ and $\bfv_1 = \bfA\t \bfR^{-1} \bfu_1$ if $\bfx_0=\bf0$. Note that $\bfM_k \in \mathbb{R}^{(k+1) \times k}$ is upper Hessenberg, $\bfS_{k+1} \in \mathbb{R}^{(k+1) \times (k+1)}$ is upper triangular, and in exact arithmetic, columns of $\bfU_{k+1} \in \mathbb{R}^{m \times (k+1)}$ and $\bfV_{k+1} \in \mathbb{R}^{n \times (k+1)}$ satisfy the orthogonality conditions,
\begin{equation}
    \bfU_{k+1}\t \bfR^{-1} \bfU_{k+1} = \bfI_{k+1} \quad {\bfV}\t_{k+1}\bfQ{\bfV}_{k+1} =  \bfI_{k+1}.
\end{equation}

In order to use the FGGK process to solve \eqref{eq:sdproblem}, we first consider the sequence of reweighted least-squares problems,
\begin{equation}
\label{eq:sdreweighted}
\min_{\bfpsi, \bfs} {\{\|\bfA \bfQ \bfpsi+\bfA \bfs-\bfb\|}^{2}_{\bfR^{-1}}+\alpha \|\bfpsi\|_{\bfQ}^2+\lambda{\|\bfW(\bfs)\bfs\|}_{2}^2\},
\end{equation}
where $\bfxi = \bfQ\bfpsi$ and $\bfW$ is defined as in \eqref{weights}.
We first project the objective function ${\|\bfA \bfQ \bfpsi+\bfA \bfs-\bfb\|}^{2}_{\bfR^{-1}}+\alpha \|\bfpsi\|_{\bfQ}^2$ onto a modified (flexible) Krylov subspace using the decomposition in \eqref{sdHybr_decomp}, then add regularization to the projected problem, such that an approximation of the solution can be found at each iteration $k$ as:
\begin{eqnarray}
        &\bfx_k& = \bfxi_k + \bfs_k = \bfQ \bfV_k \bfy_k + \bfZ_k \bfy_k \quad \text{where} \\ &\bfy_k& = \argmin_\bfy \{{\|\bfM_k \bfy-{\|\bfb\|}_{\bfR^{-1}} \bfe_1\|}_{2}^{2}+\alpha{\| \bfy \|}_{2}^{2} +\lambda{\| \bfR_{\bfW,k} \bfy \|}_{2}^{2}\}, \label{eq:flsqr_full_sd}
\end{eqnarray}
where $\bfZ_k = \bfQ_{\bfW,k} \bfR_{\bfW,k}$ is a thin QR factorization that can be computed via efficient updates \cite{chung2022hybrid}. Contrary to the additional regularization terms added to \eqref{eq:gmres_full} and \eqref{eq:flsqr_full}, a non-standard Tikhonov regularizer is included here to distinguish the contributions from the two solution components.
Notice that a decomposition of the solution $\bfx_k$ is available, i.e., $\bfxi_k$ is a reconstruction of the smooth component and $\bfs_k$ is a reconstruction of the group-sparse component (e.g., containing anomalies). Moreover, in this case, the weights depend only on a part of the approximate solution at each iteration $k$, i.e. ${\bfW}_{k} = \widetilde{\bfW}(\bfs_{k-1})=\widetilde{\bfW}(\bfZ_{k-1} \bfy_{k-1})$. A detailed explanation of the algorithm developed to obtain \eqref{sdHybr_decomp}, originally proposed under the name of solution decomposition hybrid projection approach (sdHybr), can be found in \cite{chung2022hybrid}. 

\subsubsection{Regularization parameter choice criteria}
One of the main advantages of flexible Krylov methods is the ability to incorporate regularization on the projected problem (IRW-F and hybrid-F), and hence to adaptively set the regularization parameter(s) $\bflambda_k$ at each iteration. In particular, we consider the use of the discrepancy principle (DP)\cite{Morozov1966OnTS} for the projected problem, i.e. finding the regularization parameters $\bflambda \geq 0$ such that 
\begin{equation}\label{dp}
    \| \bfA \bfx_k(\bflambda)- \bfb\|_2 = \| \bfA \bfZ_k \bfy_k(\bflambda)- \bfb\|_2= \eta \, \|\bfe\|_2.
\end{equation} 
Note that $\bflambda=\lambda$ for the minimization in \eqref{eq:gmres_full} and \eqref{eq:flsqr_full}, but $\bflambda =[\lambda,\alpha]$ for the optimization in \eqref{eq:flsqr_full_sd}.  The approximate solution computed by any given algorithm at iteration $k$ is given by $\bfx_k(\bflambda)$, with $\bfy_k$ being the corresponding coefficients in the projected space. Using the DP requires a good estimate of the noise level. Alternatively other methods can be used seamlessly, see, e.g. \cite{chung2021computational}\cite{msl2}.

\section{Numerical examples} \label{sec:numerics}

In this section three imaging examples are presented to demonstrate the performance of flexible Krylov methods for group sparsity regularization. In order to showcase the generality of the described approach, the flexible schemes are used for different grouping strategies including both overlapping and non-overlapping groups and in combination with different solvers.  

The notation used to describe the solvers in this section can be found in Table \ref{tab:table1}, where group sparsity regularization is indicated by appending a `G' suffix at the end of the method. In the case where more than one grouping strategy has been compared, different grouping strategies are denoted `G1', `G2', etc. Moreover, if both $\ell_1$ and $\ell_{2,1}$ are considered simultaneously, we use `C' to denote the combined approach. If no group suffix is appended for a flexible method, then this corresponds to $\ell_1$ regularization. Recall that the projection types are explained in Section \ref{sub:flexible}. And note that hybrid solution decompose (SD) \cite{chung2022hybrid} and restarted conjugate gradient (RCG) for group sparsity \cite{Chen2014} are used as a comparison with previous works.

\begin{table}[h!]
    \captionsetup{width=\textwidth}
    \caption{Notation for the solvers: \textbf{projection} - \textbf{ALGORITHM} - \textbf{GROUPS}}
    \label{tab:table1}
    \begin{tabular}{c|c|c} 
      \textbf{projection}  & \textbf{ALGORITHM} & \textbf{GROUPS}\\
      none   & FLSQR & none \\
      hybrid & FGMRES & G \\
      IRW & SD & G1, G2 \\
            & RCG & C \\
    \end{tabular} 
\end{table}

The first experiment exploits group sparsity patterns in wavelet coefficients by considering the natural tree structure of wavelet decompositions. Since non-leaf elements belong to more than one branch, this example showcases the use of the new scheme for regularization with overlapping groups. A deblurring problem is considered and solved using different versions of the new scheme involving FLSQR, paying particular attention to the different projection types and recalling the differences between them while showcasing their practical performance. 

The second experiment corresponds to a dynamic deblurring problem. Since this is a spatio-temporal problem displaying a spatial sparsity pattern but smoothness in time, the groups are chosen to contain each pixel across all time points, so this is an example of non-overlapping group sparsity. In particular, we focus on the difference between enforcing $\ell_1$ regularization and $\ell_{2,1}$ regularization using hybrid flexible methods, and we present results for the new method combining both regularization terms. Since this is a symmetric problem, the performance of both methods based on FLSQR and FGMRES is tested. 

The third experiment concerns an anomaly detection problem in atmospheric inverse modeling. This corresponds to a realistic example, and our aim is to show the potential of the new schemes in real-world applications. Since this experiment also involves spatio-temporal images that are sparse in space but not time, the groups are built analogously to the second experiment. Additionally, since the solution in this experiment is modelled as the sum of a smooth (background) component and a group-sparse (anomalies) component, hybrid-SD-G is used in this experiment.

For all the experiments, $\tau$ in \eqref{eq:smooth_W} is taken to be $10^{-10}$. However, empirical observations demonstrate that the algorithms are robust to this choice as long as $\tau$ is significantly smaller than the average pixel intensity and kept above machine precision. Moreover, the regularization parameter for all hybrid and iteratively-reweighted methods is chosen across the iterations using the DP defined in \eqref{dp}.

\subsection{Deblurring example with wavelet sparsity patterns}
\begin{figure}[tb]
     \centering
     \begin{subfigure}[b]{0.47\textwidth}
         \centering
         \includegraphics[trim = 4cm 3cm 4cm 2cm,clip,height=4cm]{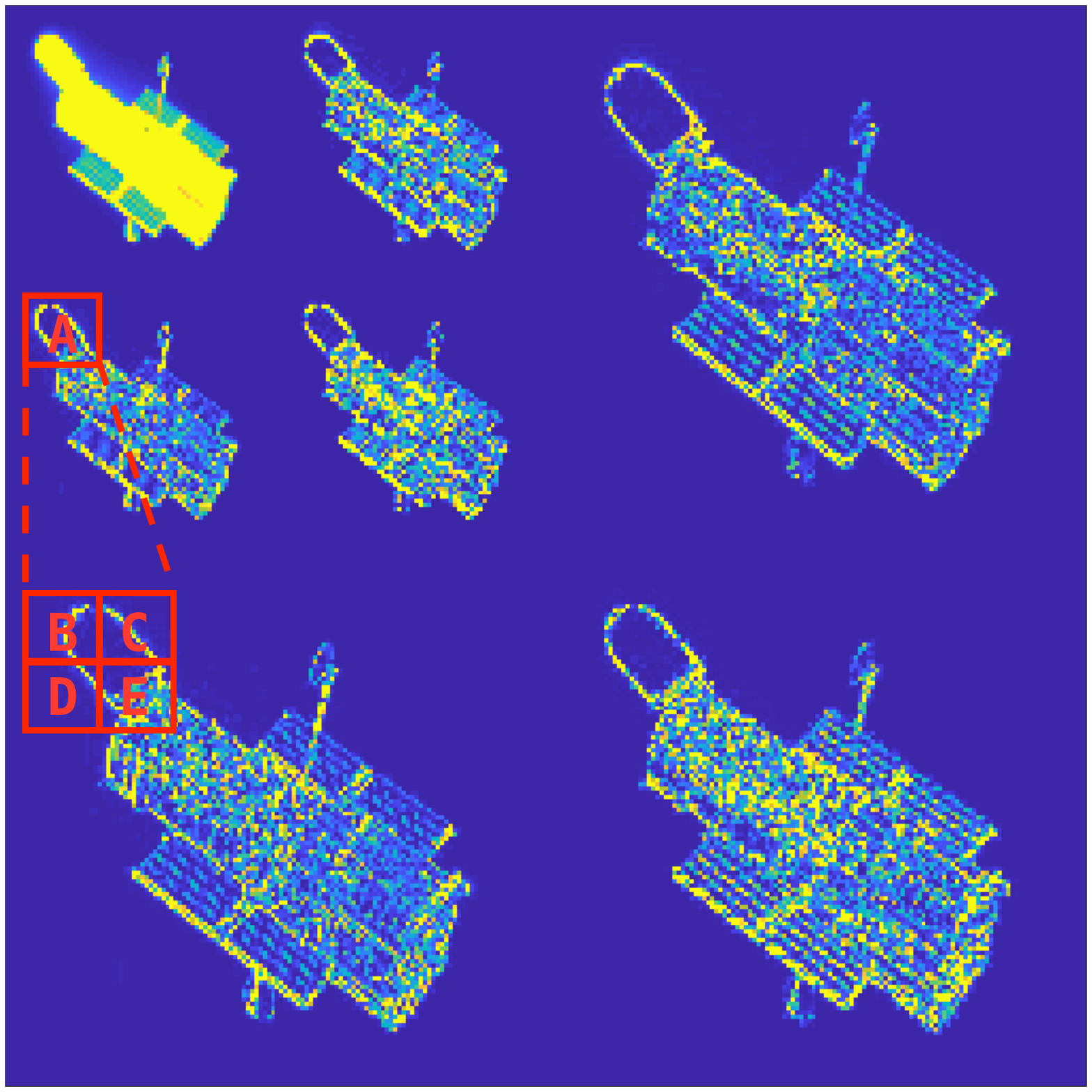}
         \caption{2 level Haar Wavelet representation}
         \label{tree1}
     \end{subfigure}
     \begin{subfigure}[b]{0.47\textwidth}
         \centering
         \includegraphics[trim = 2cm 3cm 4cm 1cm,clip,height=4cm]{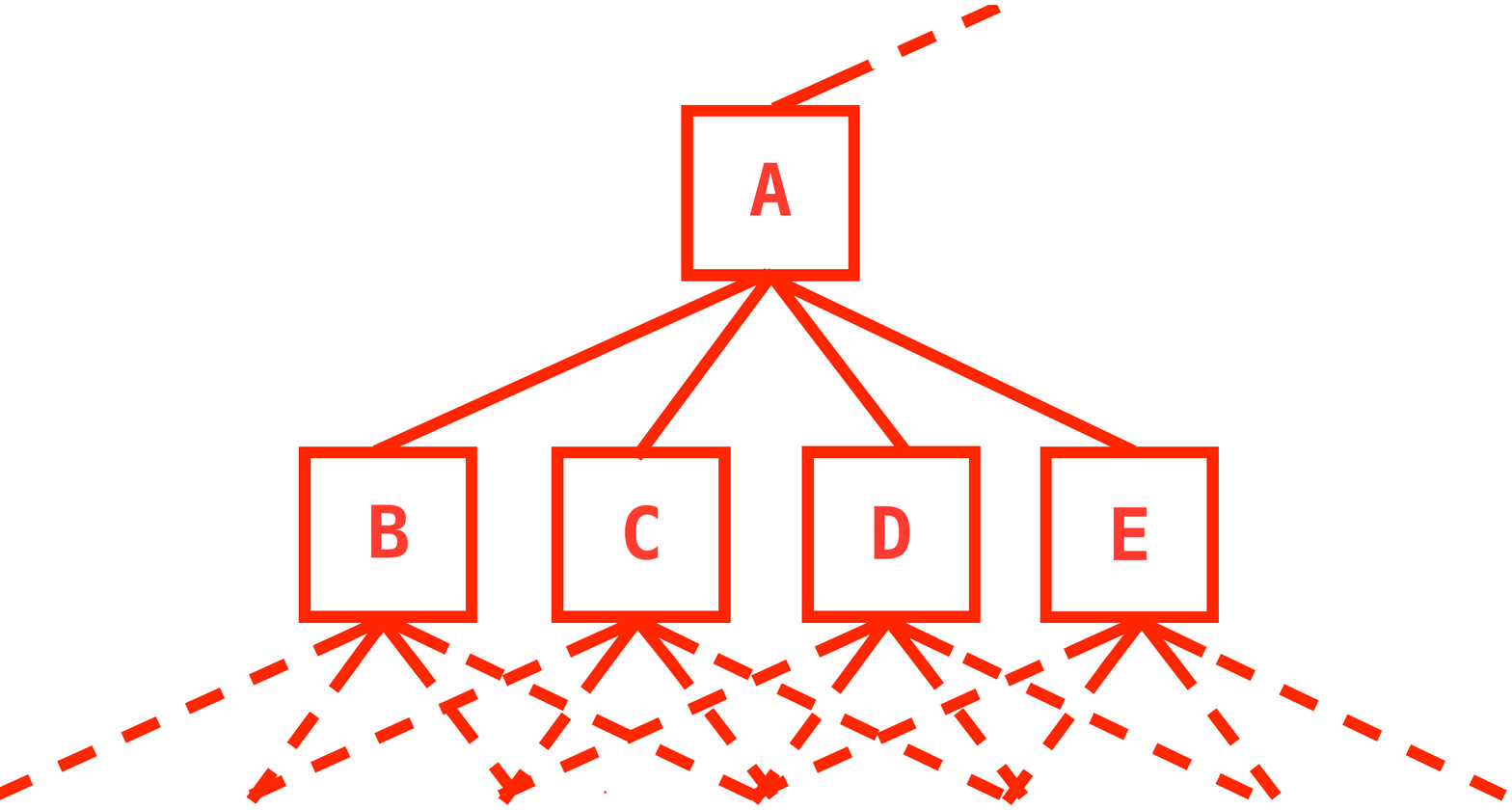}
         \caption{Tree structure for the example in (a)}
         \label{fig:three sin x}
     \end{subfigure}
        \caption{Illustrative example of the tree structure in the wavelet decomposition of the \textsc{satellite} test image used for the deblurring problem represented in Figure \ref{fig:wavelet_examples} and Haar wavelets with 2 levels. Examples of groups using strategy G1 are $\{A,B\}$, $\{A,C\}$, $\{A,D\}$, $\{A,E\}$ while an example of a group using strategy G2 is $\{A,B,C,D, E\}$.}
        \label{fig:trees_example}
\end{figure}

Natural images have been shown to be well described by wavelet representations \cite[Chapter 9]{doi:10.1137/1.9781611970104}. Moreover, wavelets have a natural tree-structure where each non-leaf wavelet coefficient for a given orientation (diagonal, vertical and horizontal) has four children coefficients corresponding to the same orientation at a finer scale. This is explained in, e.g., \cite{6115845}. There are of course different groupings than one could develop. For example, in this paper, both 2-elements groups of each children with their parent (G1) or 5-elements groups of each parent with all their children (G2) are considered as suggested in \cite{6115845}. In Figure \ref{fig:trees_example} an upscaled example for a wavelet decomposition with two levels is shown to illustrate a case of parent-child relationship in the wavelet domain. Note that Figure \ref{fig:trees_example} is intended as an illustration, but in reality the tree structure is considered at the pixel scale, i.e. where the red boxes in Figure \ref{fig:trees_example} correspond to individual pixels, to match the wavelet coefficients. Thus, at each level, each given pixel in the wavelet representation corresponding to a non-leaf element in the tree is the parent of 4 pixels in a finer level representing the same area of the original image (at that given orientation).

\begin{figure}[tb]
     \centering
     \begin{subfigure}[b]{0.4\textwidth}
         \centering
         \includegraphics[trim = 3.2cm 1cm 3.2cm 1cm,clip,height =4cm]{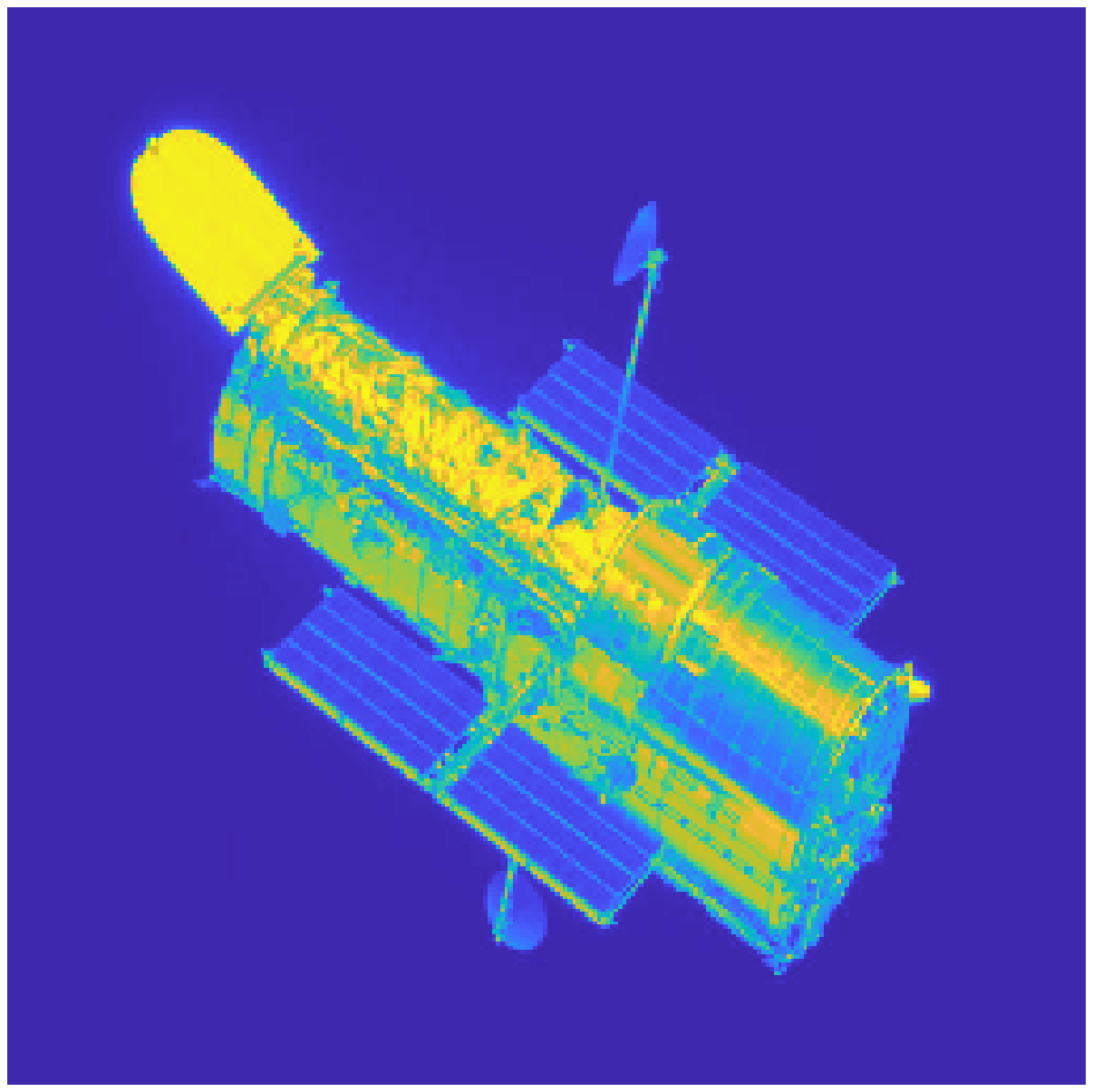}
         \caption{\textsc{satellite} $\bfx_{\rm true}$}
         \label{deblur_x}
     \end{subfigure}
     \begin{subfigure}[b]{0.4\textwidth}
         \centering
         \includegraphics[trim = 3.2cm 1cm 3.2cm 1cm,clip,height =4cm]{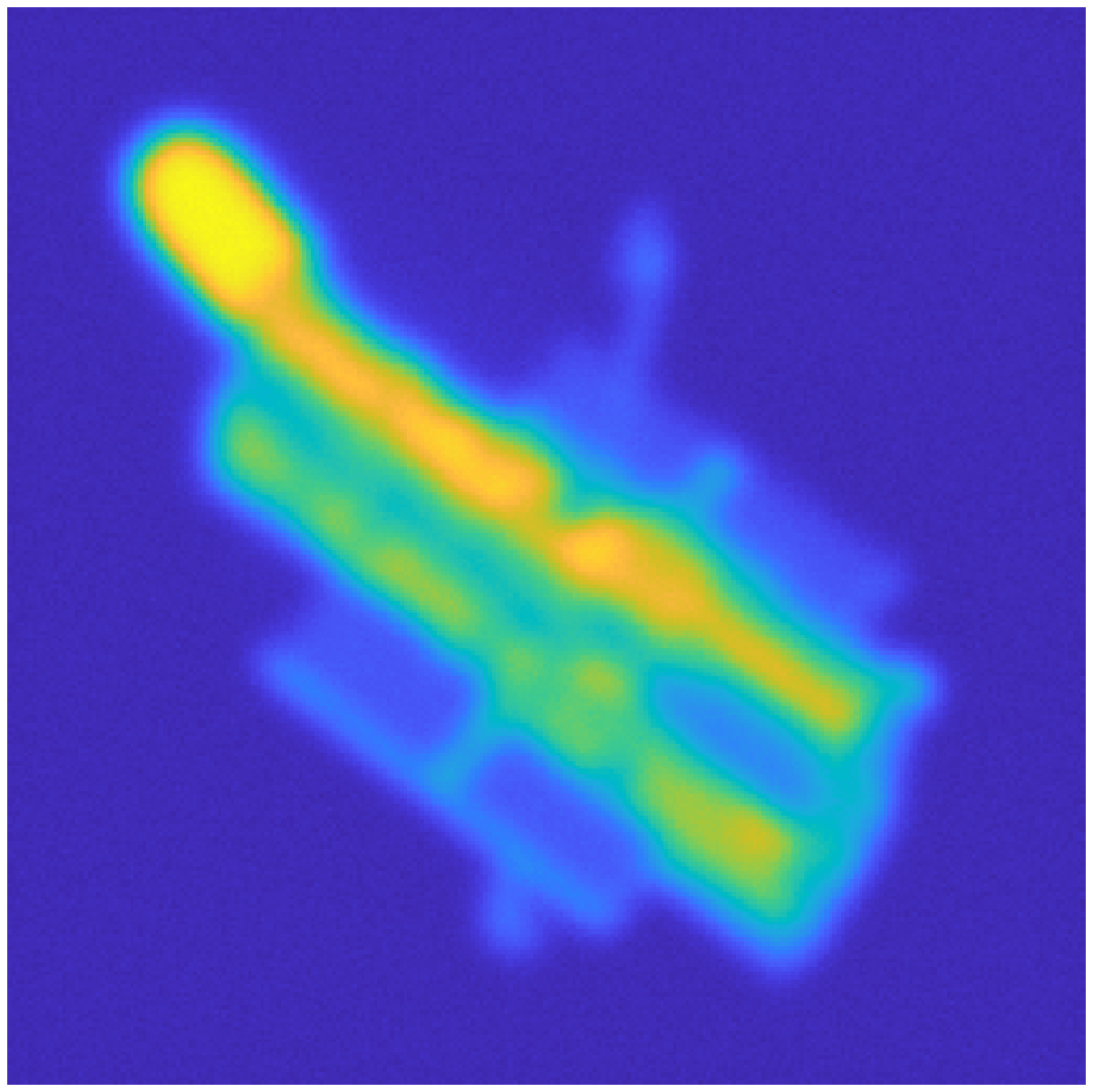}
         \caption{observed, blurred $\bfb$}
         \label{fig:deblur_b}
     \end{subfigure}
        \caption{Image deblurring problem for use with group sparsity in the wavelet coefficients.}
        \label{fig:wavelet_examples}
\end{figure}

For this deblurring example, medium Gaussian blur is generated with IR Tools \cite{gazzola2018ir} with added Gaussian noise at noise level $\|\bfe\| / \|\bfA \bfx_{\rm true}\| =0.05$. Here the true solution $\bfx_{\rm true} \in \mathbb{R}^{65,536}$ contains $256\times 256$ pixels and corresponds to the \textsc{satellite} test image in IR Tools and can be observed in Figure \ref{fig:wavelet_examples} along with the corresponding blurred and noisy measurements. One can observe that the true image contains regions with details (and edges) at different scales, while other areas of the image contain many zeroes and hence do not have information at any scale. Therefore, group sparsity in the wavelet coefficients should provide suitable regularization. 

\begin{figure}[tb]
     \centering
     \begin{subfigure}[b]{0.48\textwidth}
         \centering
         \includegraphics[height=4.7cm]{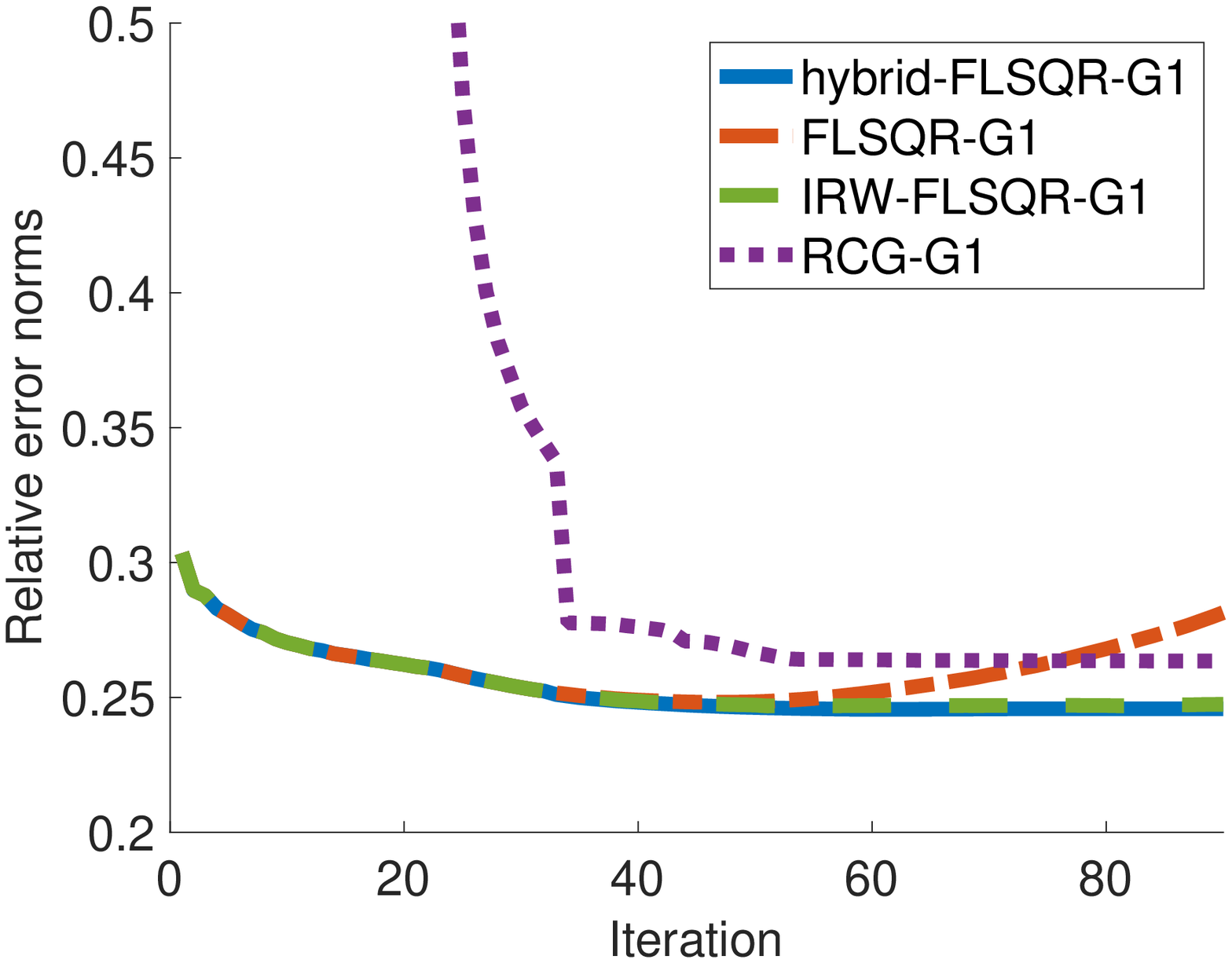}
         \caption{Error norms for G1}
         \label{fig:errors_deblur_1}
     \end{subfigure} \hspace{0.15cm}
     \begin{subfigure}[b]{0.48\textwidth}
         \centering
         \includegraphics[height =4.7cm]{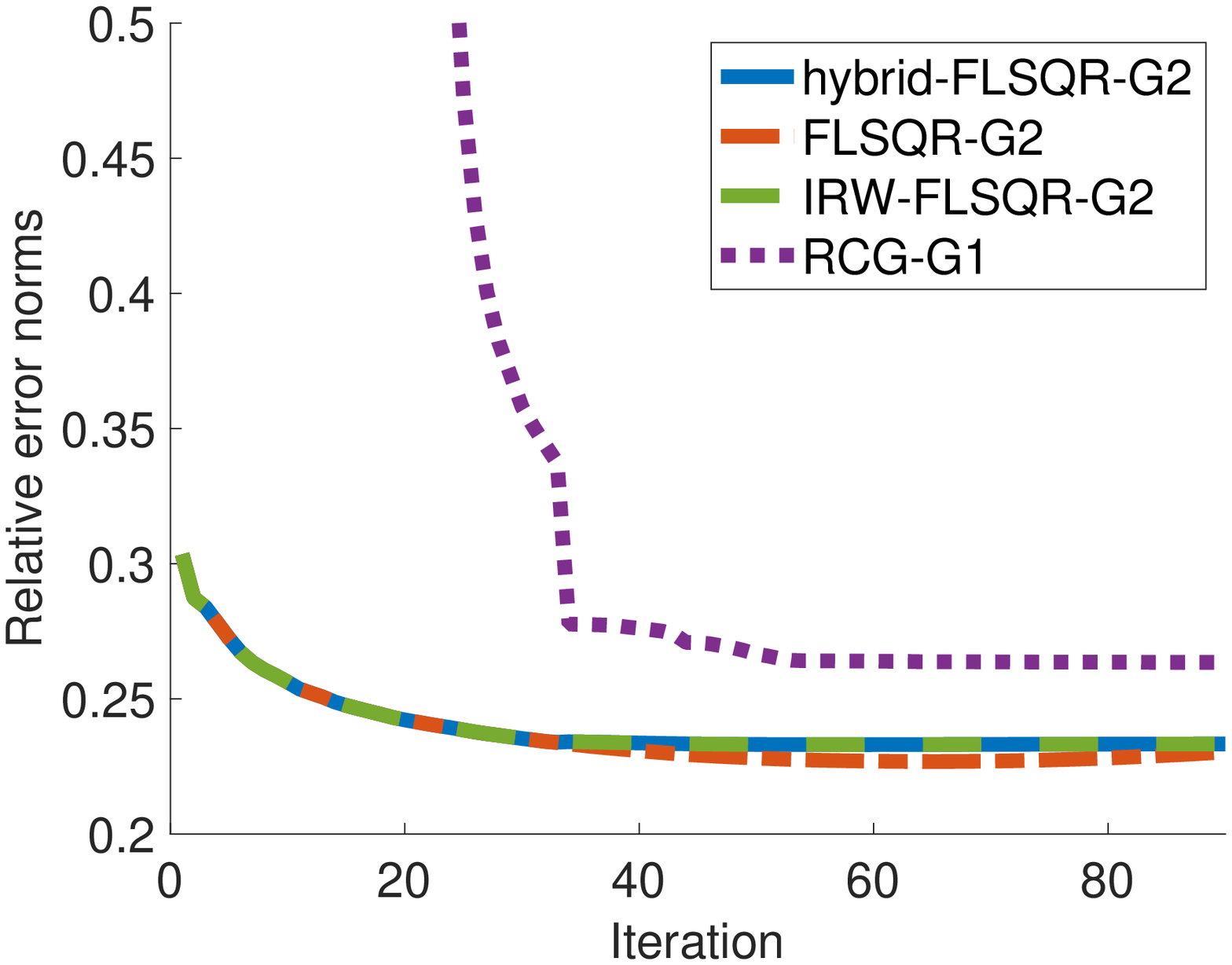}
         \caption{Error norms for G2}
         \label{fig:errors_deblur_2}
     \end{subfigure}\hspace{0.15cm}
        \caption{Relative reconstruction error norms for the methods based on FLSQR for the deblurring example presented in Figure \ref{fig:wavelet_examples} with wavelet group sparsity. The regularization parameter has been chosen at each iteration using the DP.}
        \label{wavelet_results}
\end{figure}

Exploiting group sparsity in the wavelet coefficients corresponds to choosing $\bfPsi$ in \eqref{eqn:problem_gs_1} to be an orthogonal wavelet transform and choosing the groups according to G1 or G2 as defined at the beginning of this subsection (and illustrated in Figure \ref{fig:trees_example}). Note that both grouping schemes correspond to overlapping group sparsity regularization. The results for these experiments can be found in Figure \ref{wavelet_results}, where we provide relative reconstruction error norms, computed as $\norm[2]{\bfx_k - \bfx_{\rm true}}/\norm[2]{\bfx_{\rm true}}$ where $\bfx_k$ is the reconstruction at the $k$th iteration. 

We provide a comparison against the algorithm presented in \cite{Chen2014} (following their original grouping structure G1). Note that we have not included their proposed preconditioning, nor their starting guess ($\bfx_0=\bf0$ is used instead of $\bfA\t \bfb$), since those are highly problem-dependent choices and their approach is not suitable for this example. The regularization parameter for the hybrid flexible and IRW flexible methods is chosen at each iteration using the DP \eqref{dp} with safety parameter $\eta=1.01$. For the RCG algorithm \cite{Chen2014}, the regularization parameter is chosen to be the one computed with the hybrid flexible methods at the end of the iterations. It is worth mentioning that RCG is better suited to other forward models (such as MRI), but flexible Krylov methods have a comparable performance and allow the (semi-automatic) computation of the regularization parameter on-the-fly.

\begin{figure}[tb]
         \centering
         \includegraphics[height =4.7cm]{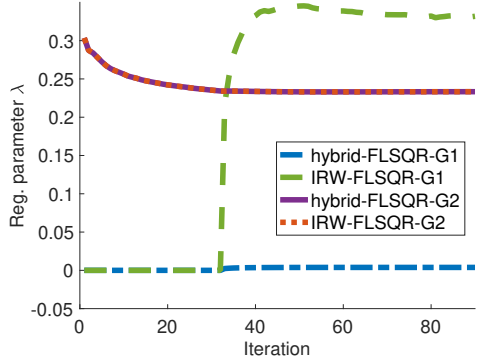}
        \caption{Regularization parameter $\lambda_k$ for the deblurring example presented in Figure \ref{fig:wavelet_examples}, chosen at each iteration using the DP.}
        \label{fig:regparam_deblur}
\end{figure}

One can observed that flexible Krylov methods display a fast convergence in comparison to RCG. Different grouping strategies, displayed in Figure \ref{fig:errors_deblur_1} and Figure \ref{fig:errors_deblur_2} respectively, have similar performances (note that both error norm plots have the same scale for the axis). Moreover, one can observe that flexible Krylov methods without regularization in the projected problem (FLSQR-G1 and FLSQR-G2) display semi-convergence, and that the severity of it depends on the grouping strategy. To be able to use this algorithms in practice, one should use appropriate early stopping. On the contrary, methods with explicit regularization in the projected problem (hybrid-FLSQR-G and IRW-FLSQR-G for both grouping strategies), using the DP as a parameter choice criterion at each iteration, display a stabilization of the error norm throughout the iterations. Is is worth noting that the regularization parameter for hybrid and IRW methods does not necessarily need to be the same as that for the original problem, as explained in Section \ref{sub:flexible}.  In fact, for the grouping strategy G1, very similar error norms are achieved for the two methods using very different regularization parameters. This can be observed in Figure \ref{fig:regparam_deblur}. Recall that hybrid methods are cheaper per iteration than iteratively-reweighted methods, but iteratively-reweighted methods are equipped with theoretical guarantees of convergence. In practice, as can be observed for this experiment, both methods typically display very similar behaviours in the error norm plots.

\subsection{Combined regularization for spatio-temporal image deblurring}
\label{sec:combined} 

\begin{figure}[htb]
\begin{center}
\begin{tabular}{c c c c c c c c c}
\raisebox{1.5\normalbaselineskip}[0pt][0pt]{\rotatebox[origin=c]{90}{}}  
\raisebox{1.5\normalbaselineskip}[0pt][0pt]{\rotatebox[origin=c]{90}{$\bfx_{\rm true}$}}  
\raisebox{1.5\normalbaselineskip}[0pt][0pt]{\rotatebox[origin=c]{90}{}}   \hspace{-0.4cm} &
\includegraphics[width=1.5cm]{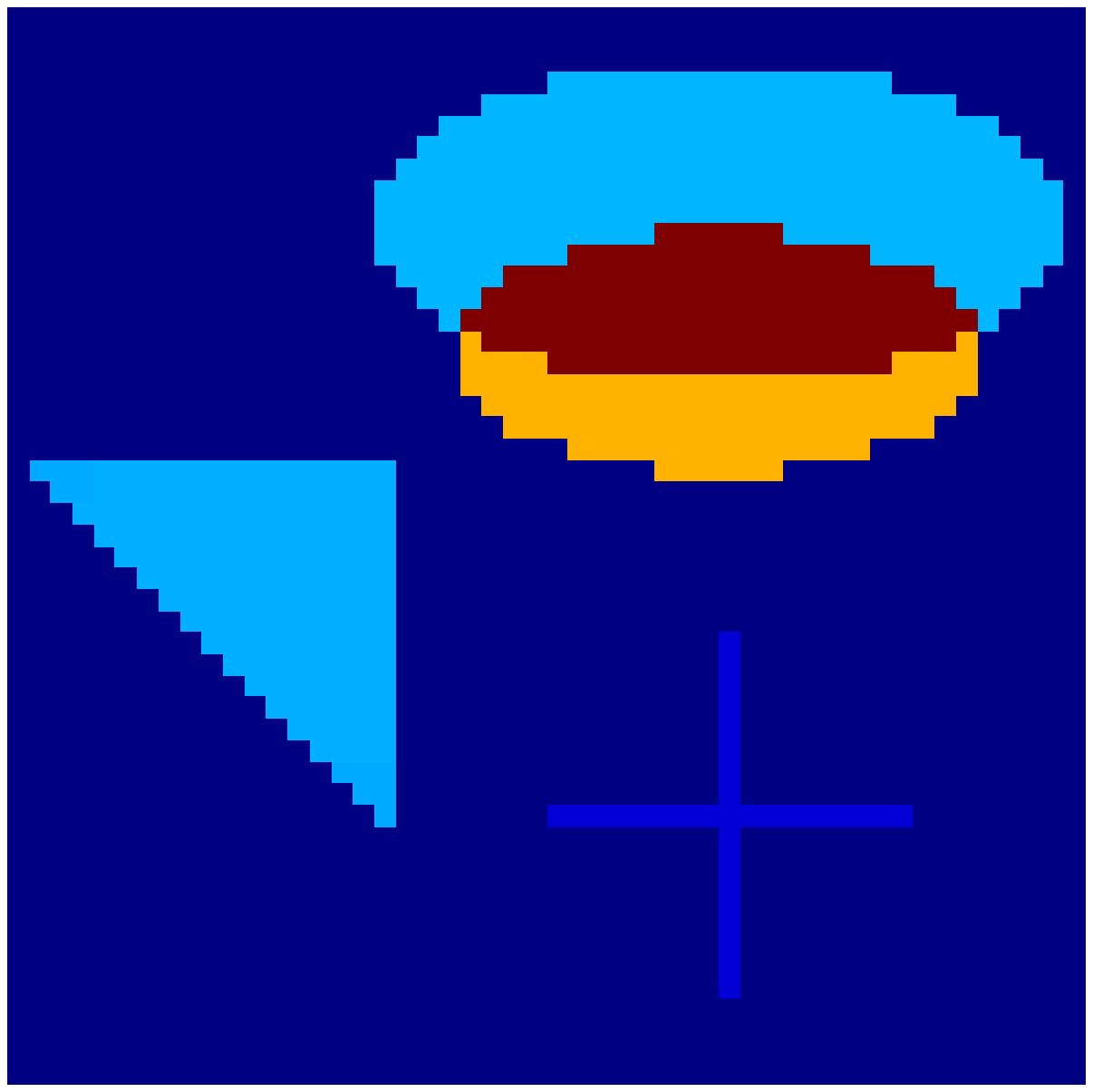} \hspace{-0.5cm} & \includegraphics[width=1.5cm]{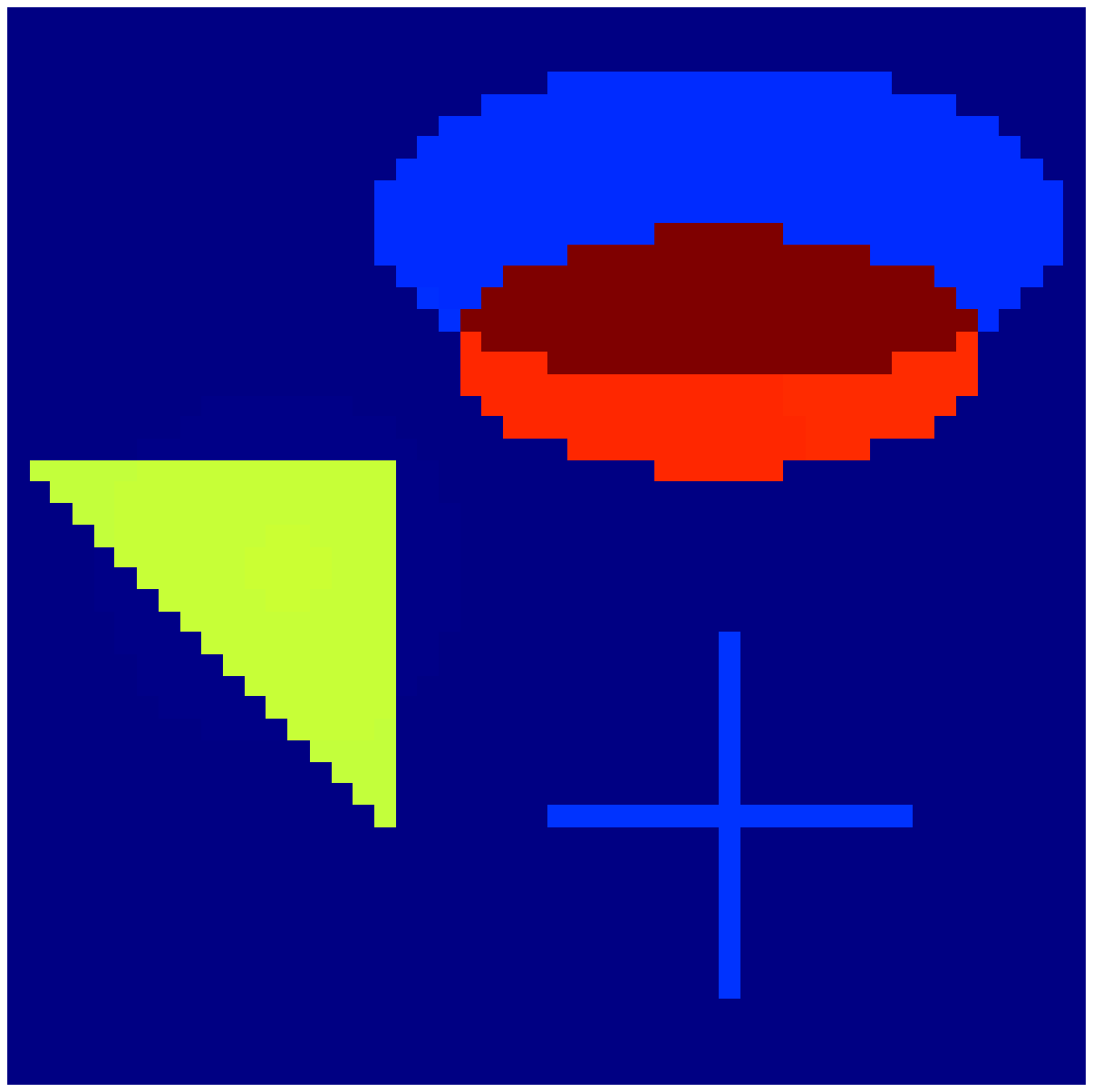} \hspace{-0.5cm} & \includegraphics[width=1.5cm]{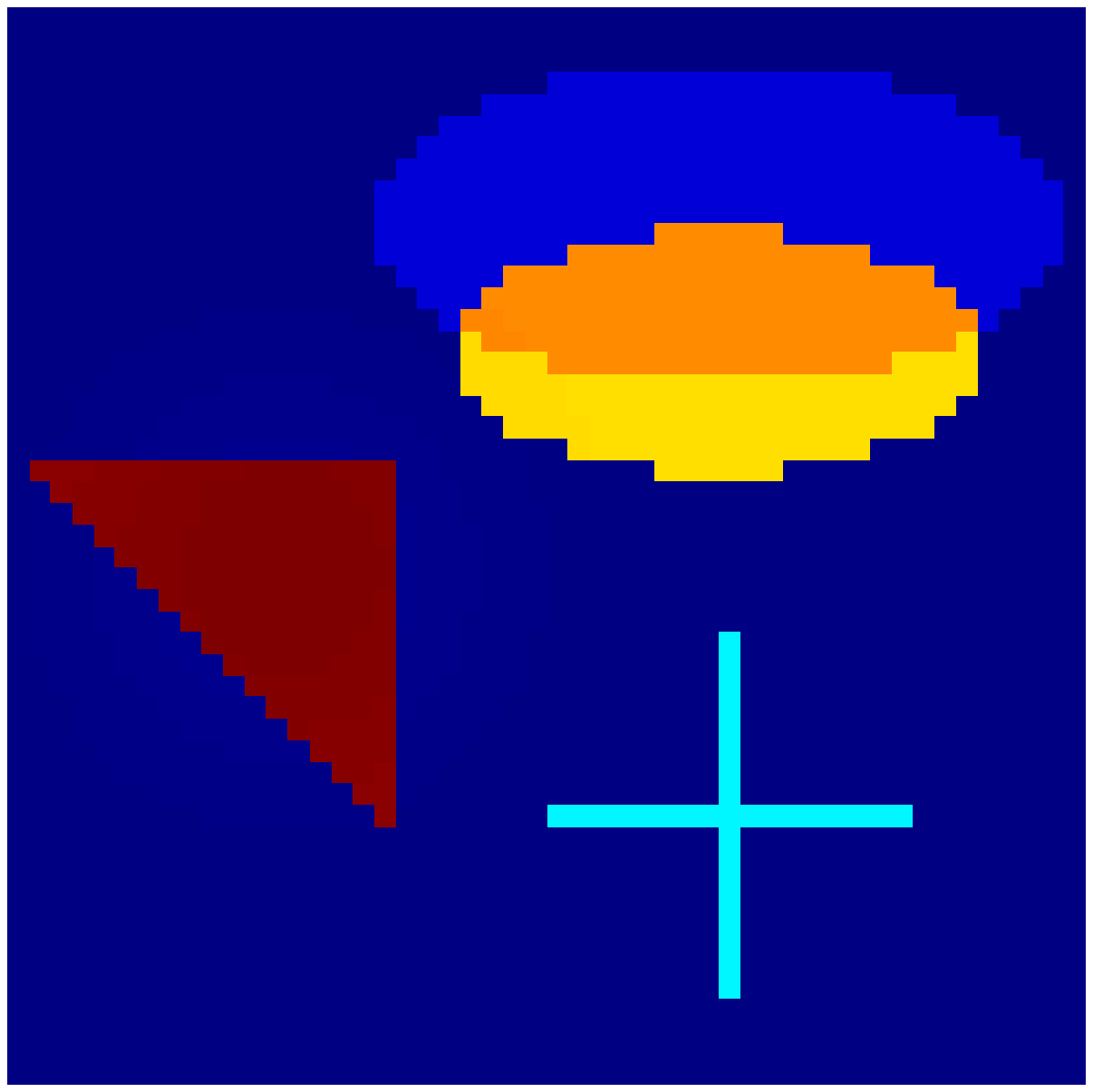} \hspace{-0.5cm} & \includegraphics[width=1.5cm]{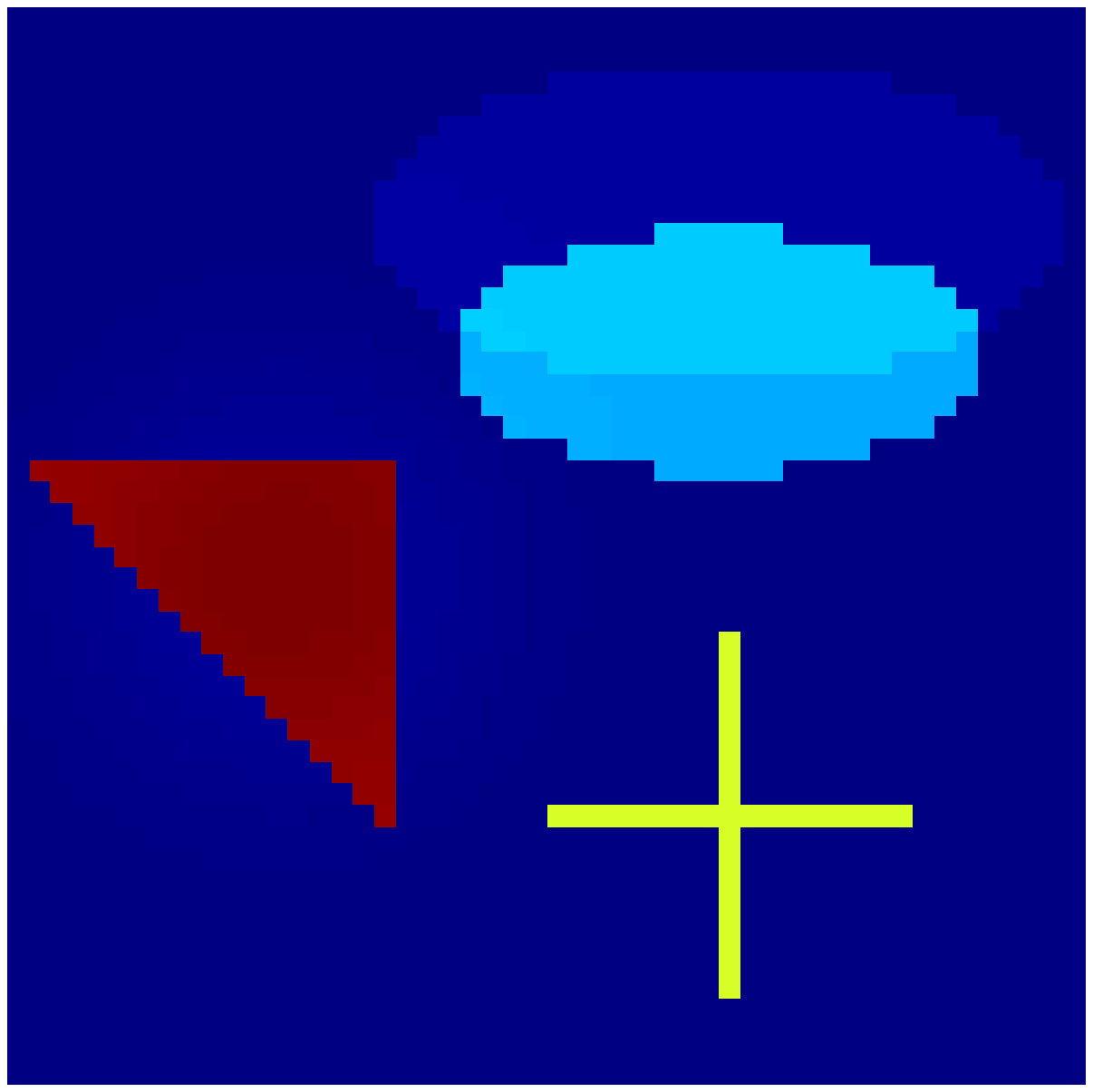} \hspace{-0.5cm} &  \includegraphics[width=1.5cm]{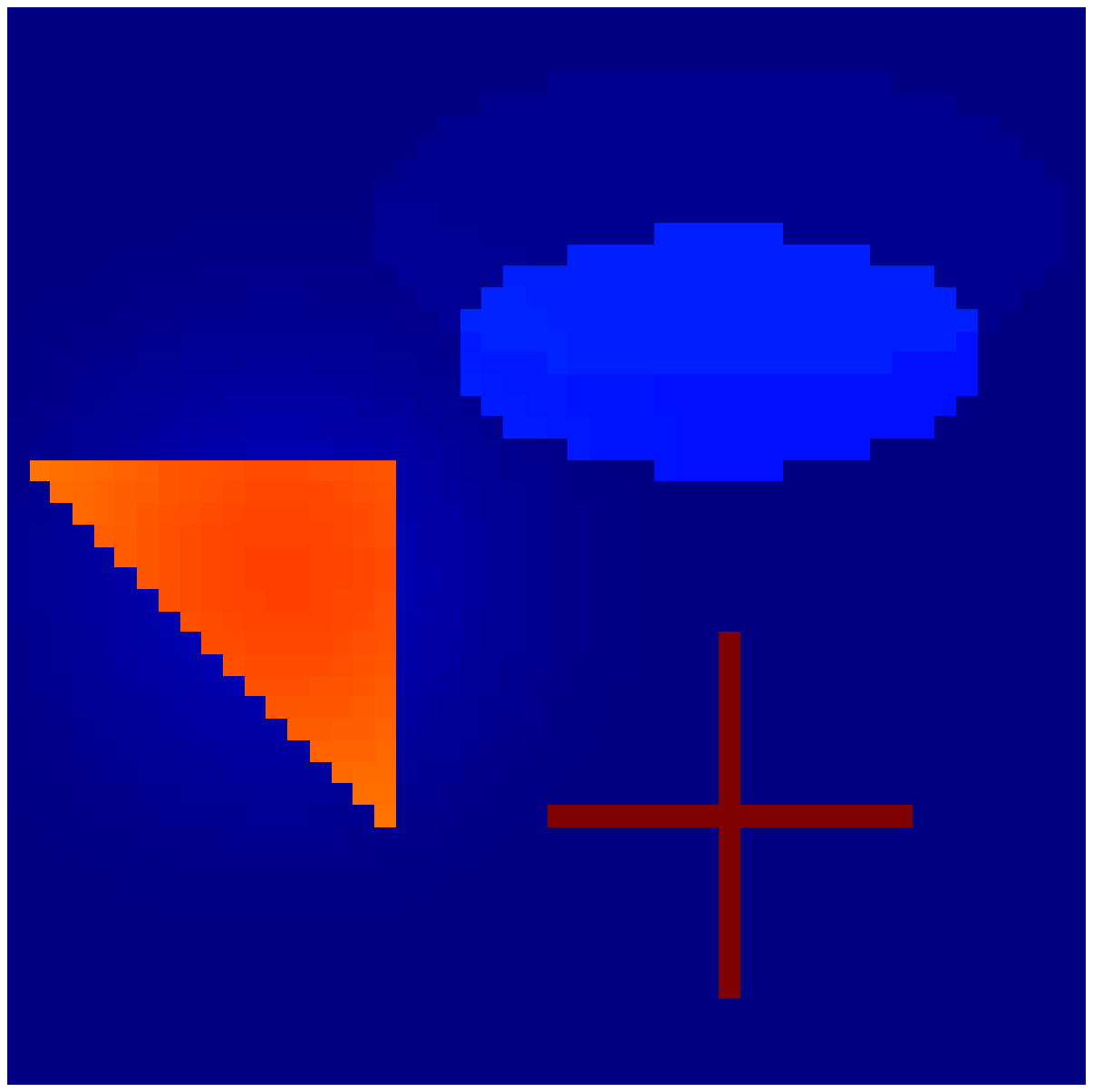} \hspace{-0.5cm} & 
\includegraphics[width=1.5cm]{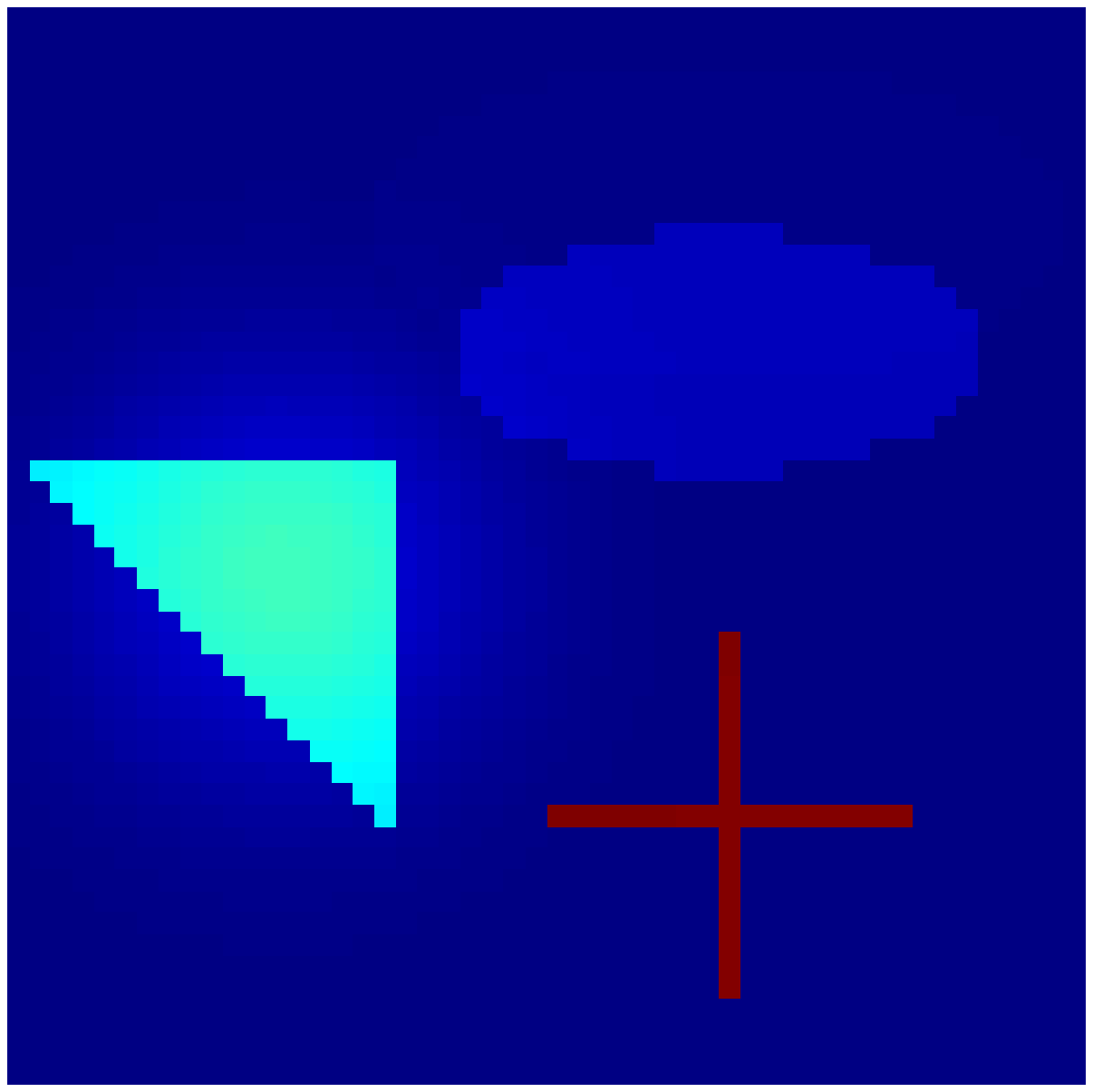} \hspace{-0.5cm} & 
\includegraphics[width=1.5cm]{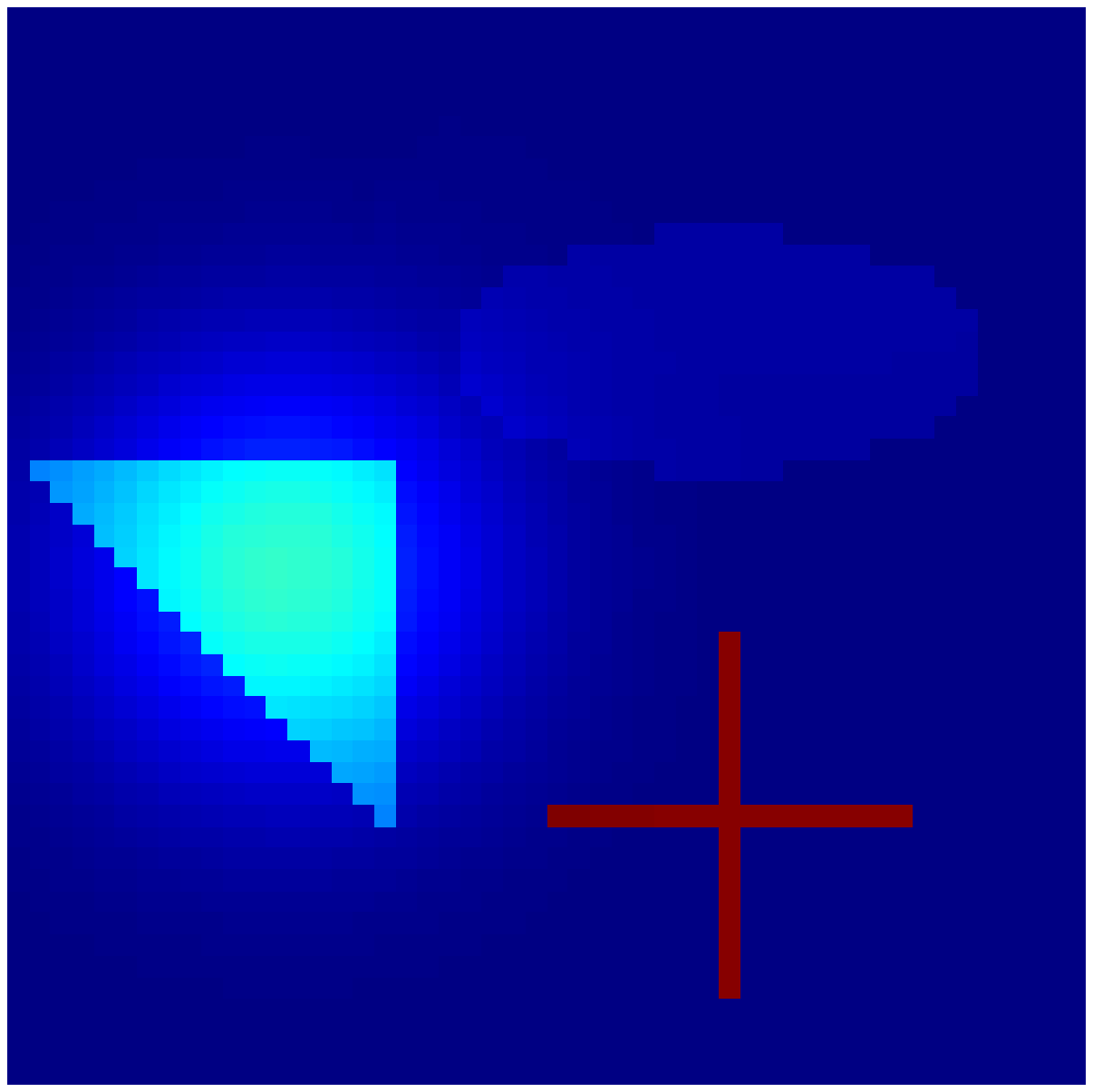} \hspace{-0.5cm} & \includegraphics[width=1.5cm]{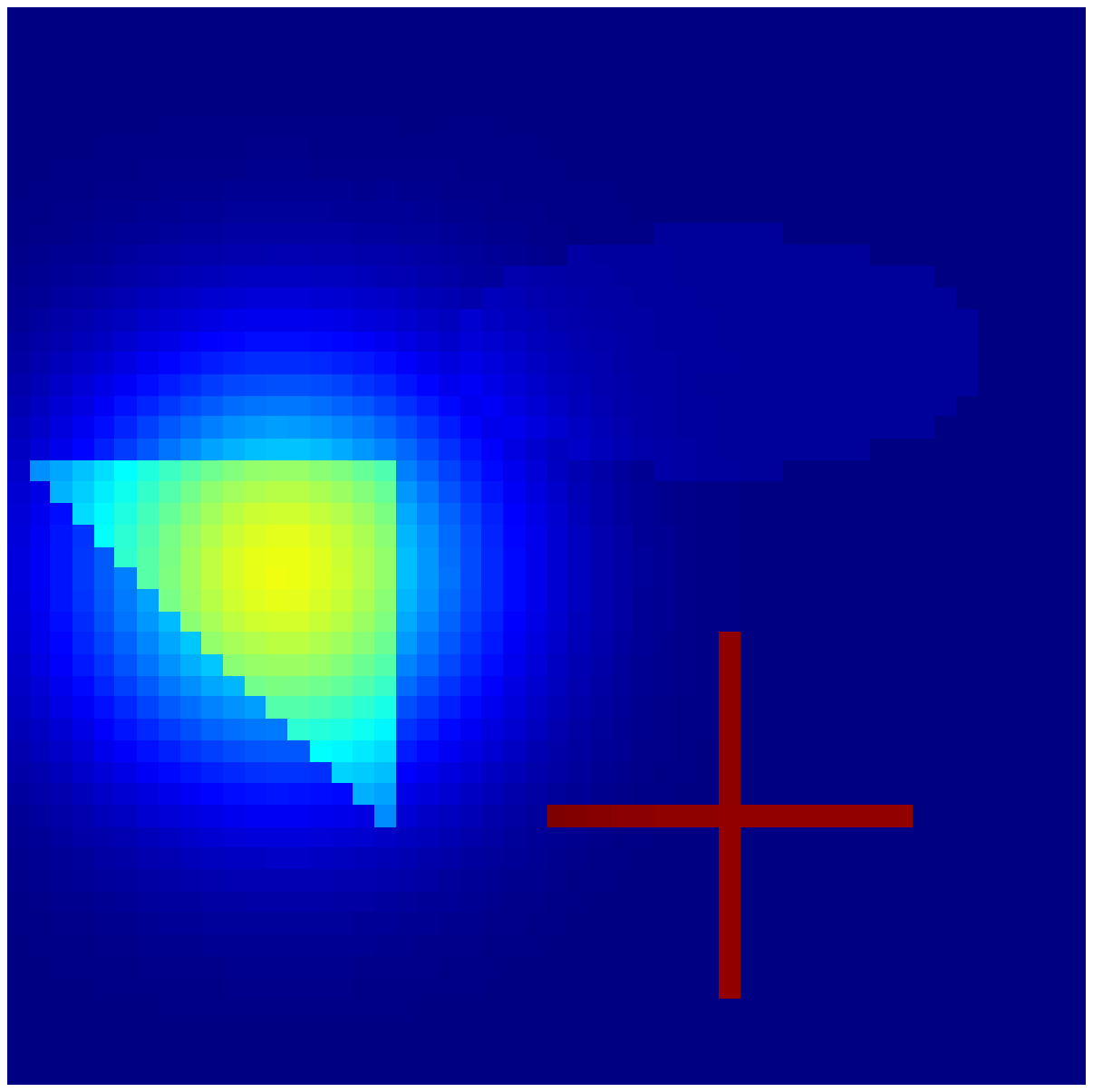}  \\ 
\raisebox{1.5\normalbaselineskip}[0pt][0pt]{\rotatebox[origin=c]{90}{}}  
\raisebox{1.5\normalbaselineskip}[0pt][0pt]{\rotatebox[origin=c]{90}{$\bfb$}}  
\raisebox{1.5\normalbaselineskip}[0pt][0pt]{\rotatebox[origin=c]{90}{}}   \hspace{-0.4cm} &
\includegraphics[width=1.5cm]{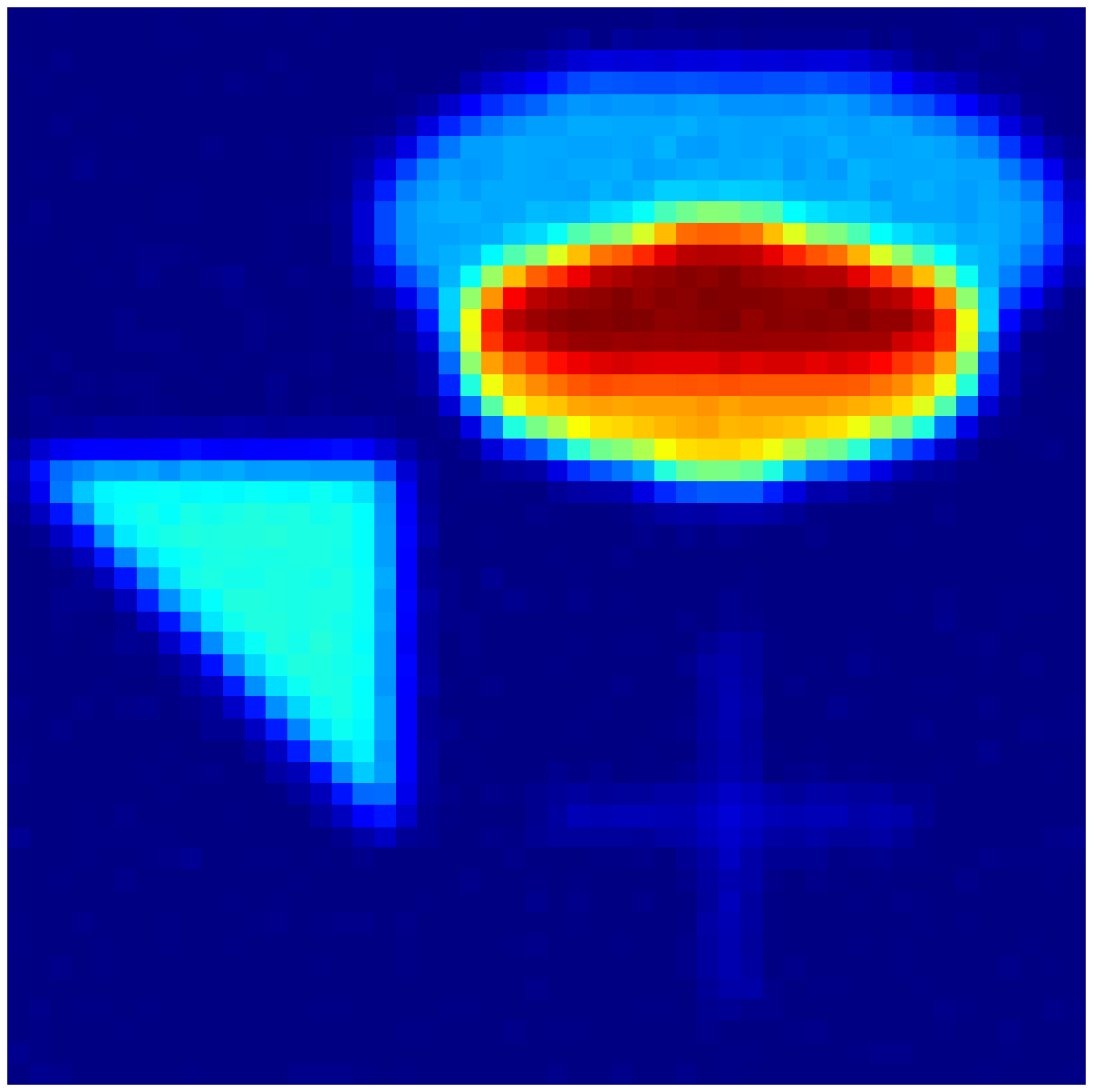} \hspace{-0.5cm} & \includegraphics[width=1.5cm]{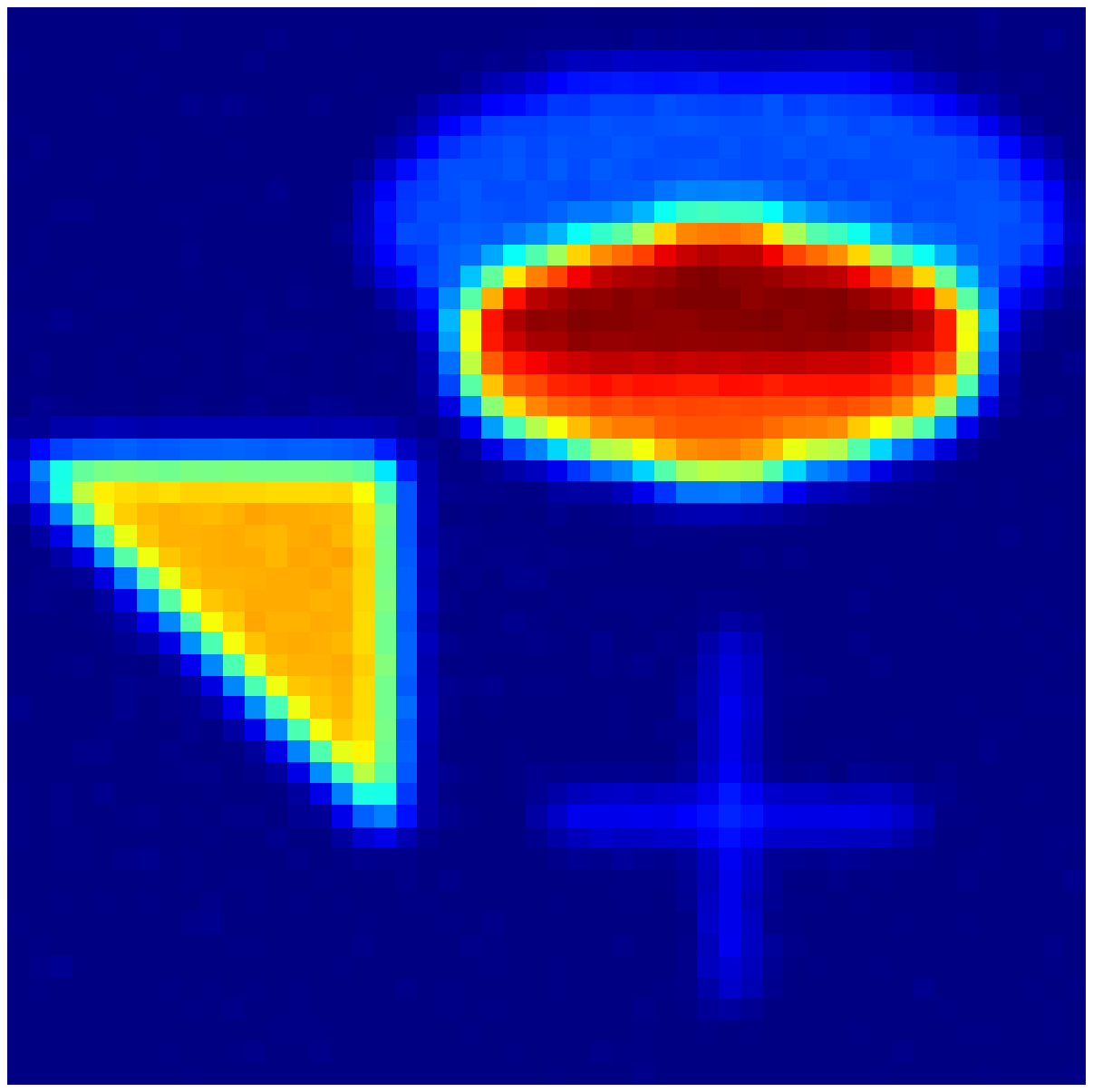} \hspace{-0.5cm} & \includegraphics[width=1.5cm]{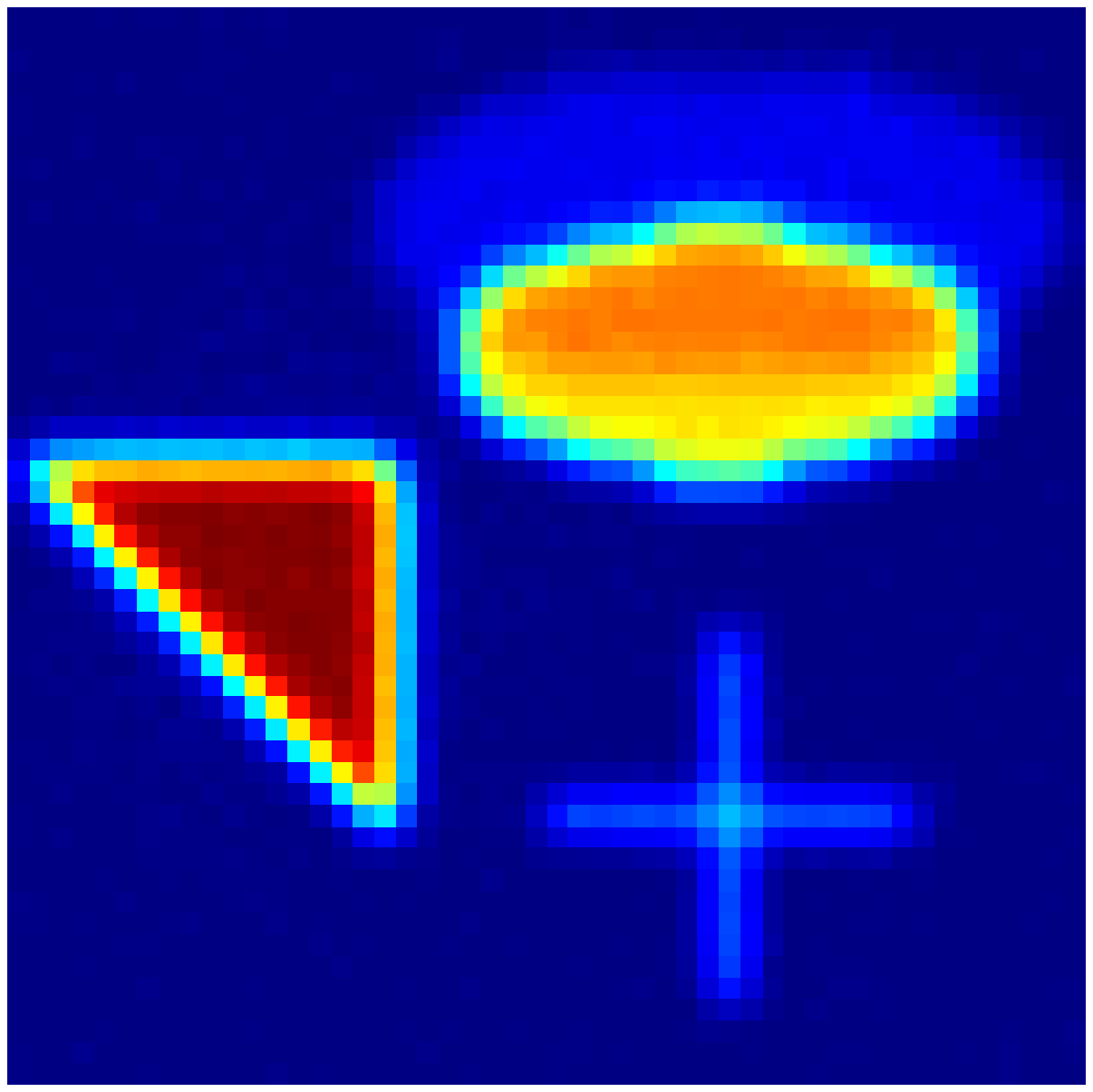} \hspace{-0.5cm} & \includegraphics[width=1.5cm]{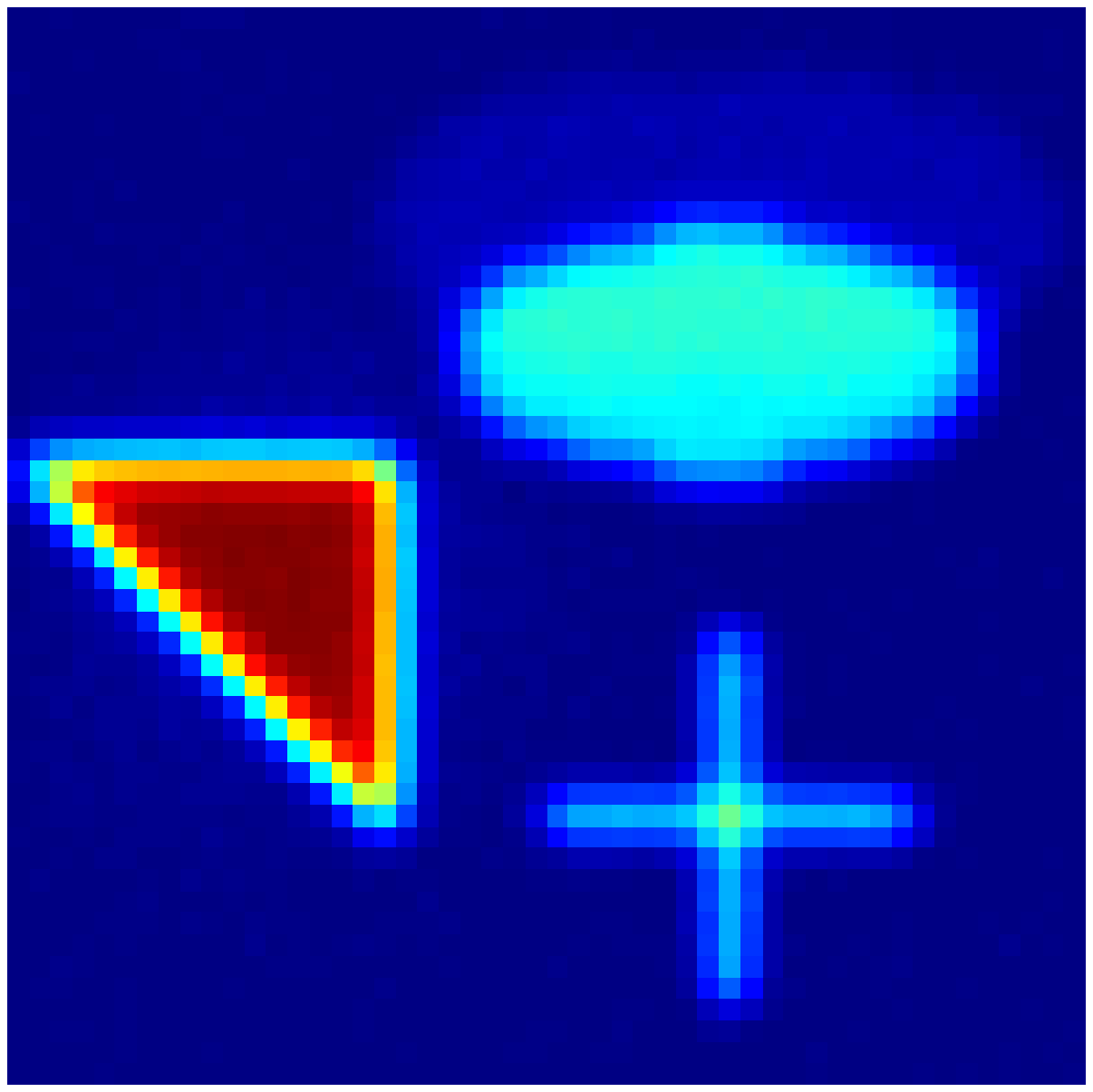} \hspace{-0.5cm} &  \includegraphics[width=1.5cm]{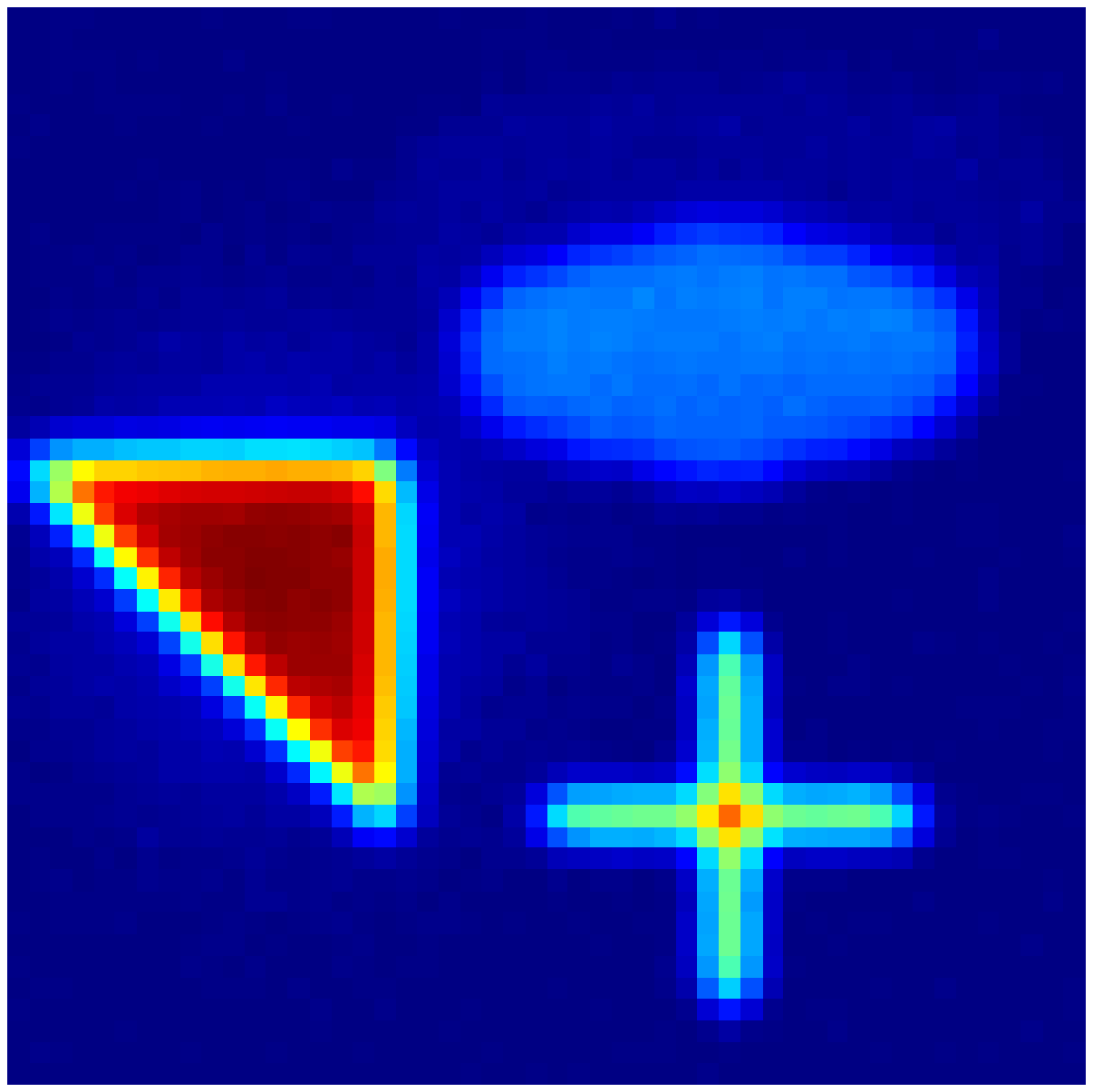} \hspace{-0.5cm} & 
\includegraphics[width=1.5cm]{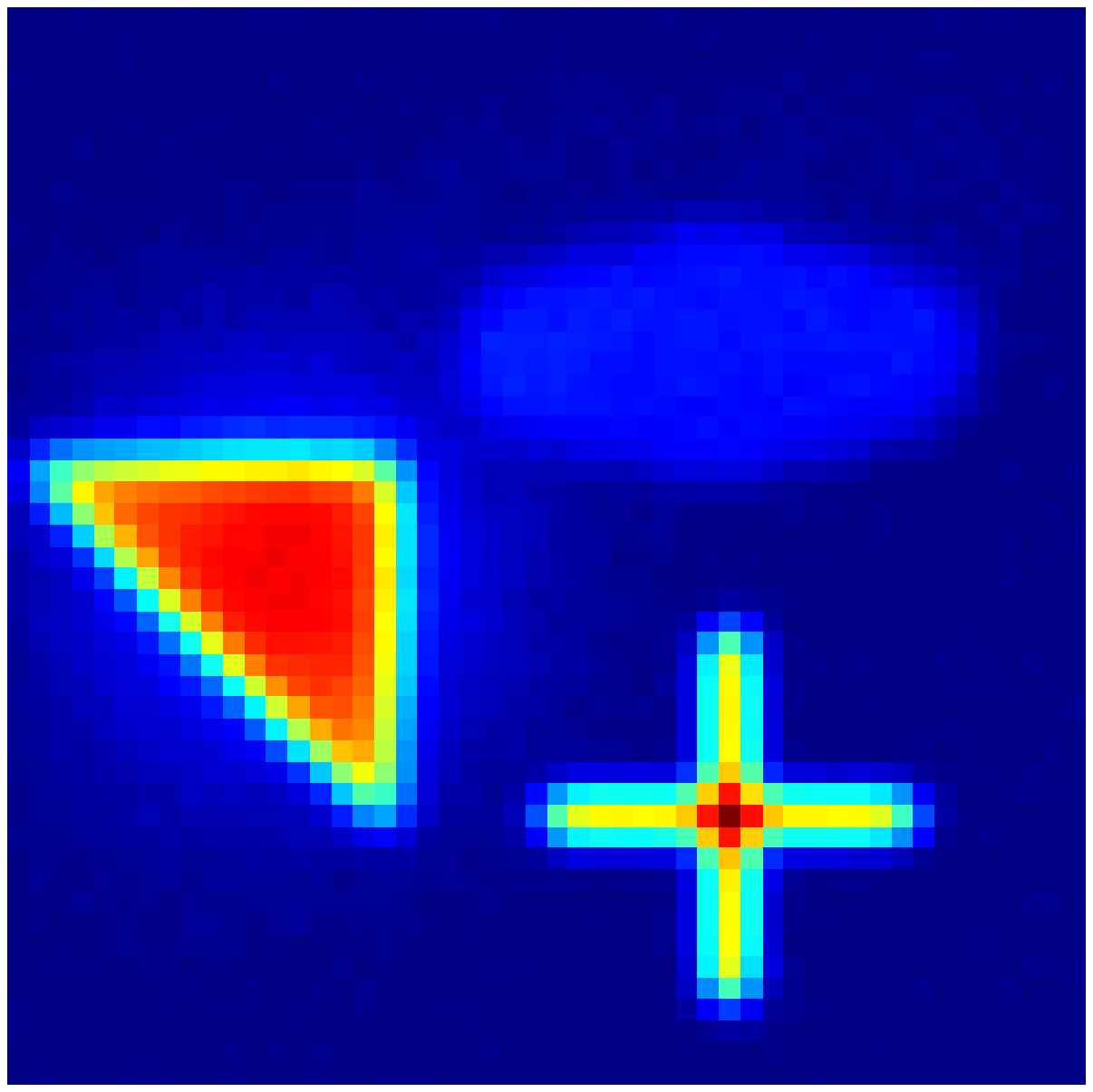} \hspace{-0.5cm} & 
\includegraphics[width=1.5cm]{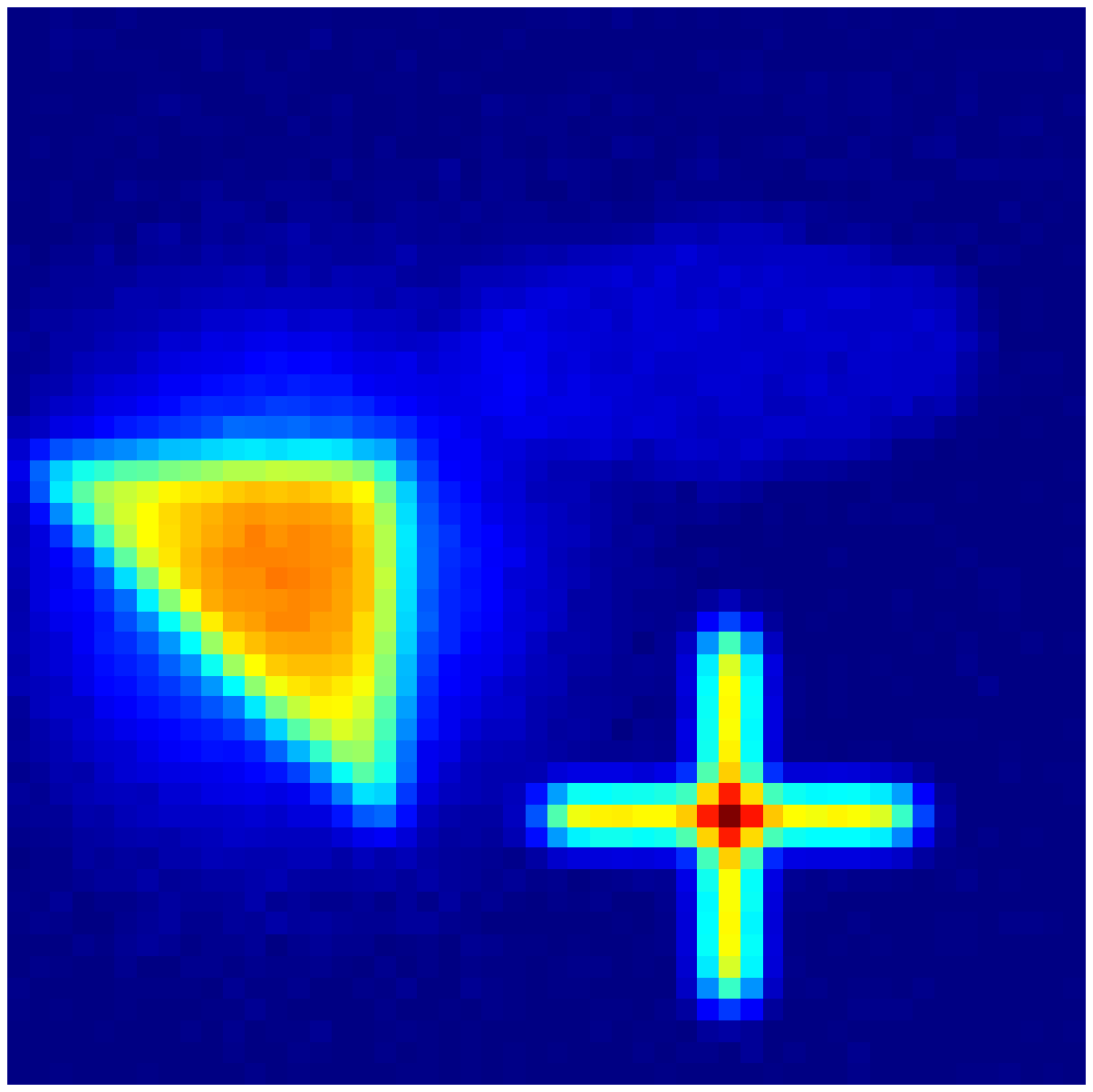} \hspace{-0.5cm} & \includegraphics[width=1.5cm]{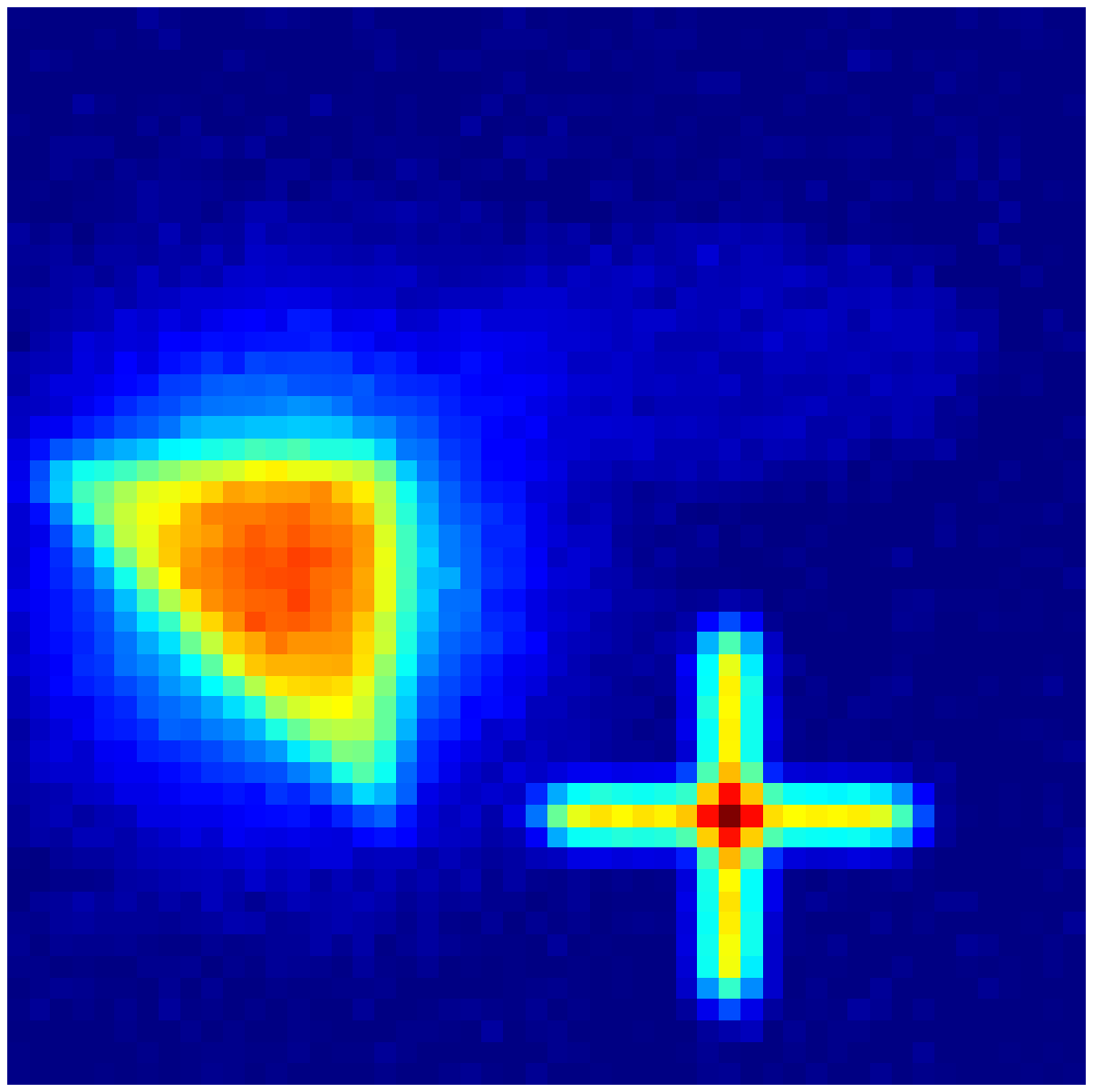}  
\end{tabular}
\caption{Dynamic image deblurring problem, where the true images are provided in the top row and the corresponding blurred observations are provided on the bottom row.  The first time point has been omitted for visualization ease.}\label{fig:DB_problem}
\end{center}
\end{figure}

This experiment concerns a synthetic dynamic image deblurring problem, where the goal is to reconstruct a sequence of images from the sequence of their corresponding blurred and noisy counterparts \cite{hansen1994regularization,chung2018efficient}. The blur matrix for this example is $\bfA = \bfA_t \kron \bfA_s$ where $\bfA_s$ represents a 2D Gaussian point spread function with spread parameter $\sigma = 1$ and bandwidth $4$ and $\bfA_t$ represents a 1D Gaussian blur with spread parameter $\sigma =1$ and bandwidth $3$. Gaussian white noise has been added to the measurements with noise level $\|\bfe\| / \|\bfA \bfx_{\rm true}\| =0.02$.  The sequence of images forming $\bfx_{\rm true}$ is displayed in the first row of Figure \ref{fig:DB_problem}, while the sequence of their corresponding blurred and noisy counterparts $\bfb$ can be found in the second row. Note that $\bfx$ and $\bfb$ are spatio-temporal images of size $50$$\times$$50$$\times$$9$ (i.e., $50$$\times$$50$ pixels at 9 time points).

\begin{figure}[htb!]
     \centering
    \begin{subfigure}[b]{0.45\textwidth}
         \centering
         \includegraphics[height =4.7cm]{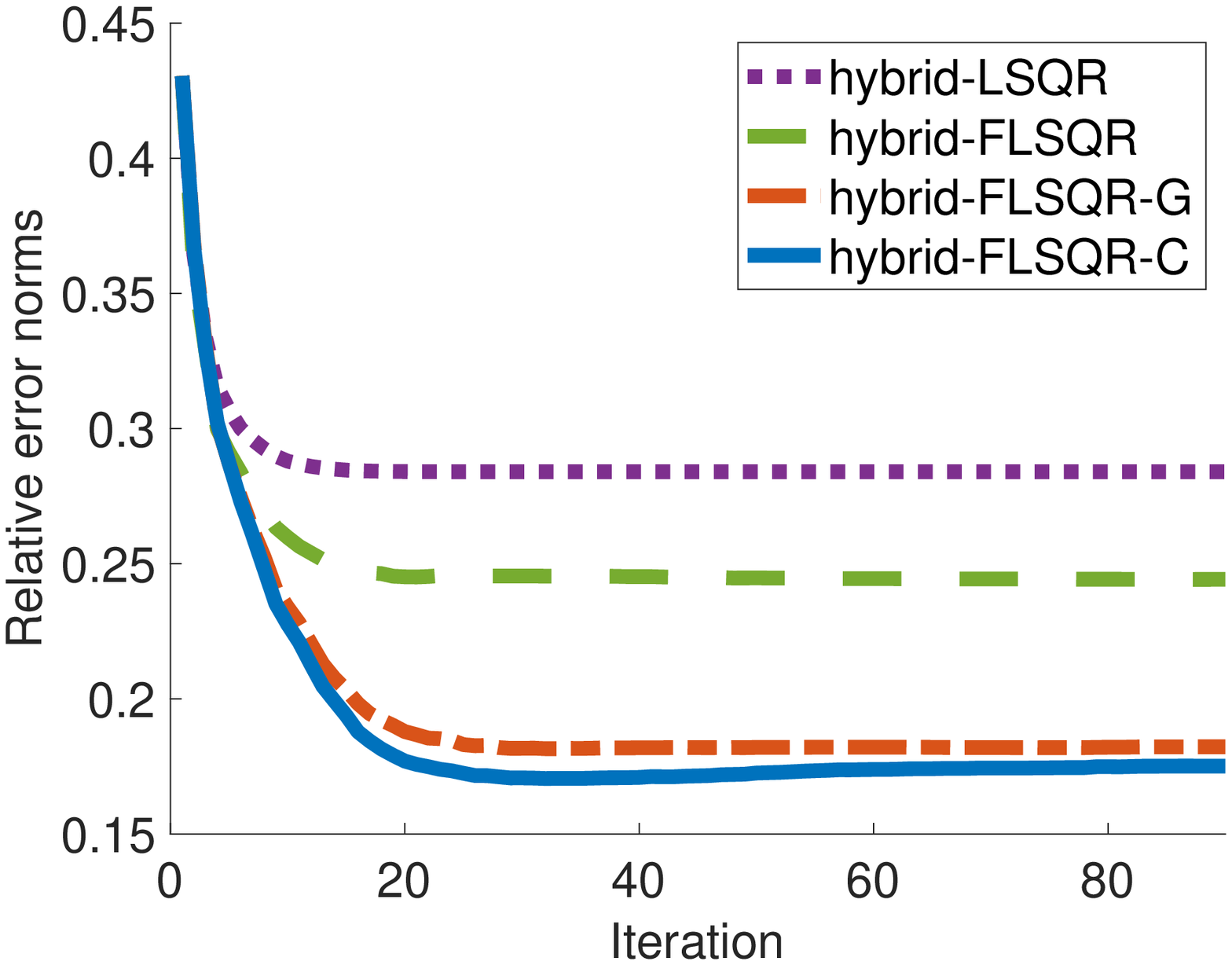}
         \caption{Error norms for methods based on LSQR}
         \label{fig:errors_deblur_dynamic_LSQR}
     \end{subfigure}
     \hspace{1cm}
     \begin{subfigure}[b]{0.45\textwidth}
         \centering
         \includegraphics[height =4.7cm]{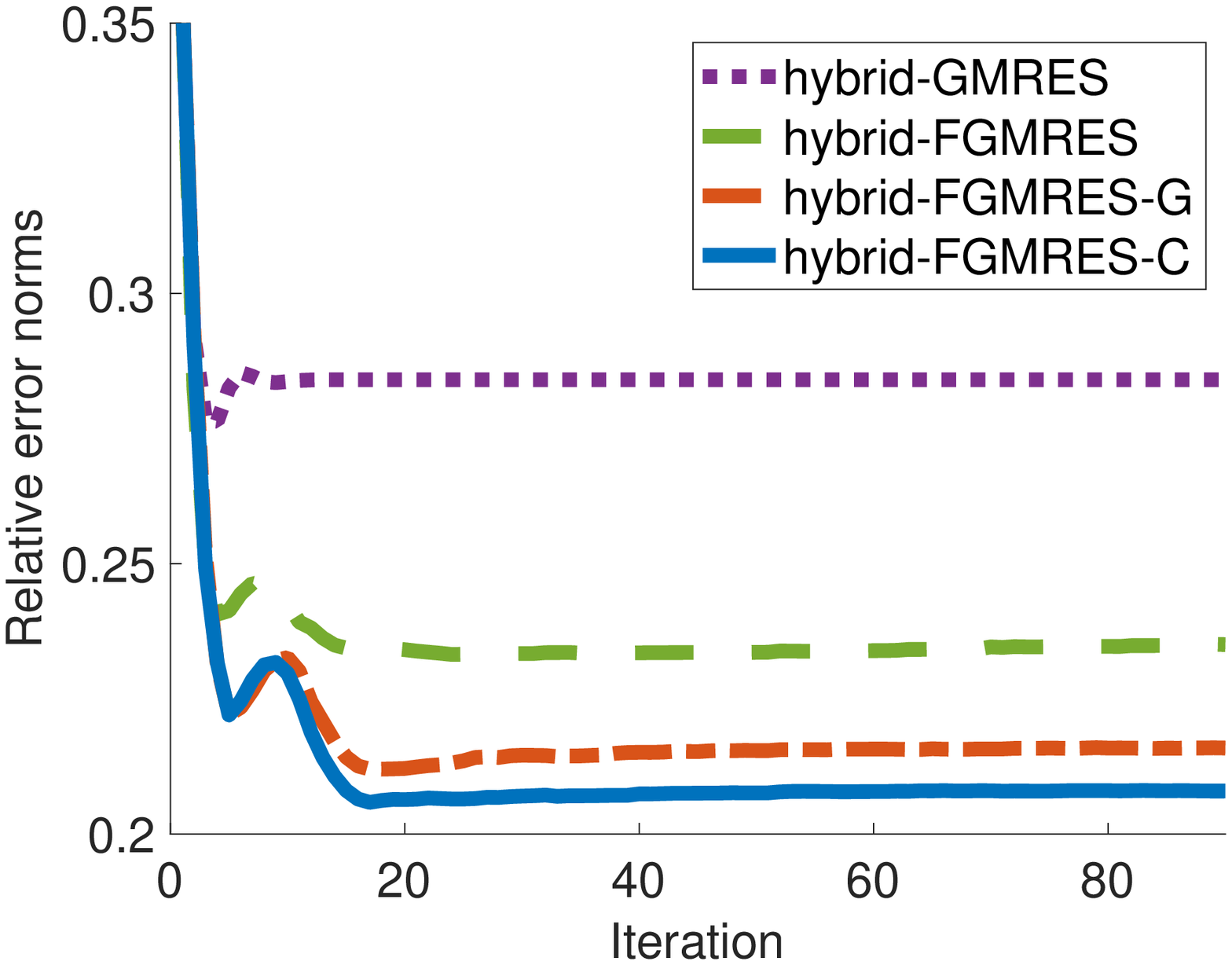}
         \caption{Error norms for methods based on GMRES}
         \label{fig:errors_deblur_dynamic_GMRES}
     \end{subfigure}
    \caption{Relative error norms for the methods based on FLSQR and FGMRES for the deblurring example presented in Figure \ref{fig:DB_problem}. The regularization parameter $\lambda_k$ has been chosen at each iteration using the DP.}
        \label{fig:DB_err}
\end{figure}

In this example, some regions in all of the images in the sequence are always zero while other regions have pixels of changing intensity. It is then appropriate to use group sparsity regularization where each spatial location, over all time points, constitutes a group.
We compare the results for $\ell_2$, $\ell_1$ and $\ell_{2,1}$ regularization using hybrid flexible methods, and we present results for the new method combining an $\ell_1$ and an $\ell_{2,1}$  regularization term. Since the forward model matrix $\bfA$ is square, methods based on both LSQR and GMRES are tested. The error norm plots for all methods can be found in Figure \ref{fig:DB_err}. All of the presented approaches used the DP as defined in \eqref{dp} to select the regularization parameter, where $\bfx_k$ in \eqref{dp} corresponds to the approximated solution of \eqref{flsqr-C} computed at iteration $k$. Note that $\tau_\lambda$ in \eqref{eq:tau_lambda} must be set ahead of the iterations. In particular, $\tau_\lambda = 1.2$ for methods based on FLSQR and $\tau_\lambda = 0.8$ for methods based on FGMRES.

It can be observed from the plot of relative reconstruction error norms, see Figure \ref{fig:DB_err}, that $\ell_1$ regularization (hybrid-FLSQR and hybrid-FGMRES) outperforms Tikhonov regularization (hybrid-LSQR and hybrid-GMRES) for this problem. Moreover, the new algorithms enforcing group sparsity regularization (hybrid-FLSQR-G and hybrid-FGMRES-G) produce significantly better reconstructions (i.e., smaller reconstruction errors), with the combined approaches (hybrid-FLSQR-C and hybrid-FGMRES-C) performing marginally better than group sparsity regularization on its own. 

\begin{figure}[ht]
\begin{center}
\begin{tabular}{ c c c c c c c c c c c}
\raisebox{1.5\normalbaselineskip}[0pt][0pt]{\rotatebox[origin=c]{90}{\footnotesize  hybrid}}  
\raisebox{1.5\normalbaselineskip}[0pt][0pt]{\rotatebox[origin=c]{90}{\footnotesize  FLSQR}} 
\includegraphics[width=1.5cm]{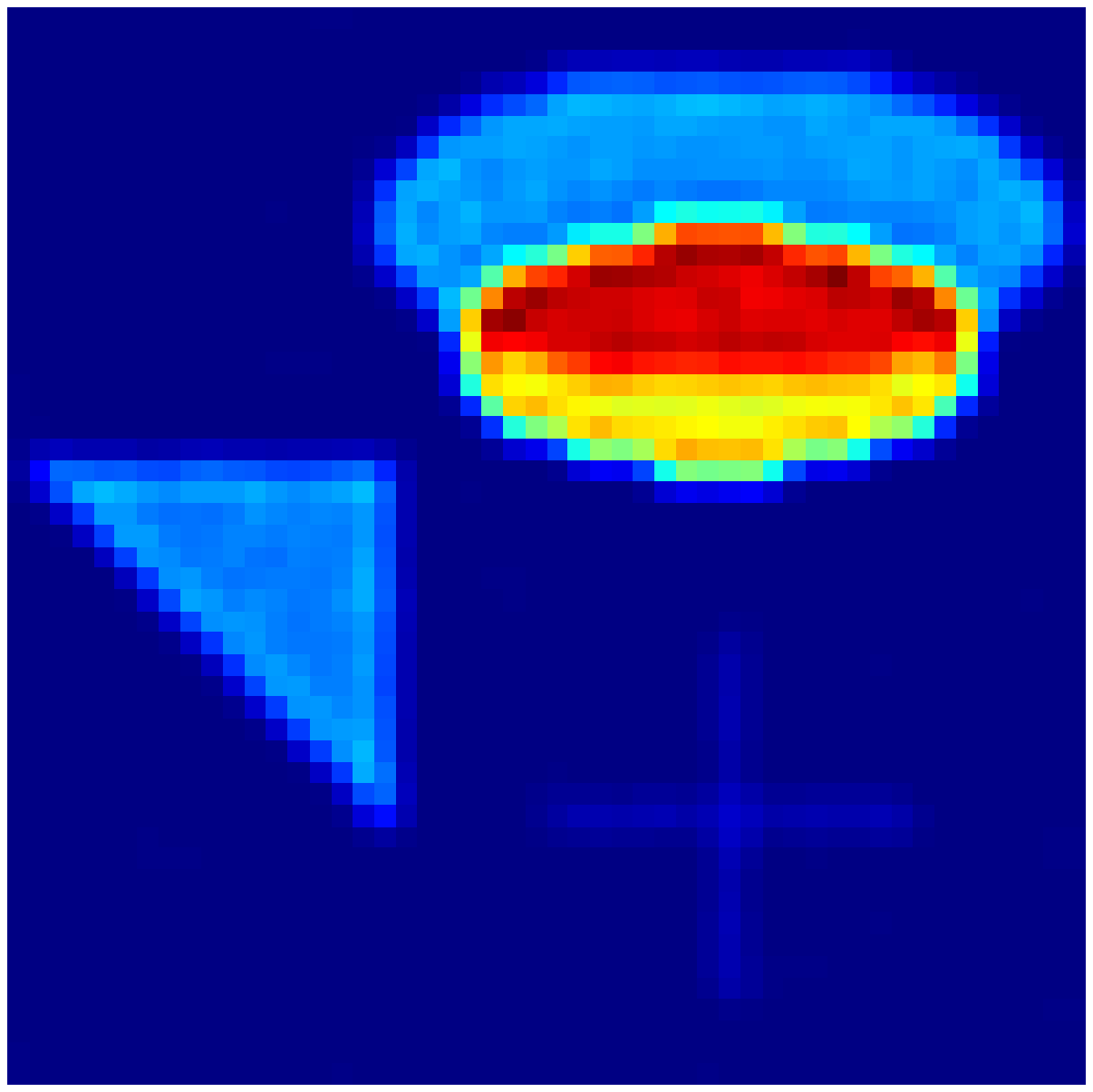} \hspace{-0.5cm} & \includegraphics[width=1.5cm]{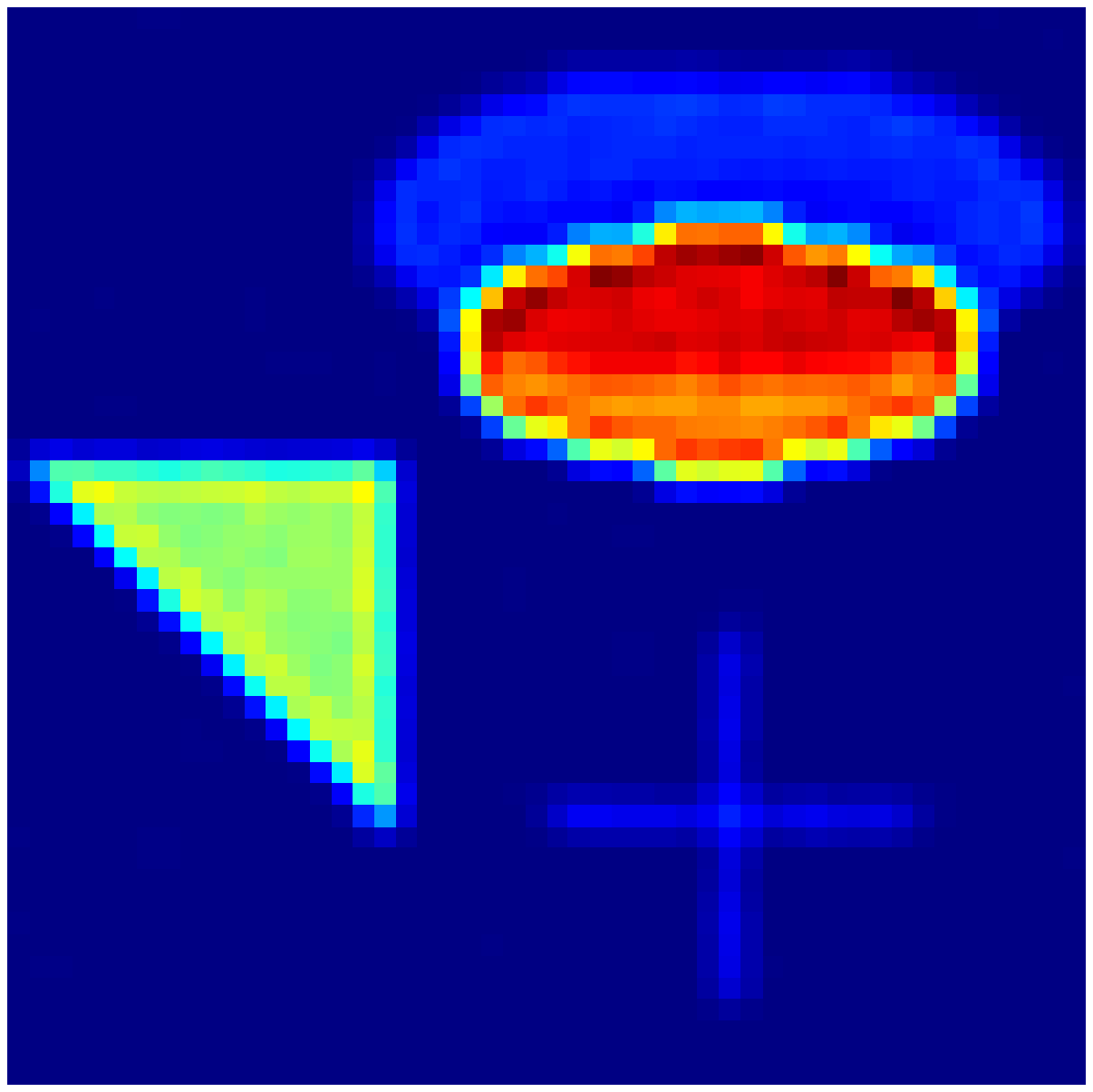} \hspace{-0.5cm} & \includegraphics[width=1.5cm]{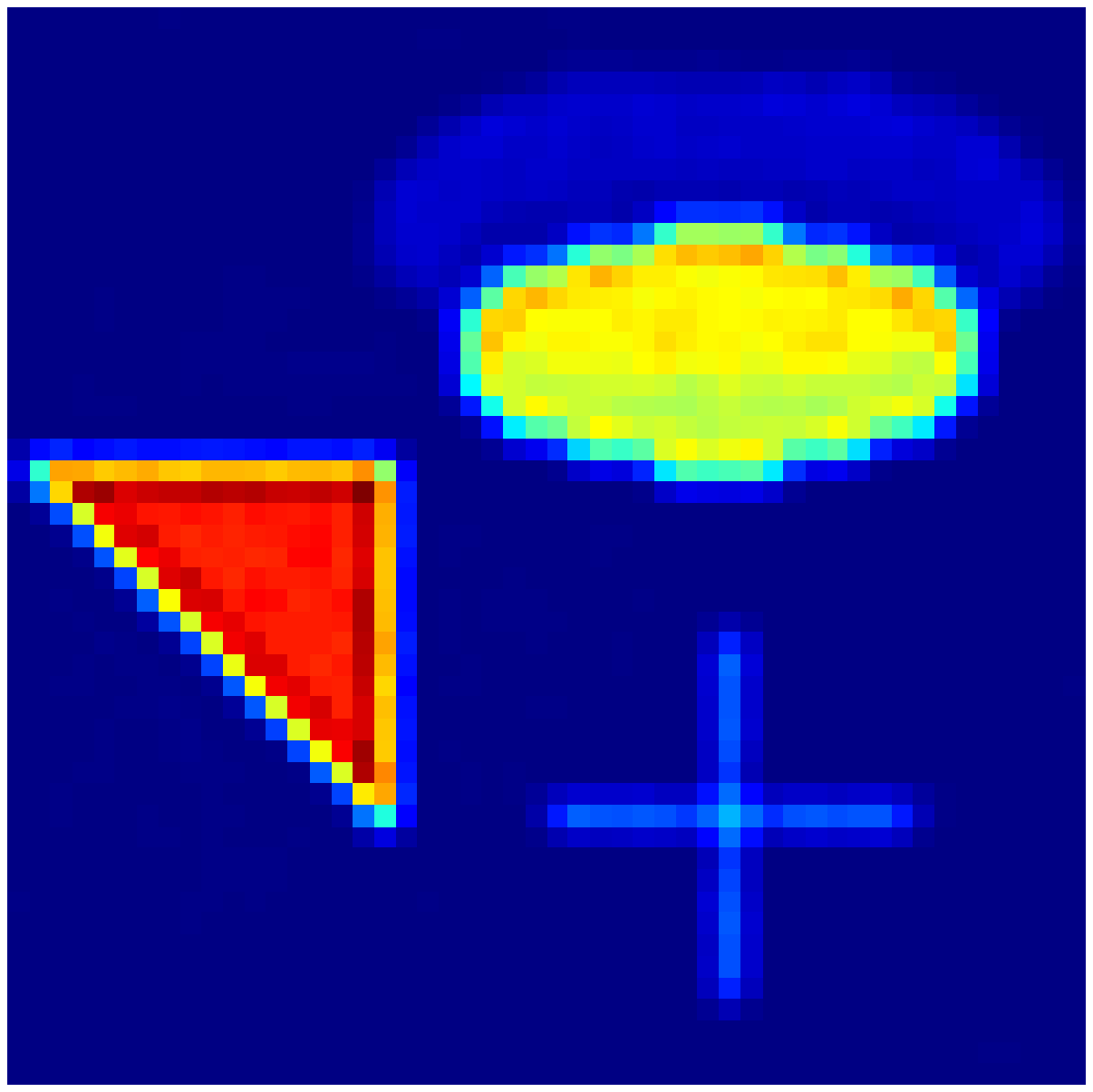} \hspace{-0.5cm} & \includegraphics[width=1.5cm]{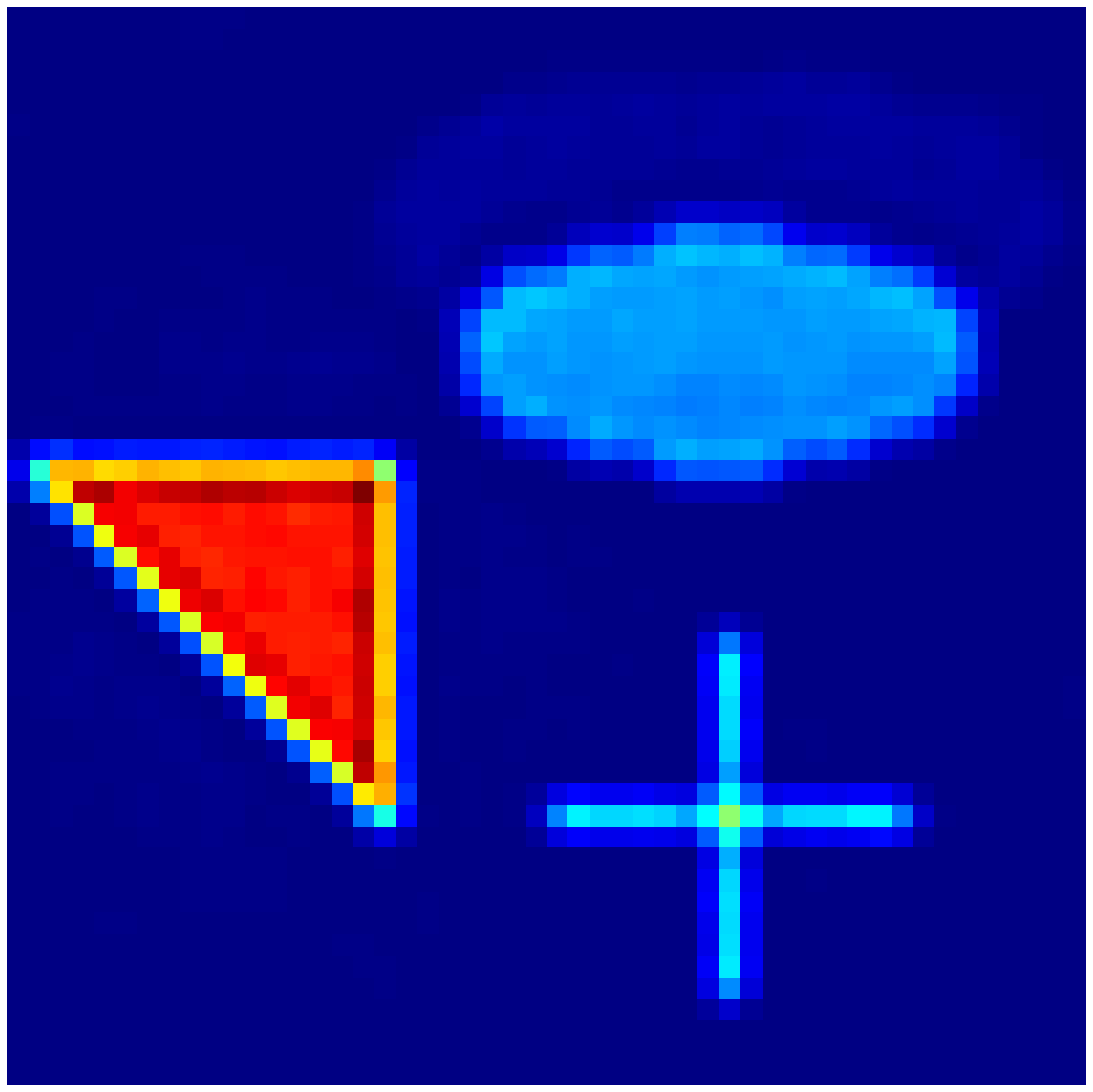} \hspace{-0.5cm} &  \includegraphics[width=1.5cm]{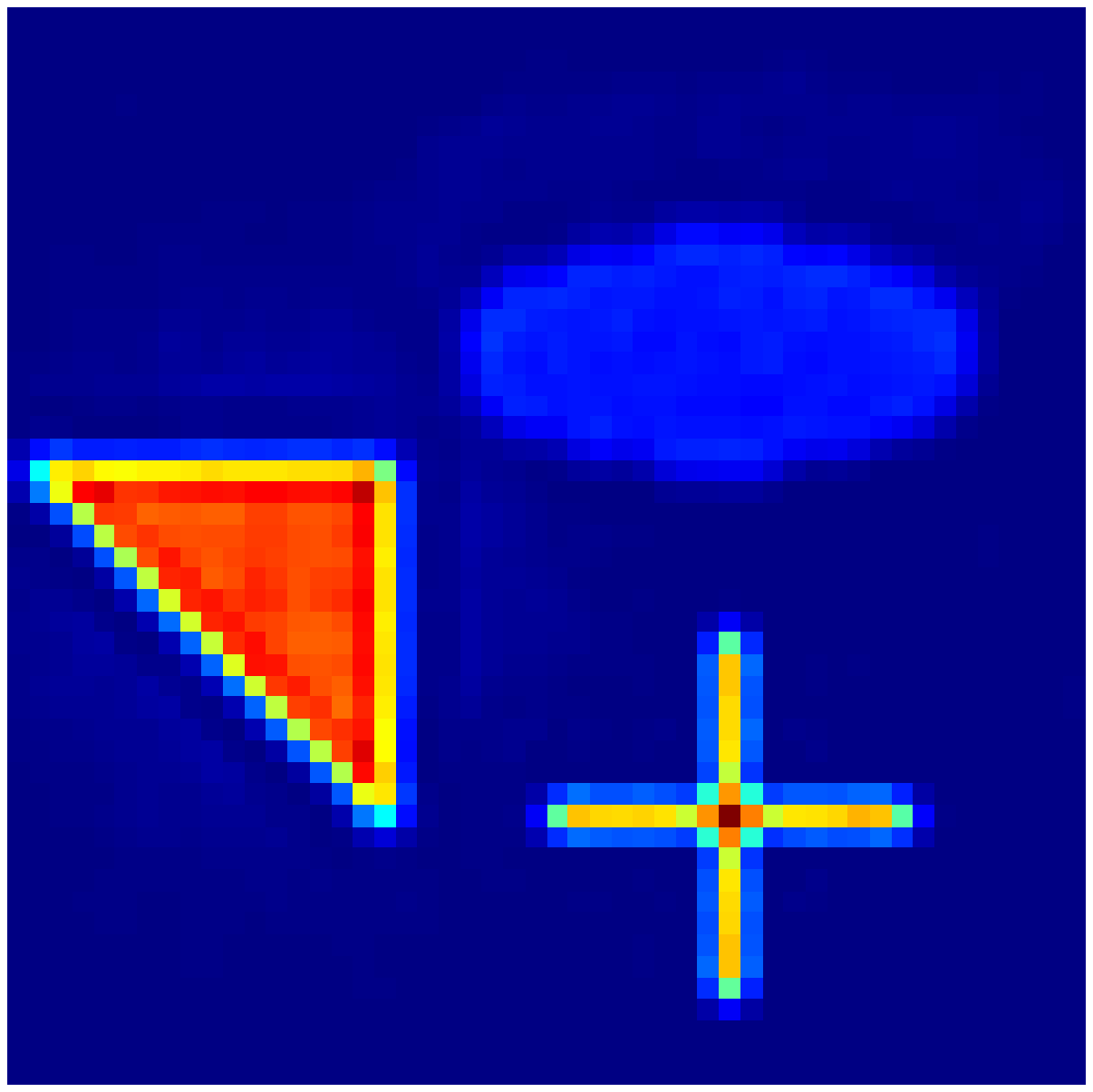} \hspace{-0.5cm} & 
\includegraphics[width=1.5cm]{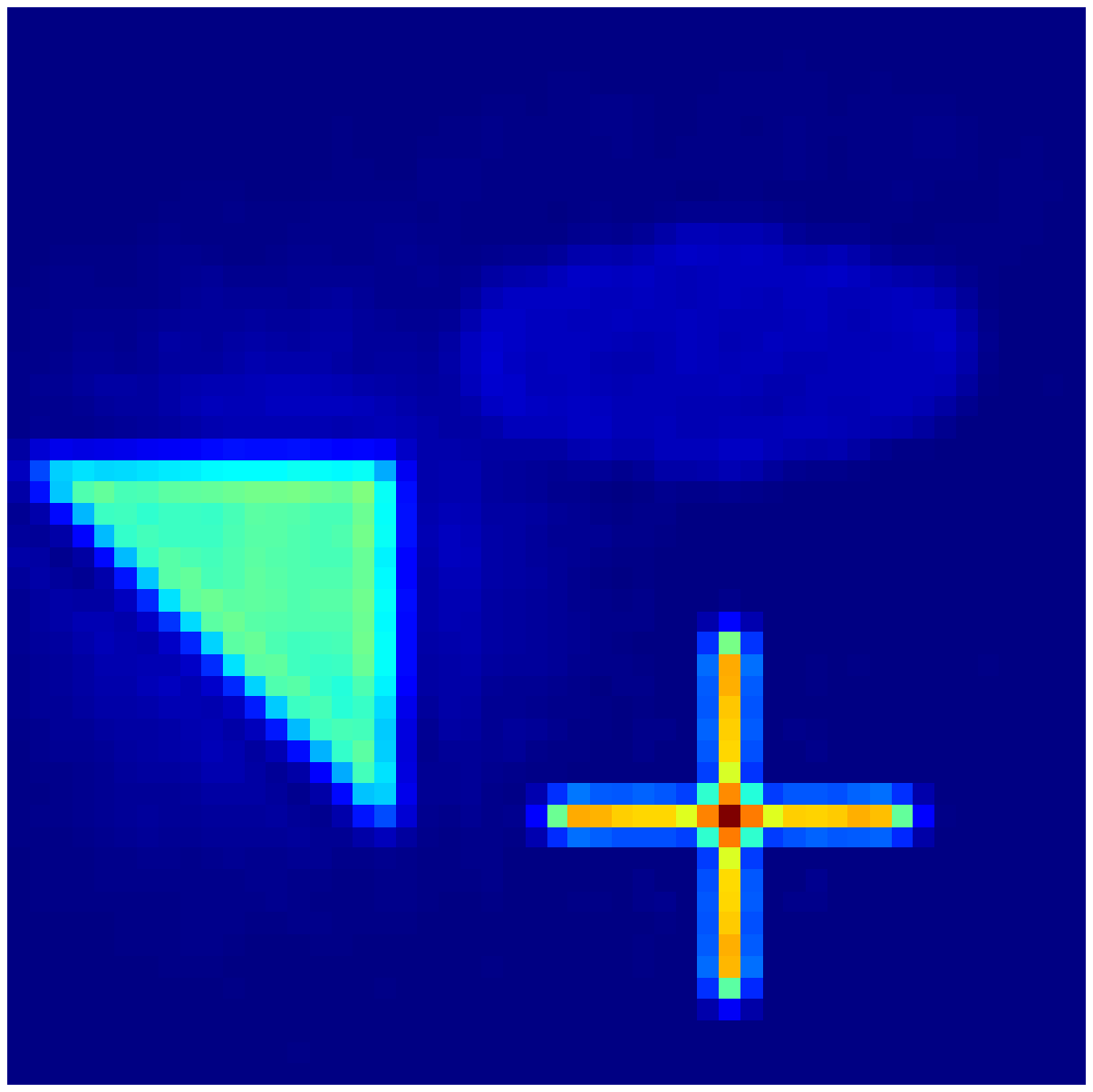} \hspace{-0.5cm} & 
\includegraphics[width=1.5cm]{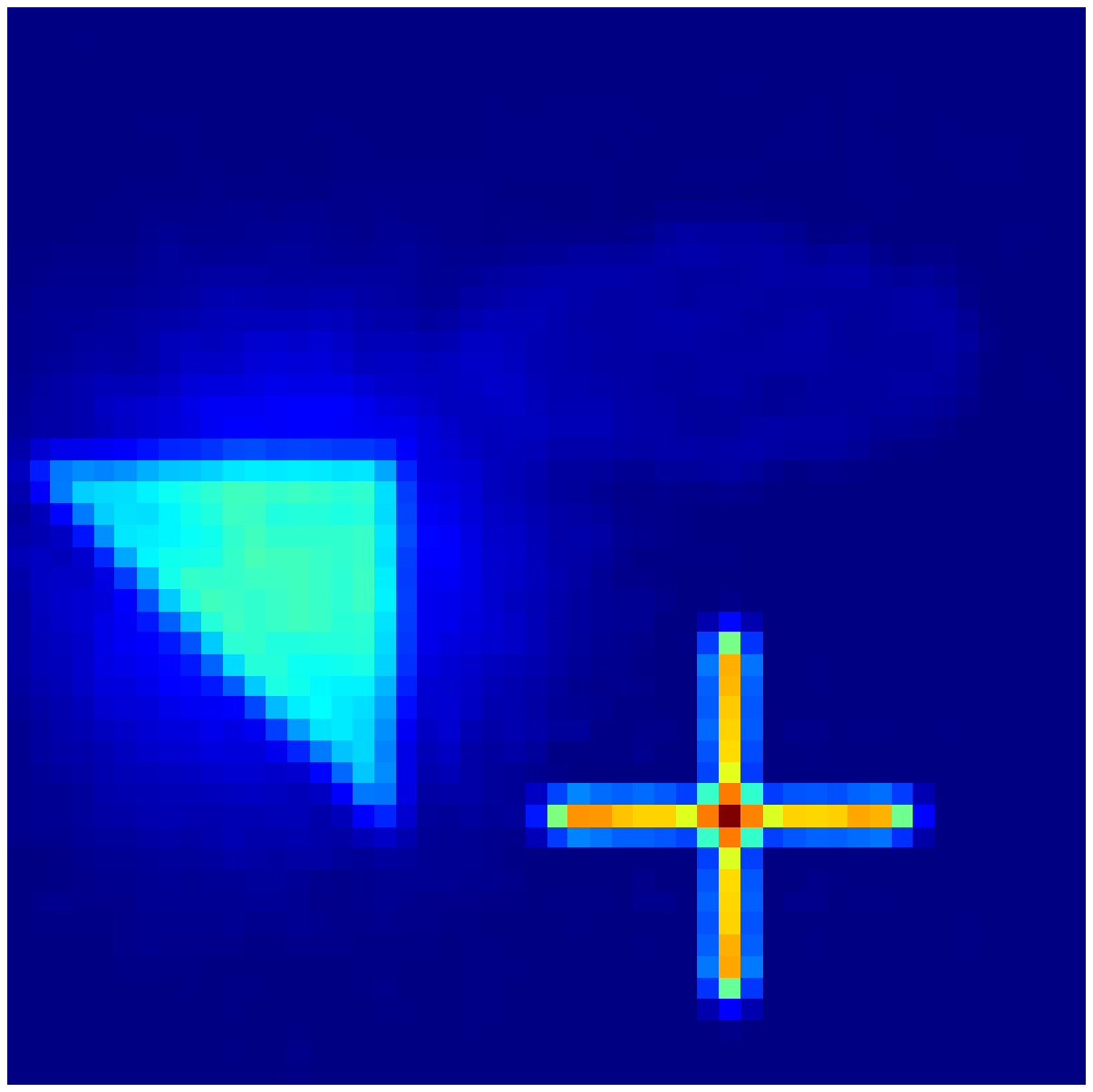} \hspace{-0.5cm} & \includegraphics[width=1.5cm]{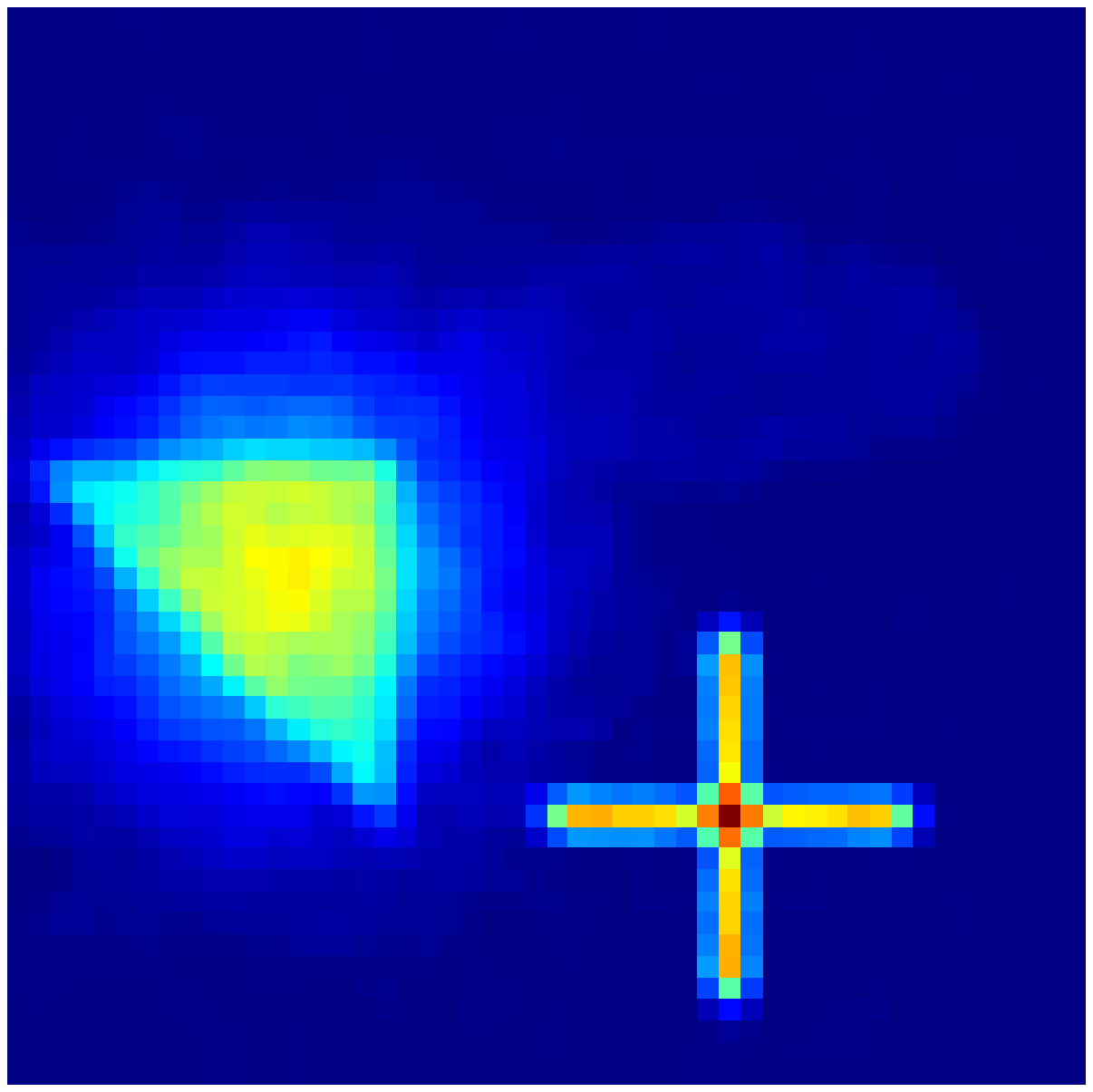}  \\ 
\raisebox{1.5\normalbaselineskip}[0pt][0pt]{\rotatebox[origin=c]{90}{\footnotesize  hybrid}}  
\raisebox{1.5\normalbaselineskip}[0pt][0pt]{\rotatebox[origin=c]{90}{\footnotesize FLSQR-C}} 
\includegraphics[width=1.5cm]{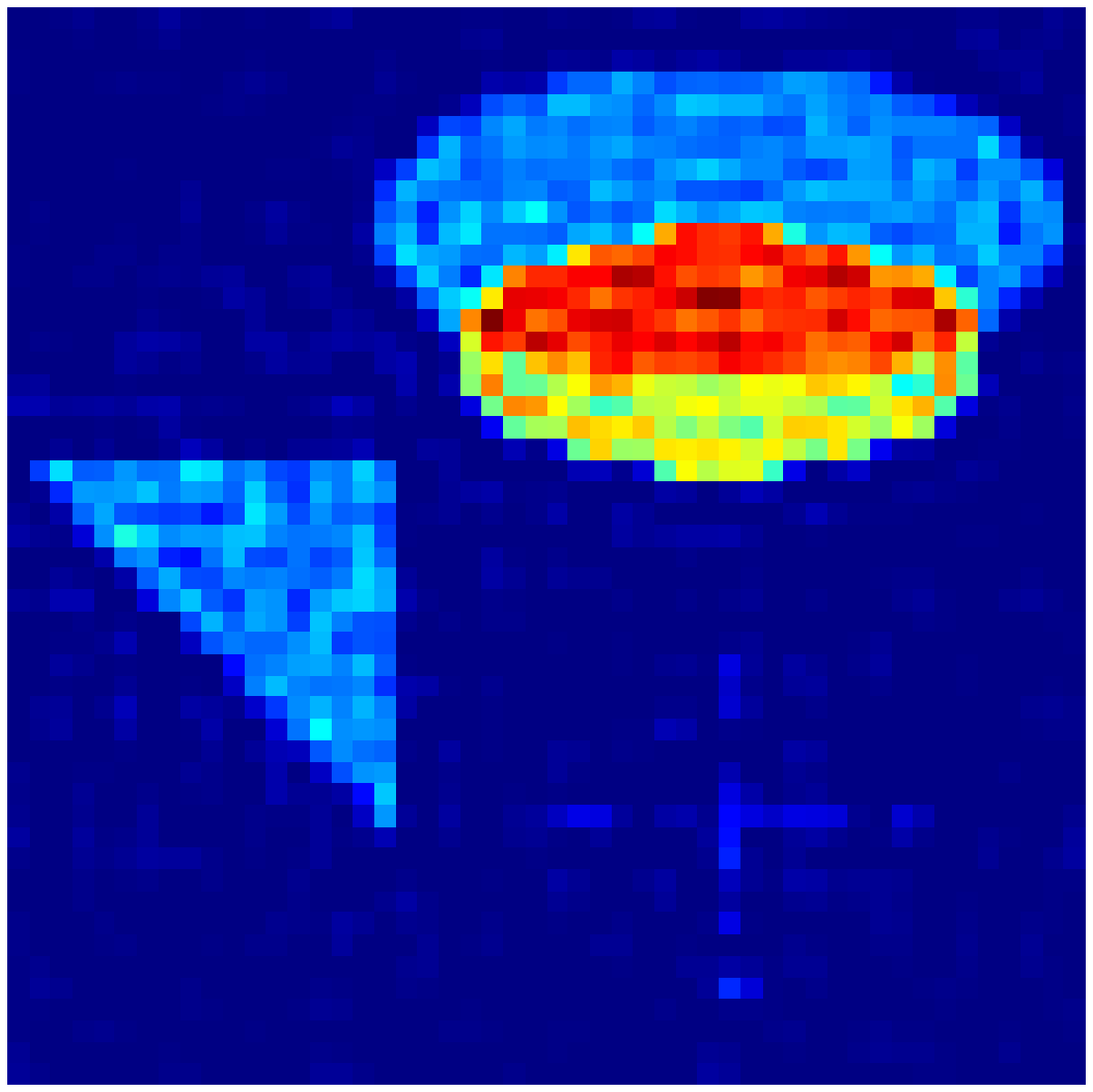} \hspace{-0.5cm} & \includegraphics[width=1.5cm]{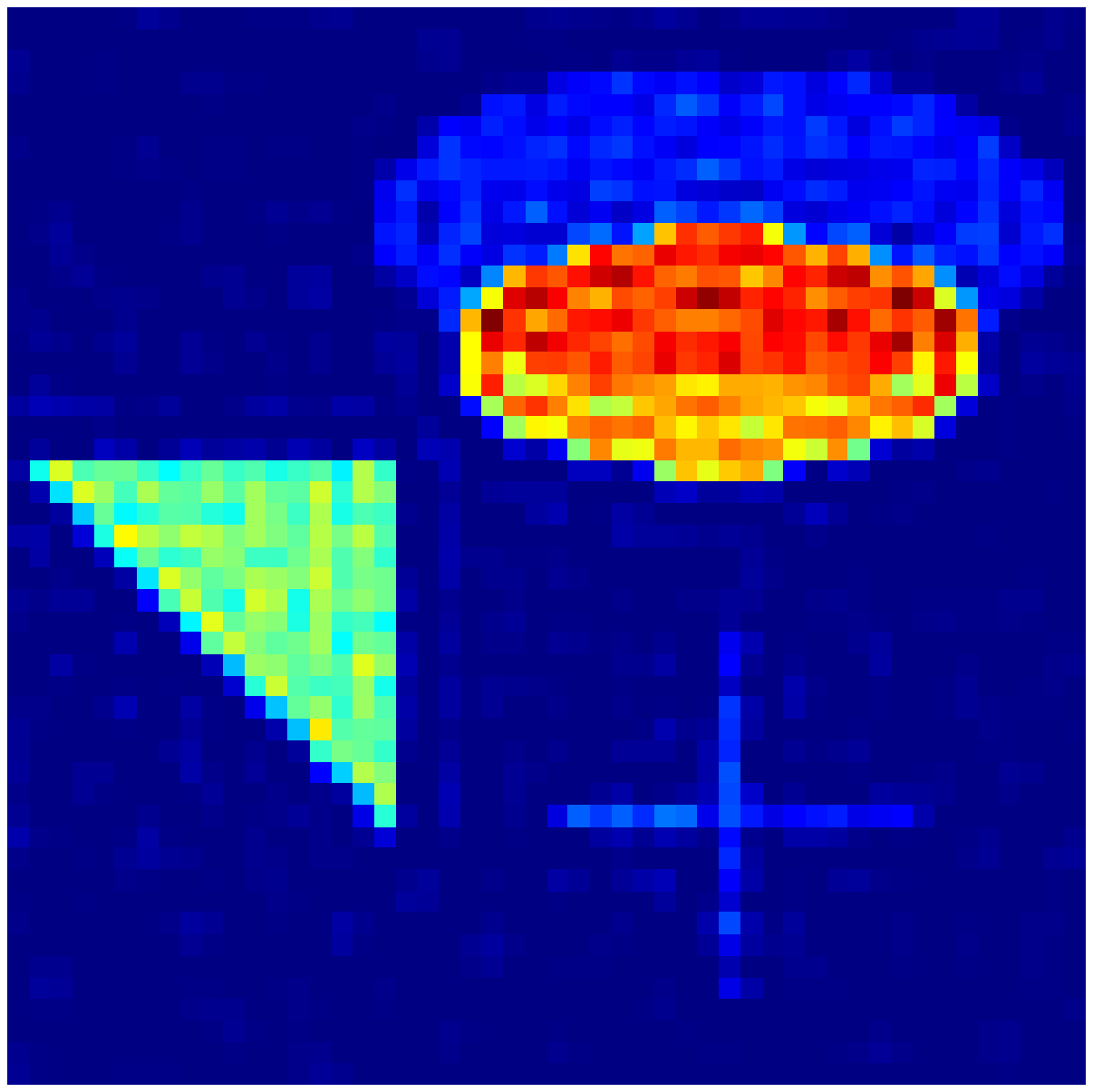} \hspace{-0.5cm} & \includegraphics[width=1.5cm]{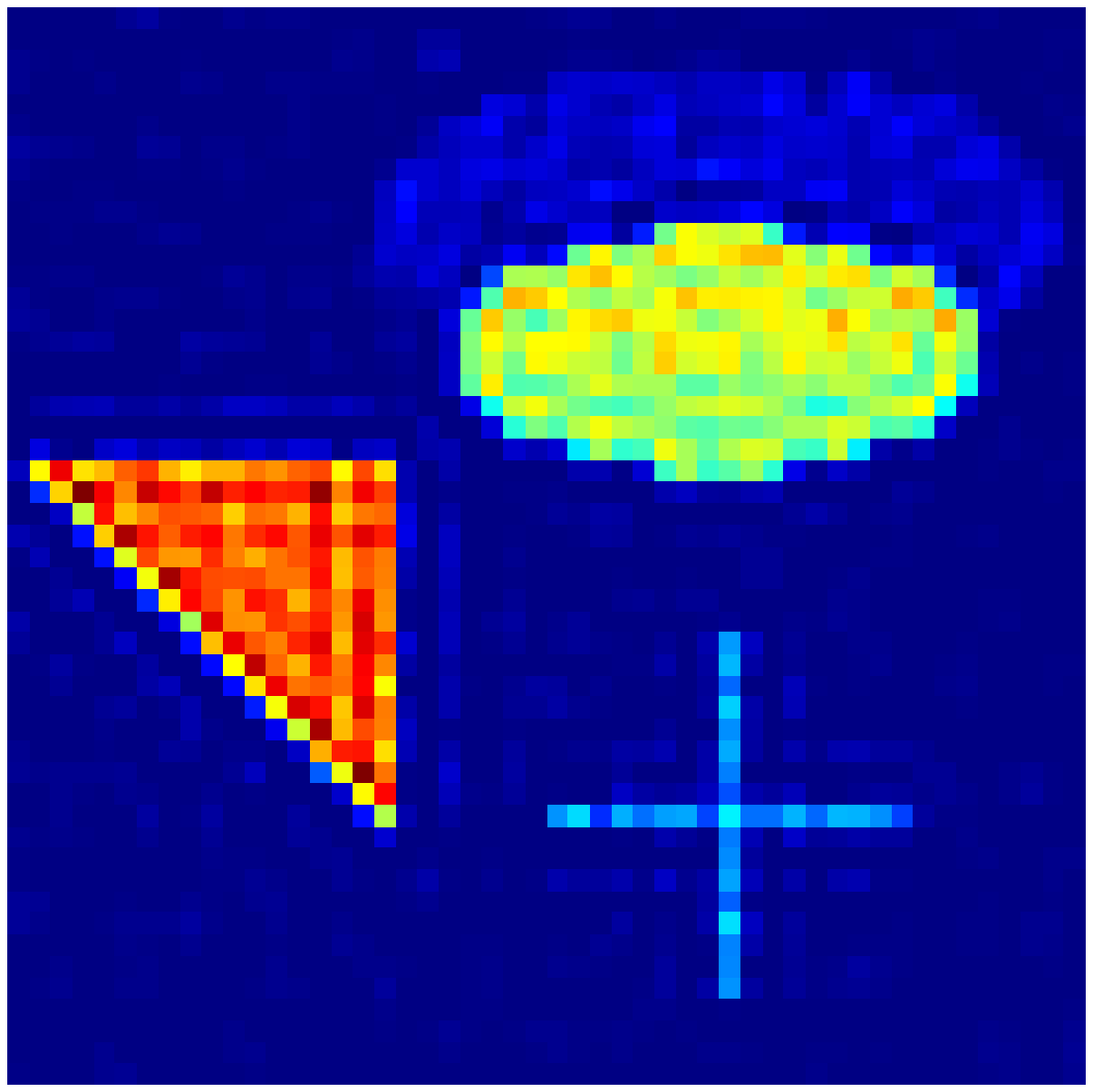} \hspace{-0.5cm} & \includegraphics[width=1.5cm]{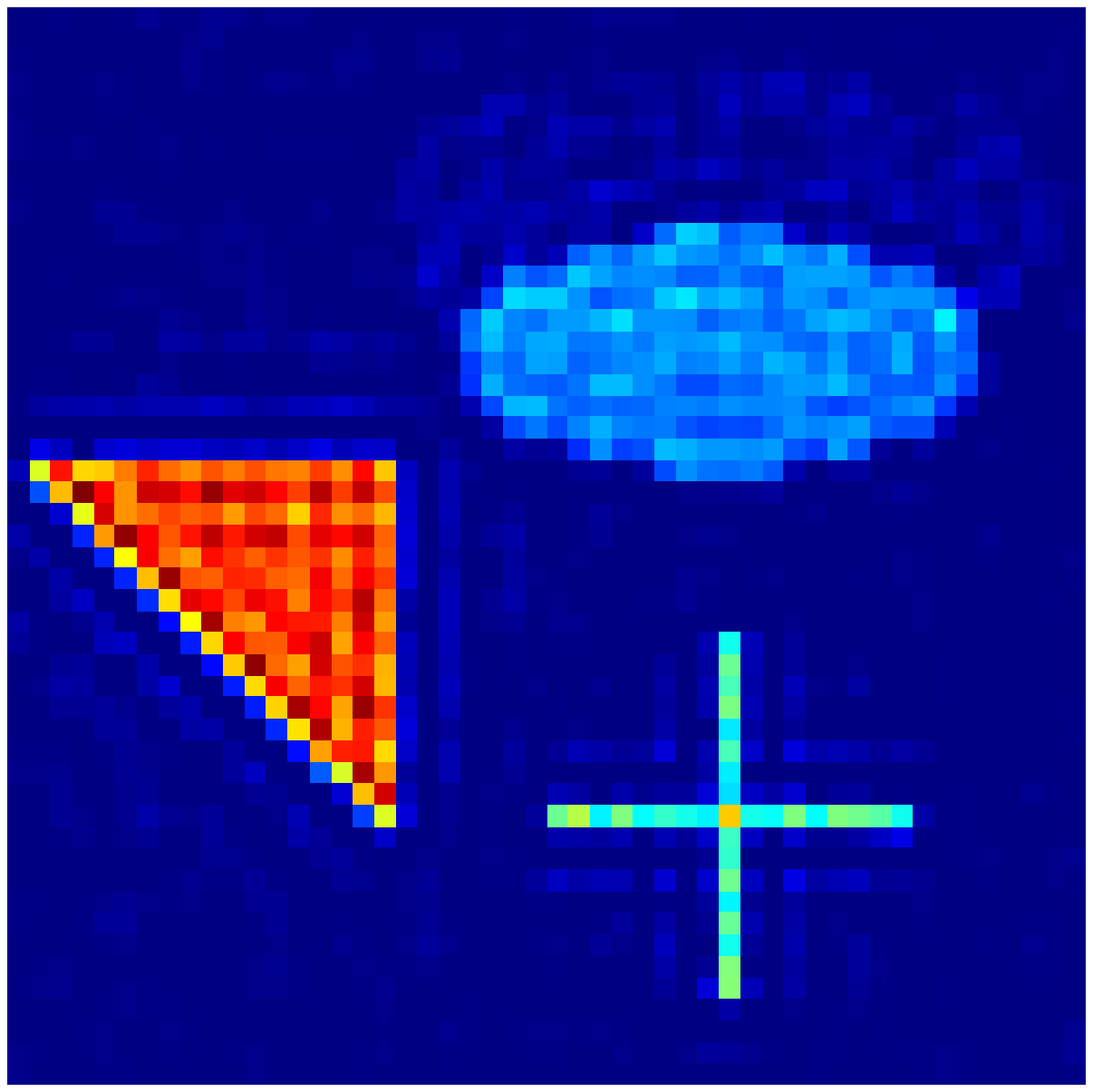} \hspace{-0.5cm} &  \includegraphics[width=1.5cm]{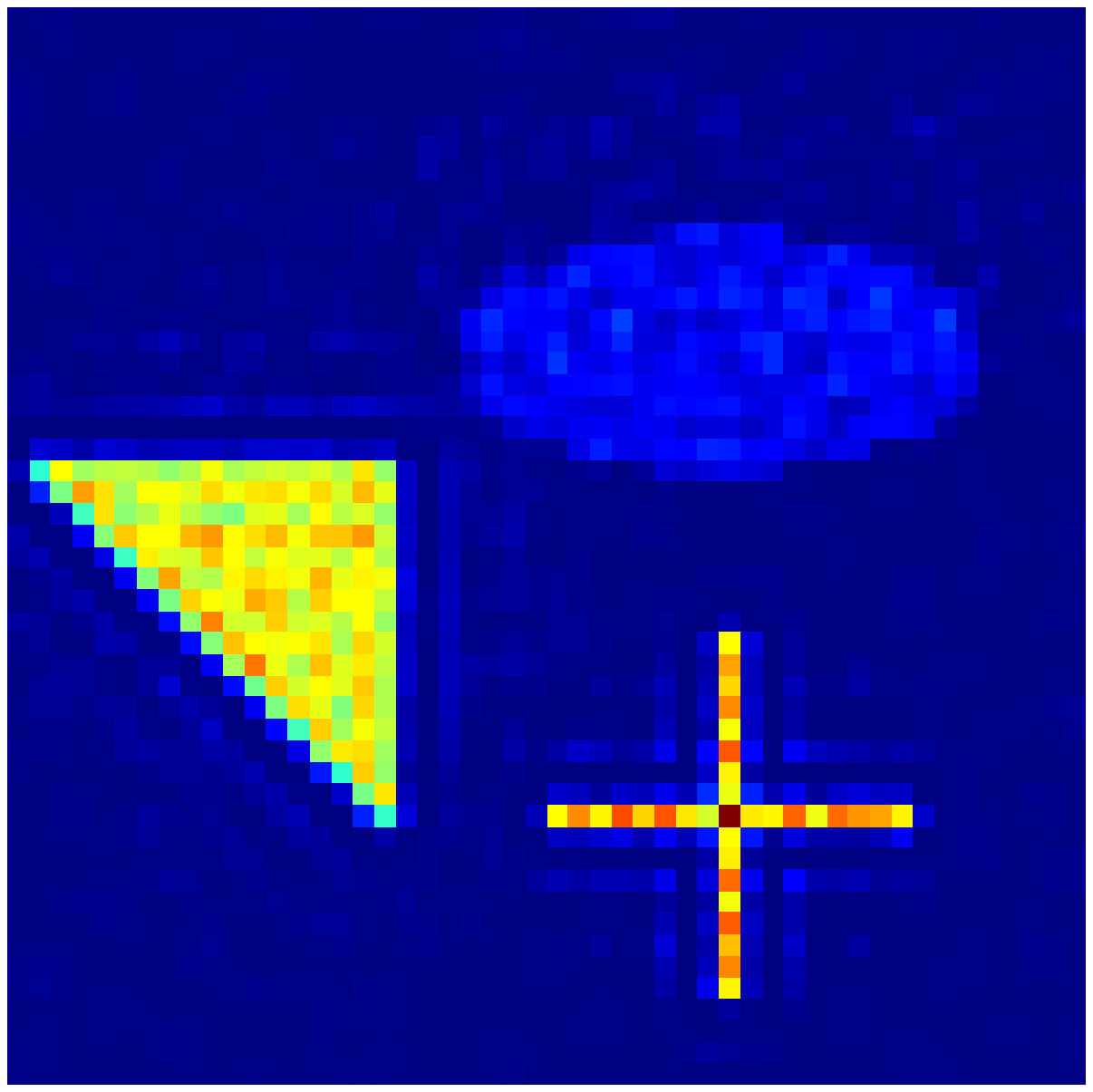} \hspace{-0.5cm} & 
\includegraphics[width=1.5cm]{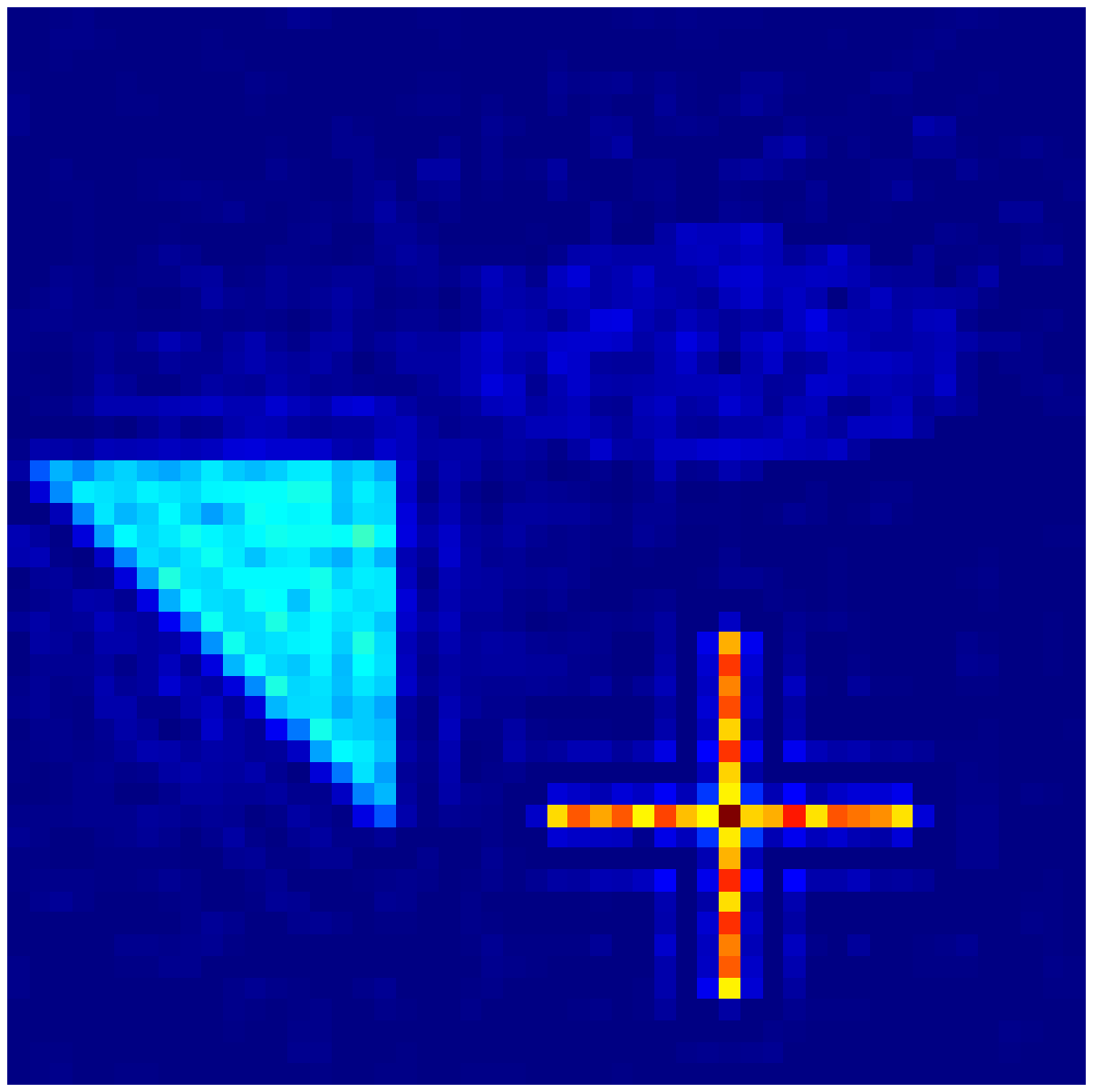} \hspace{-0.5cm} & 
\includegraphics[width=1.5cm]{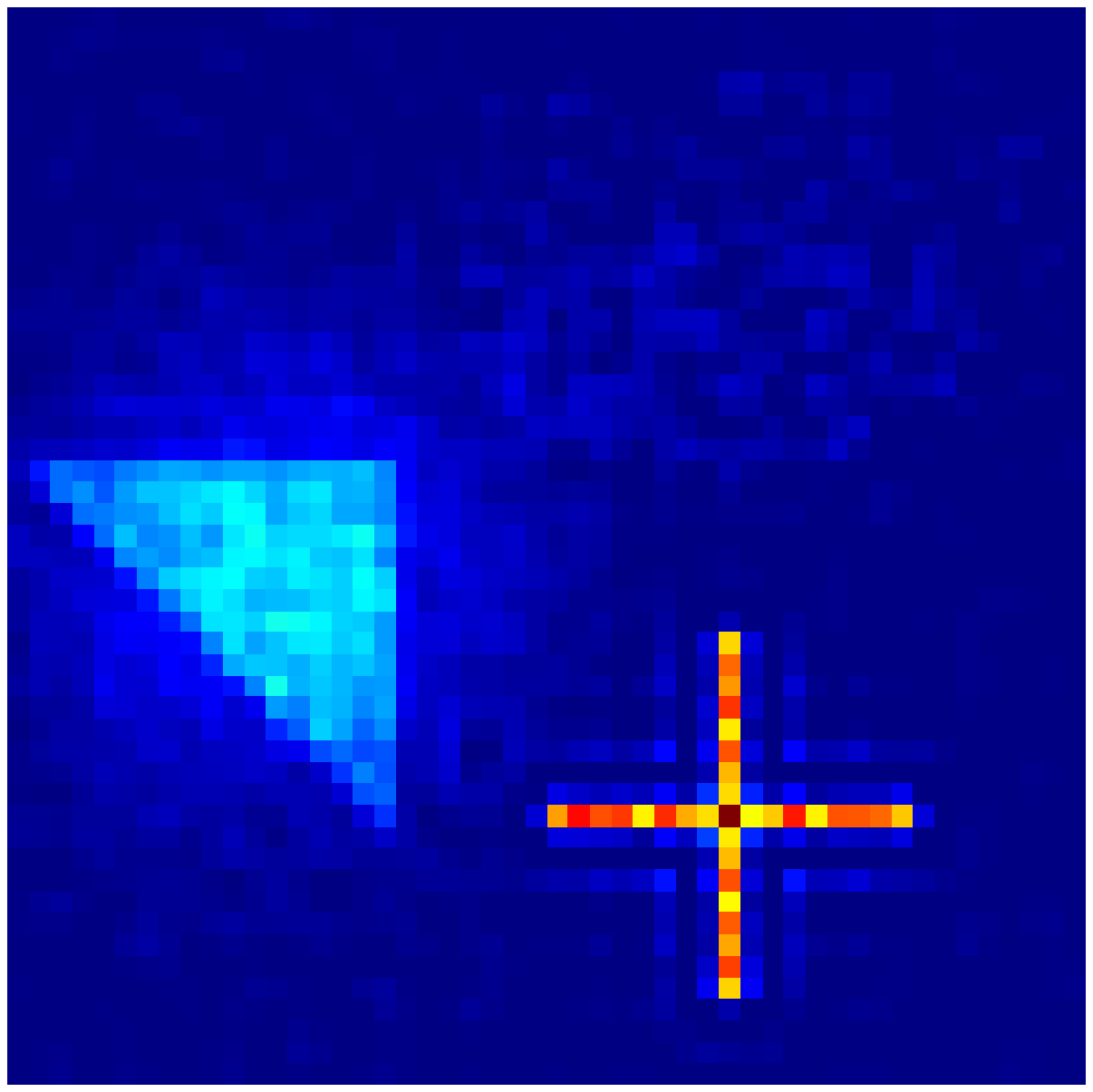} \hspace{-0.5cm} & \includegraphics[width=1.5cm]{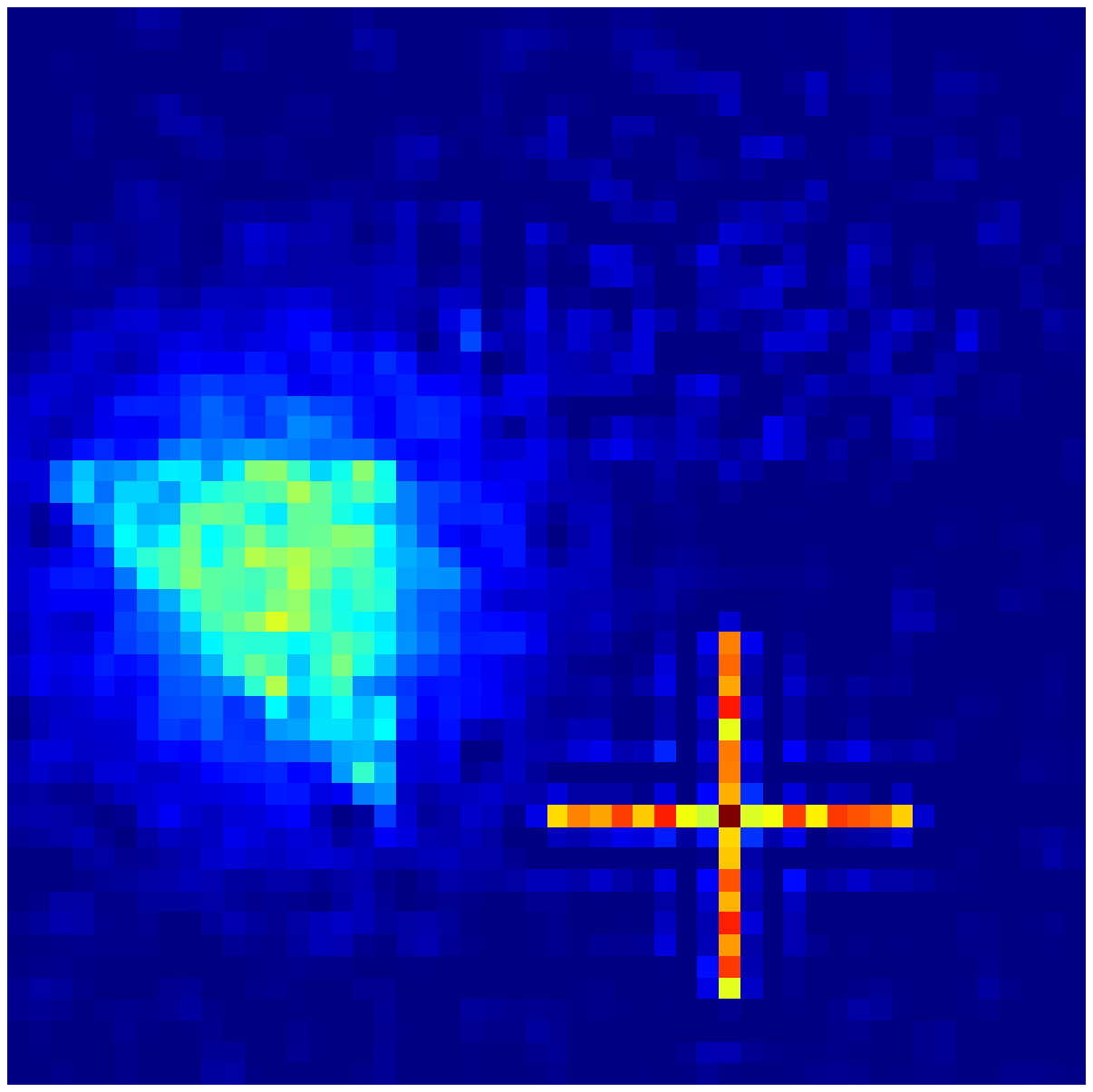}  
\end{tabular}
\caption{Reconstructions with $\ell_1$ regularization and with a combination of $\ell_1$ and $\ell_{2,1}$ regularization for the dynamic image deblurring problem.  These methods are based on FLSQR, and the first time point has been omitted for visualization ease.}\label{fig:DB_recon1}
\end{center}
\end{figure}

\begin{figure}[ht]
\begin{center}
\begin{tabular}{ c c c c c c c c c c c}
\raisebox{1.5\normalbaselineskip}[0pt][0pt]{\rotatebox[origin=c]{90}{\footnotesize hybrid}}  
\raisebox{1.5\normalbaselineskip}[0pt][0pt]{\rotatebox[origin=c]{90}{\footnotesize FGMRES}} 
\includegraphics[width=1.5cm]{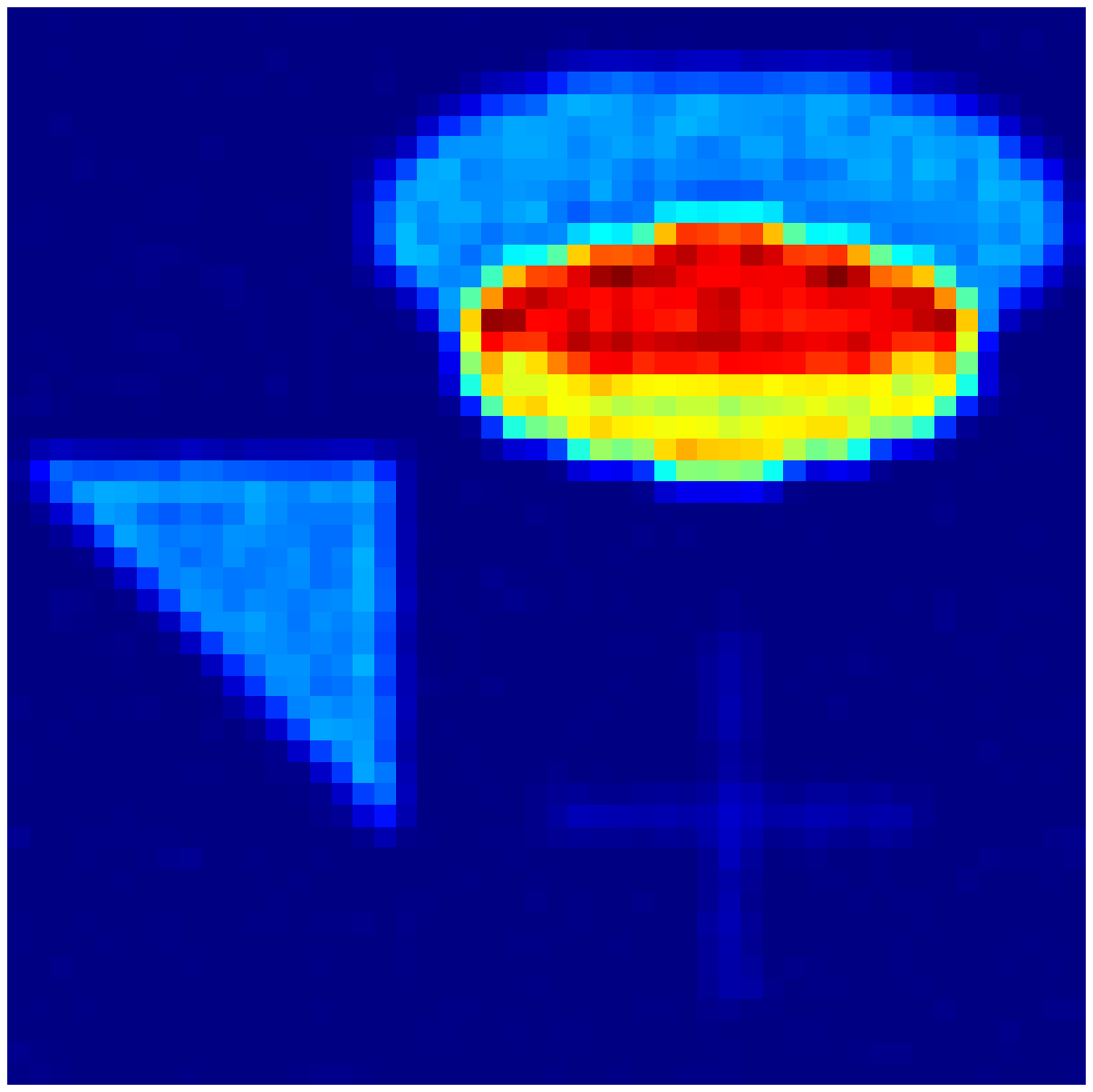} \hspace{-0.5cm} & \includegraphics[width=1.5cm]{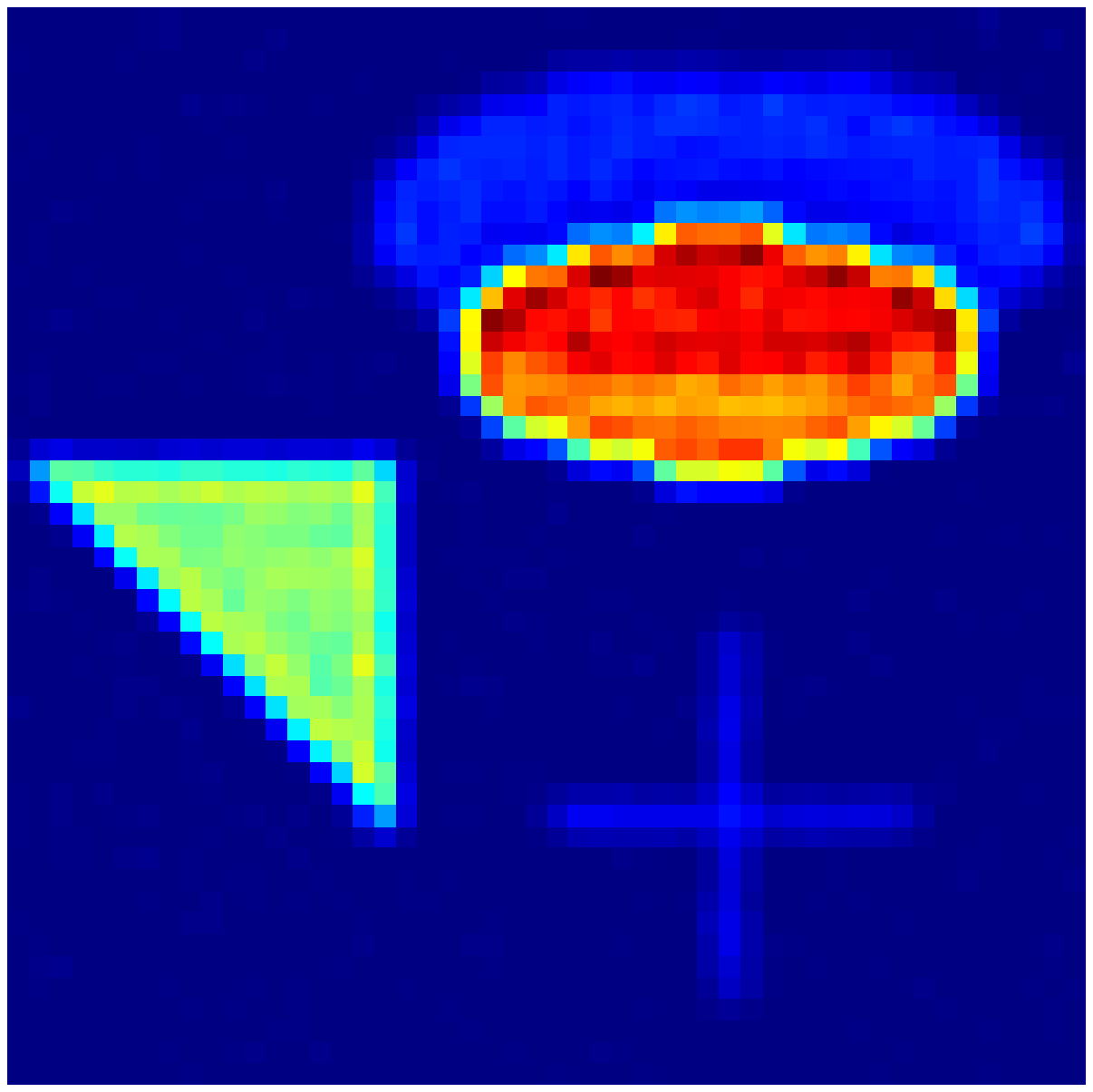} \hspace{-0.5cm} & \includegraphics[width=1.5cm]{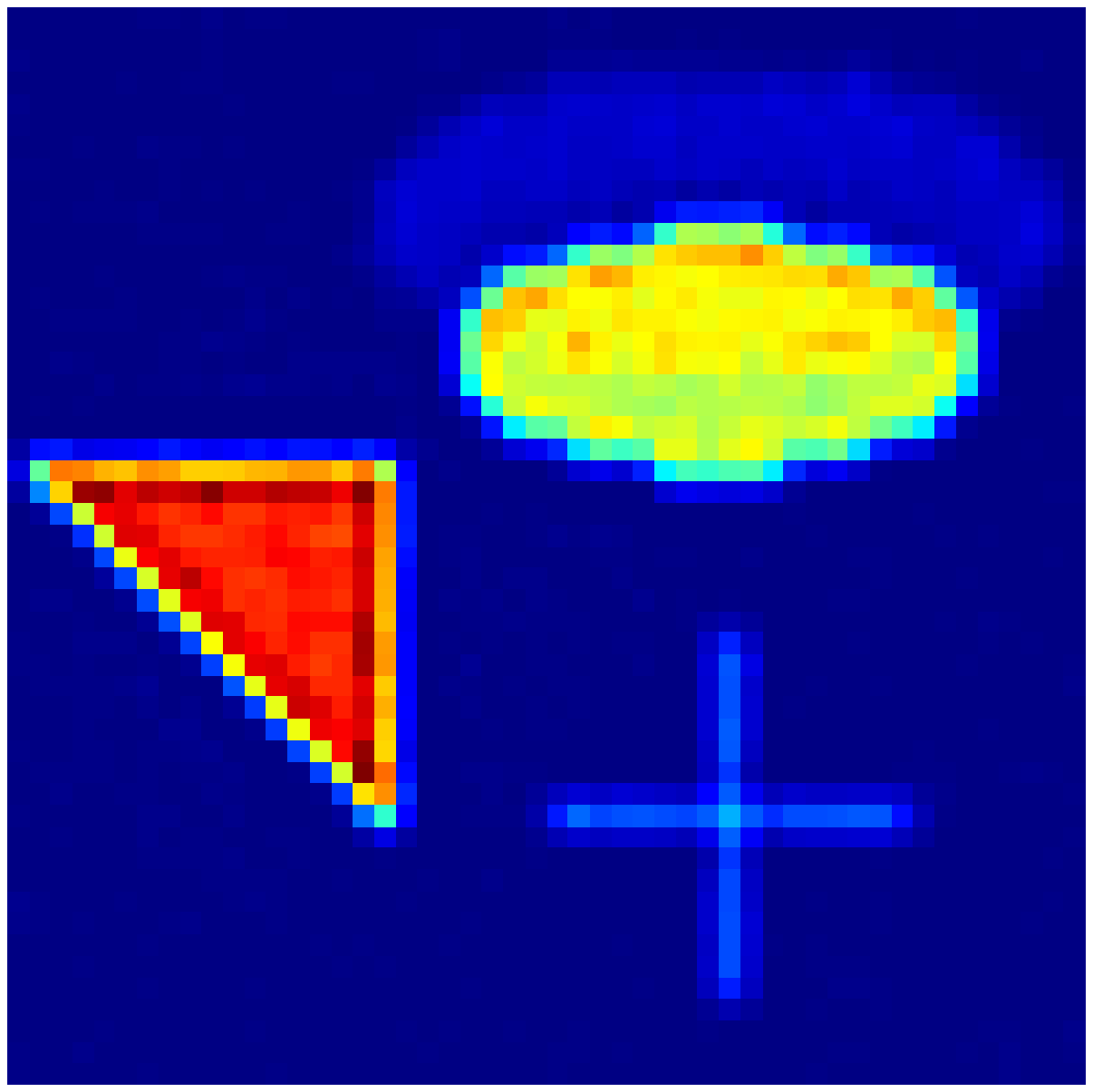} \hspace{-0.5cm} & \includegraphics[width=1.5cm]{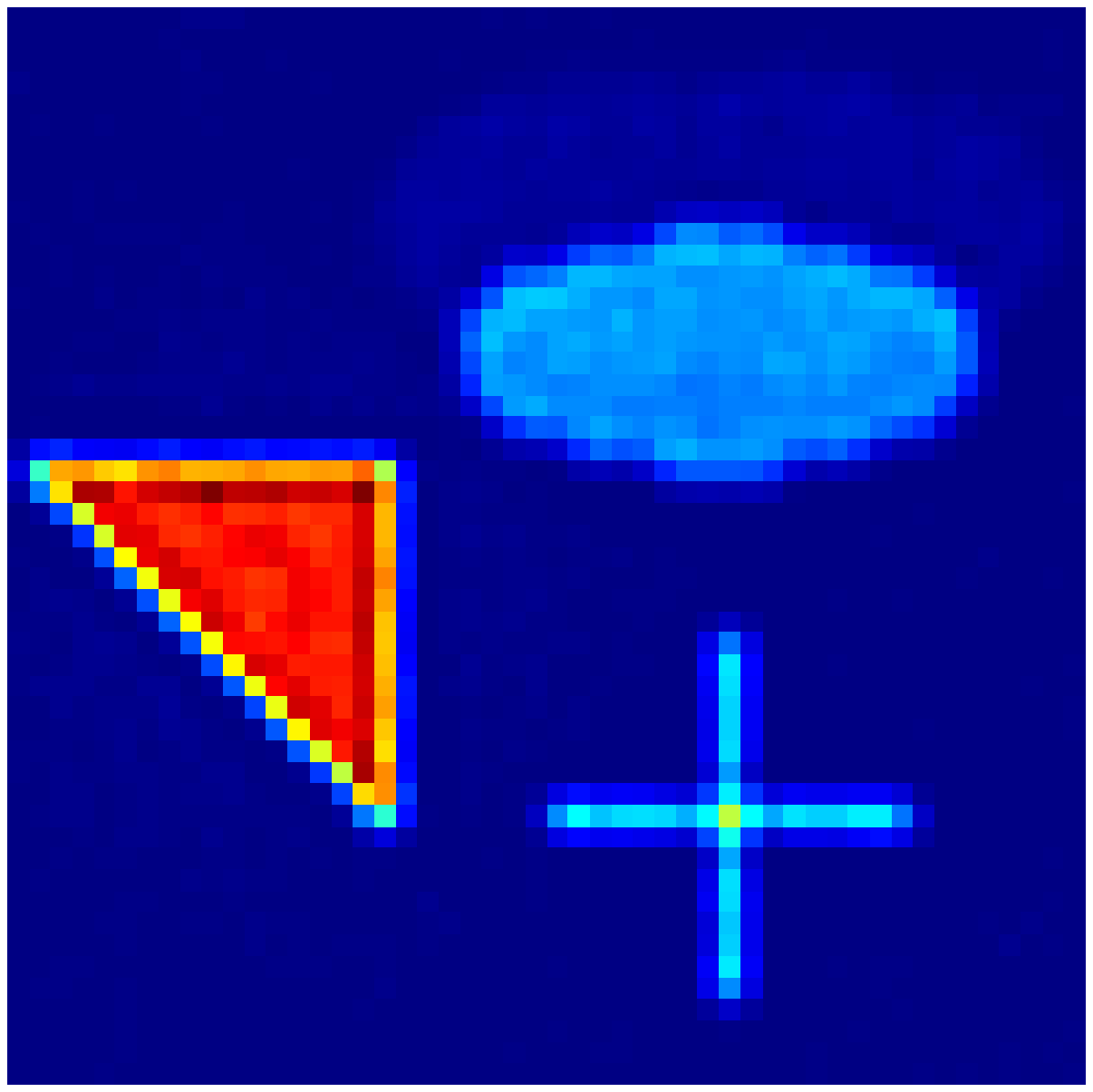} \hspace{-0.5cm} &  \includegraphics[width=1.5cm]{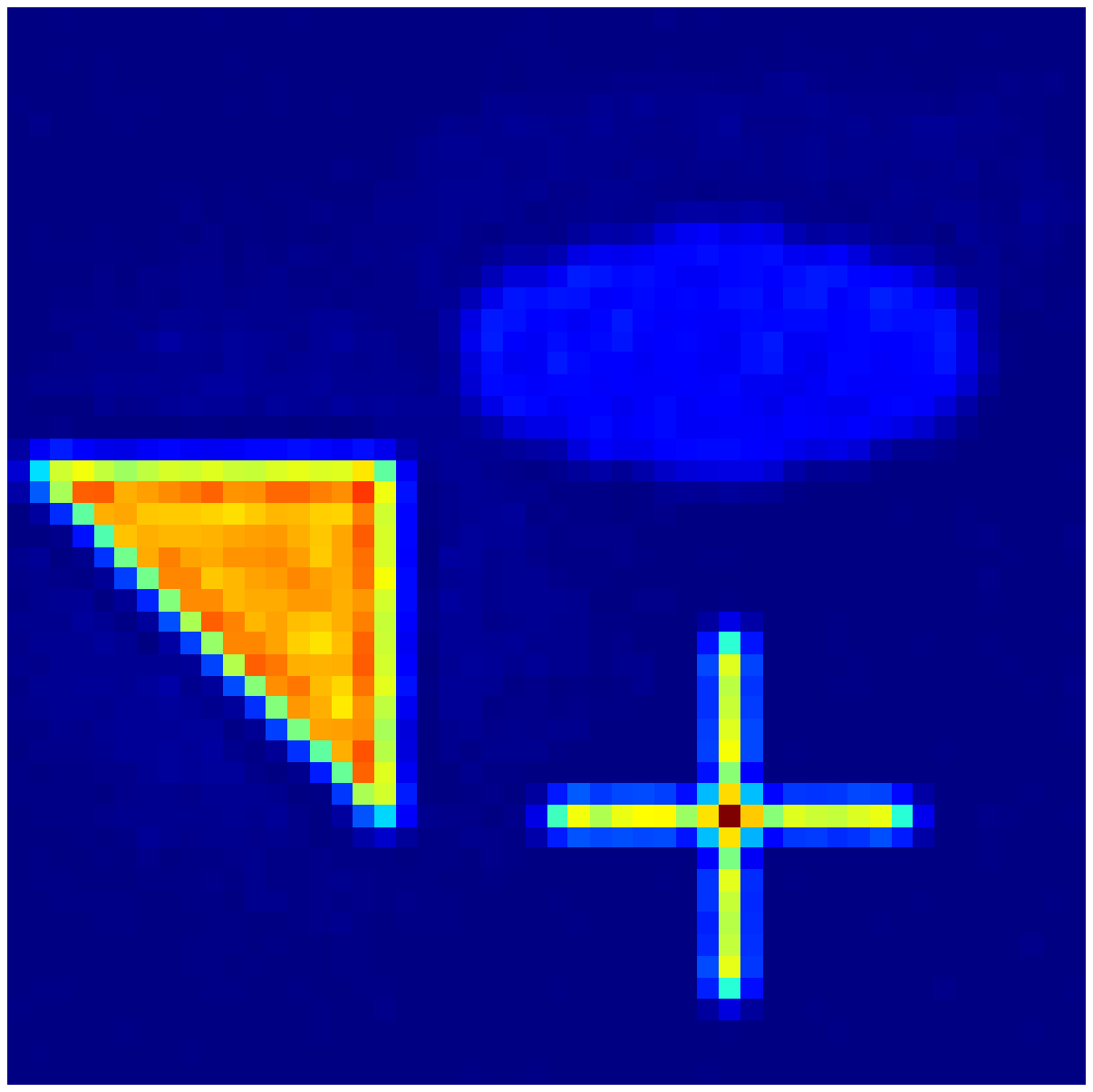} \hspace{-0.5cm} & 
\includegraphics[width=1.5cm]{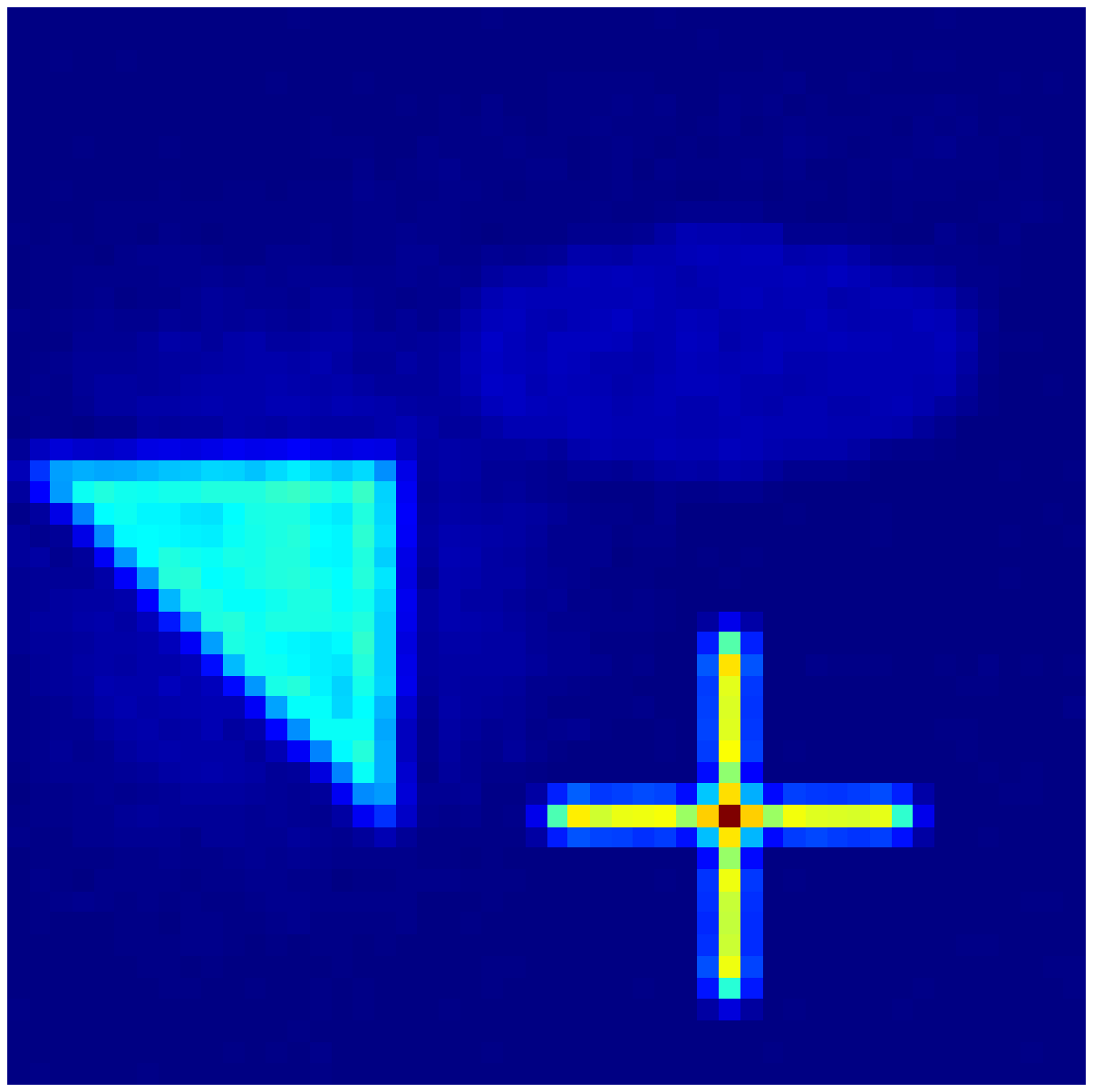} \hspace{-0.5cm} & 
\includegraphics[width=1.5cm]{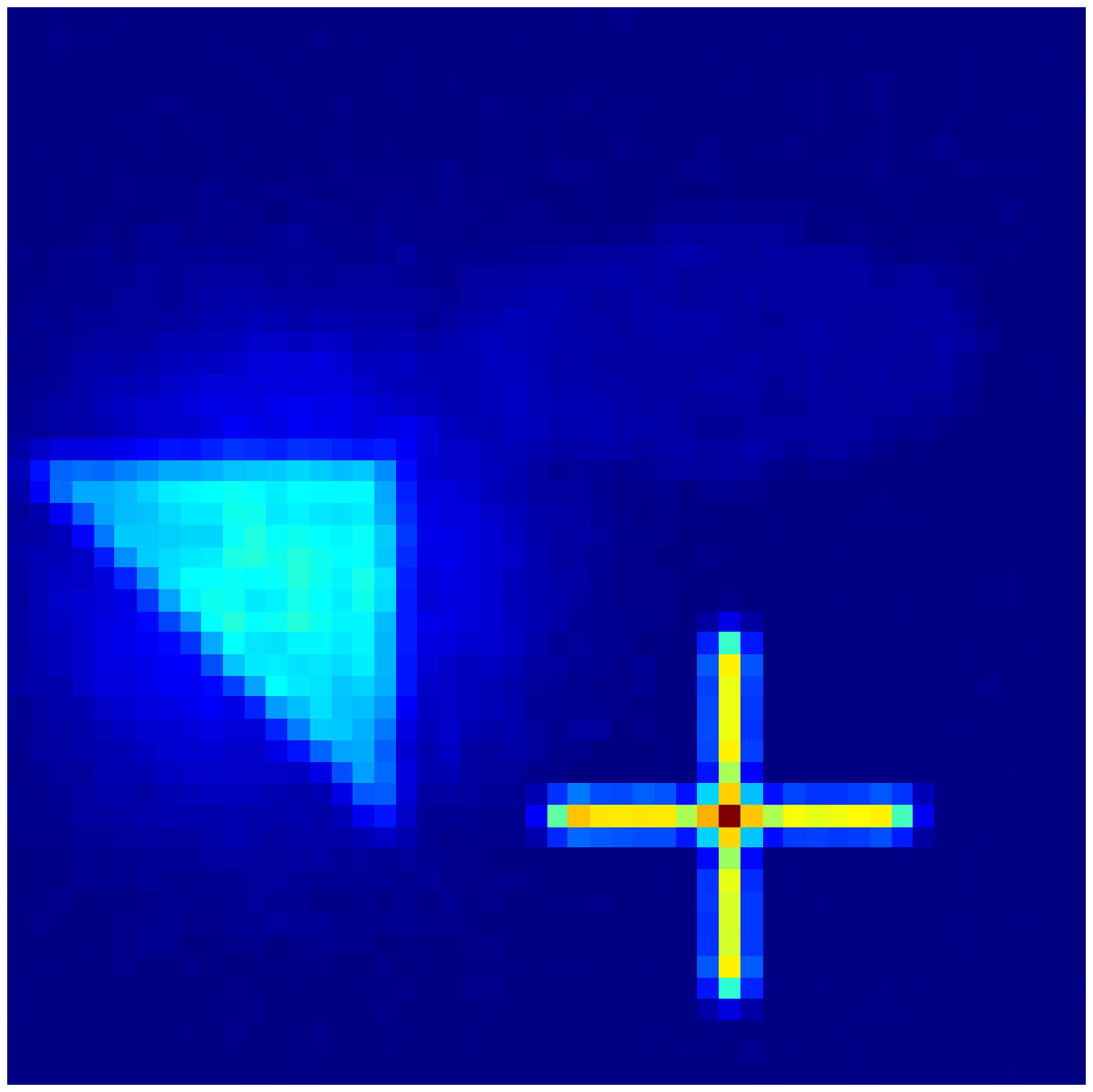} \hspace{-0.5cm} & \includegraphics[width=1.5cm]{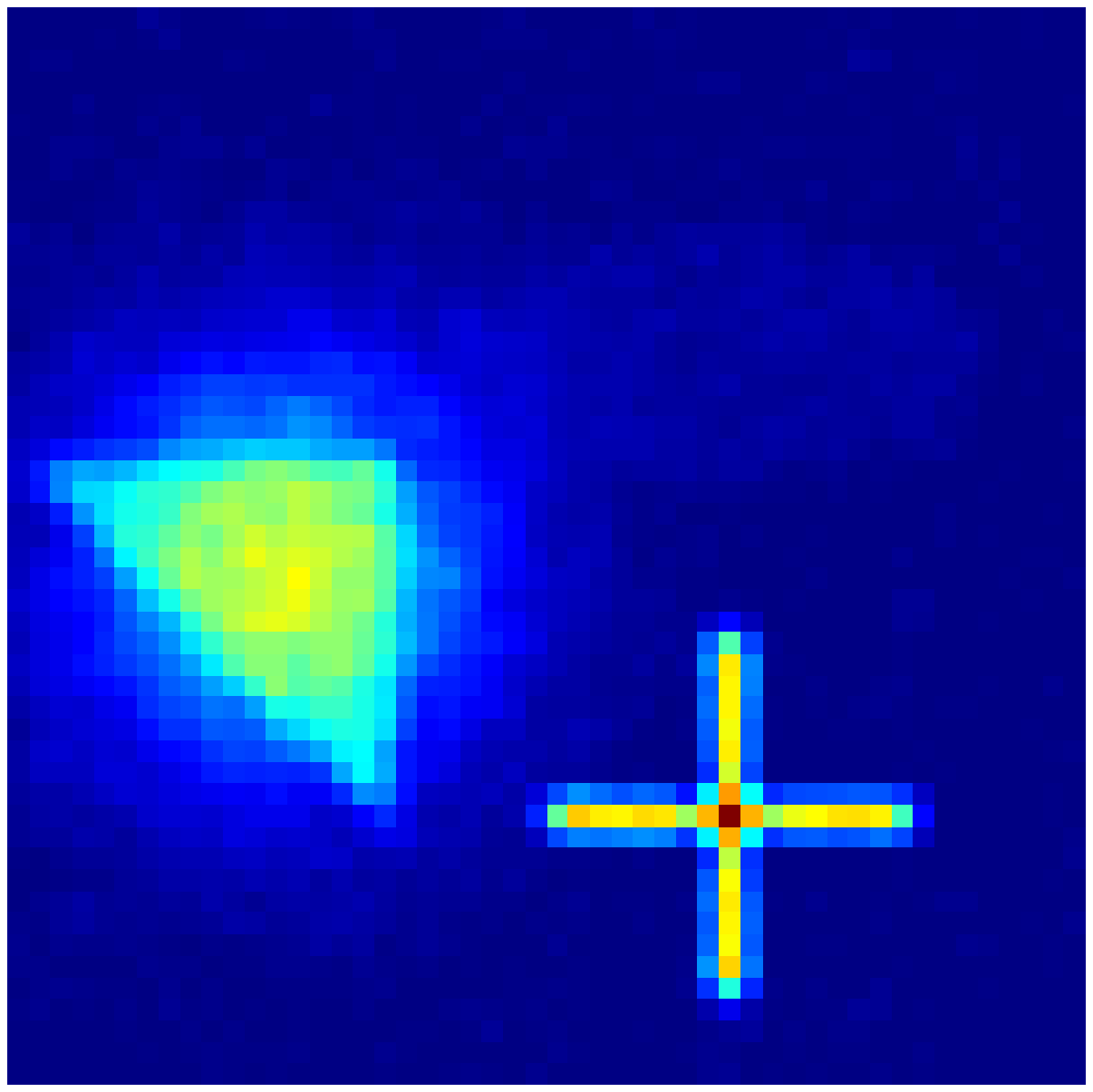}  \\ 
\raisebox{1.5\normalbaselineskip}[0pt][0pt]{\rotatebox[origin=c]{90}{\footnotesize hybrid}}  
\raisebox{1.5\normalbaselineskip}[0pt][0pt]{\rotatebox[origin=c]{90}{\footnotesize FGMRES-C}}  
\includegraphics[width=1.5cm]{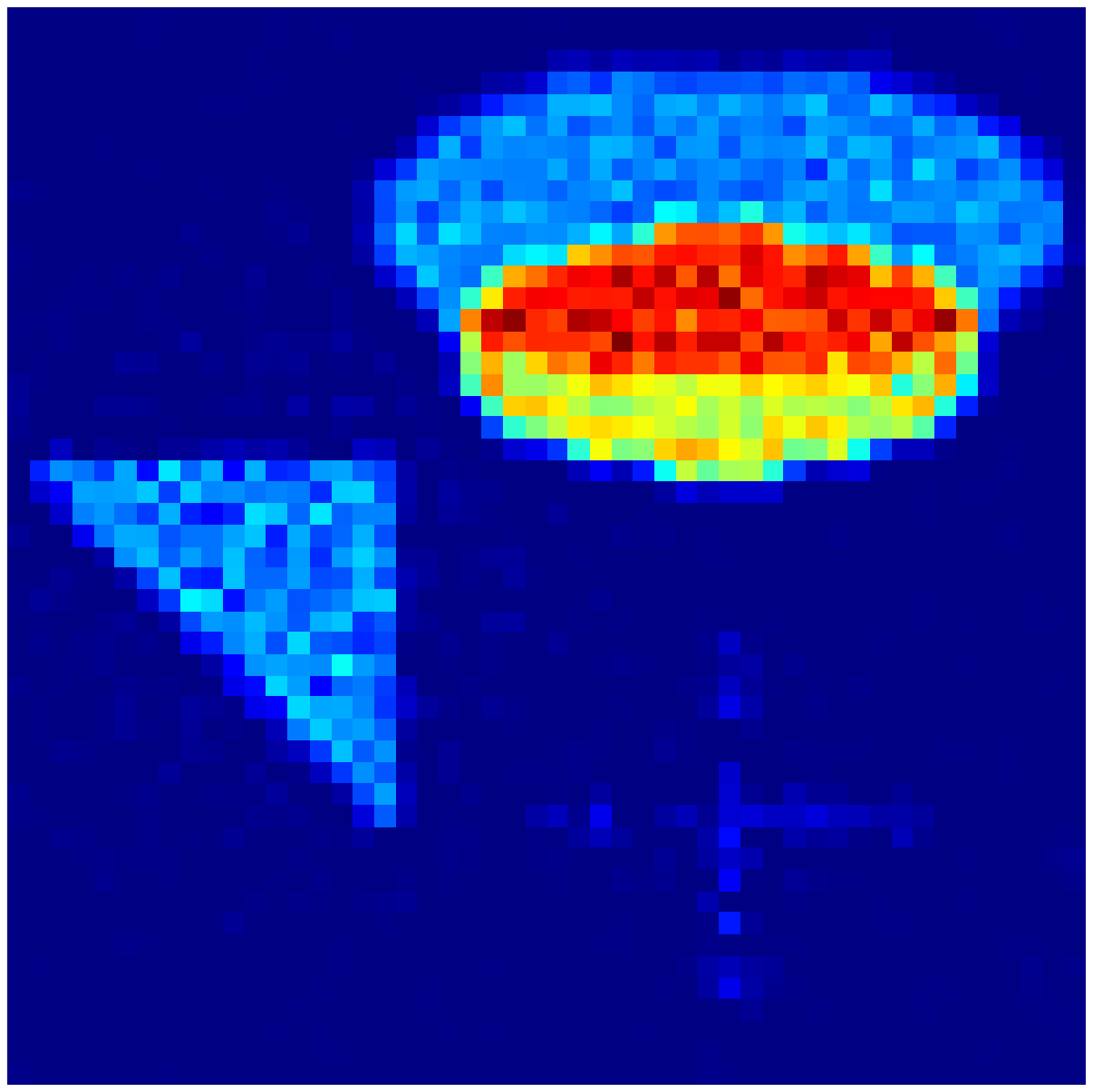} \hspace{-0.5cm} & \includegraphics[width=1.5cm]{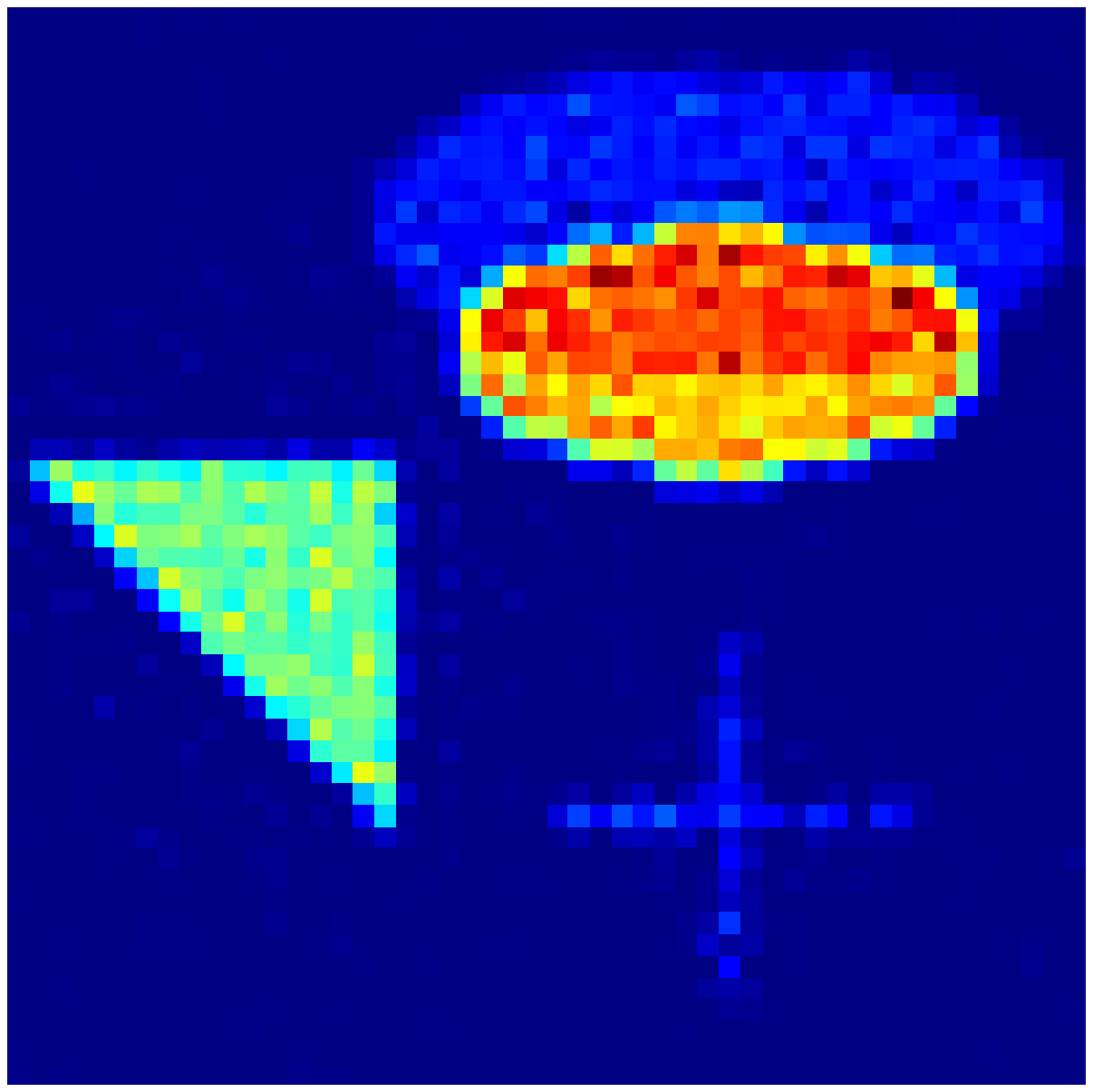} \hspace{-0.5cm} & \includegraphics[width=1.5cm]{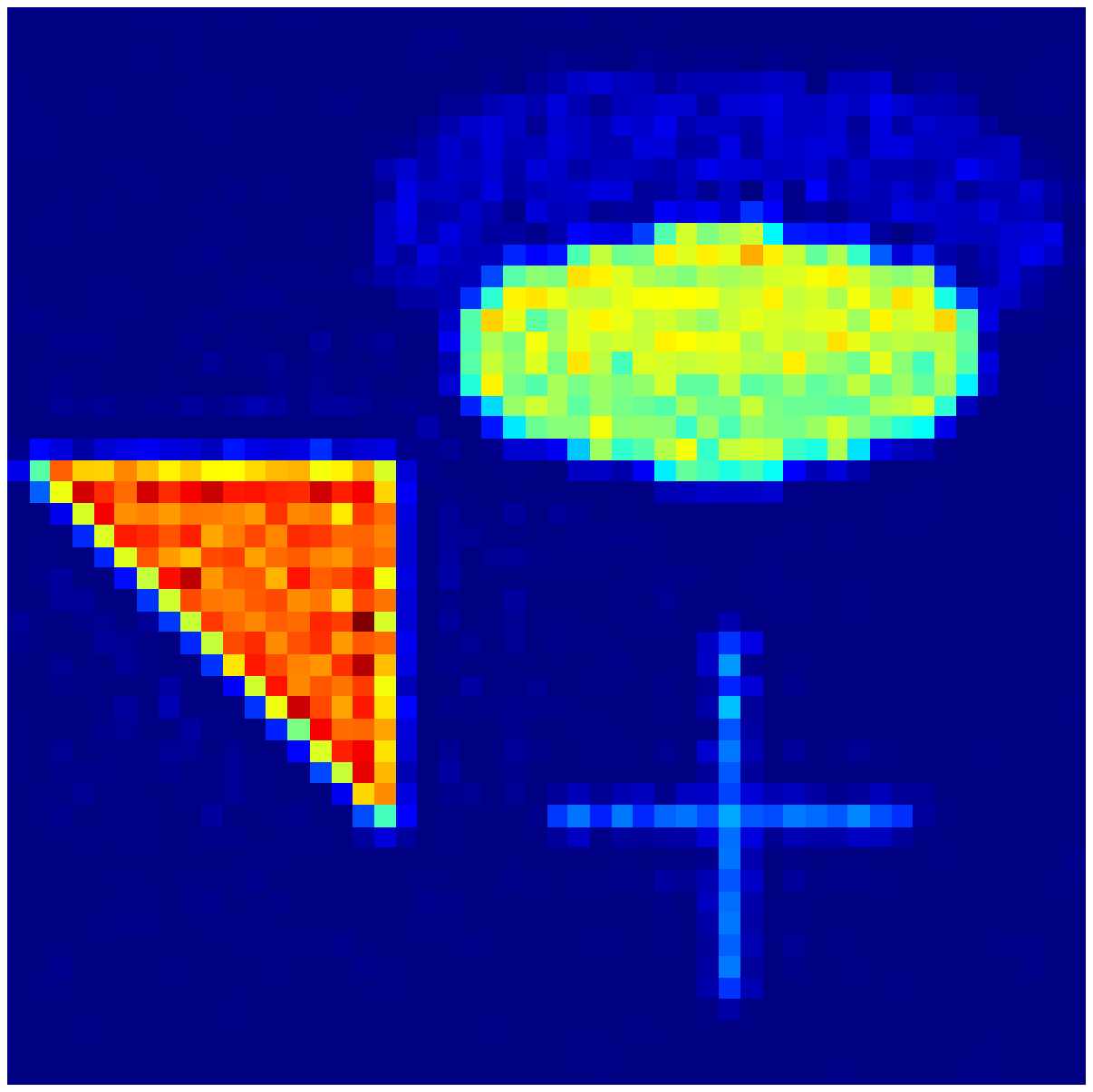} \hspace{-0.5cm} & \includegraphics[width=1.5cm]{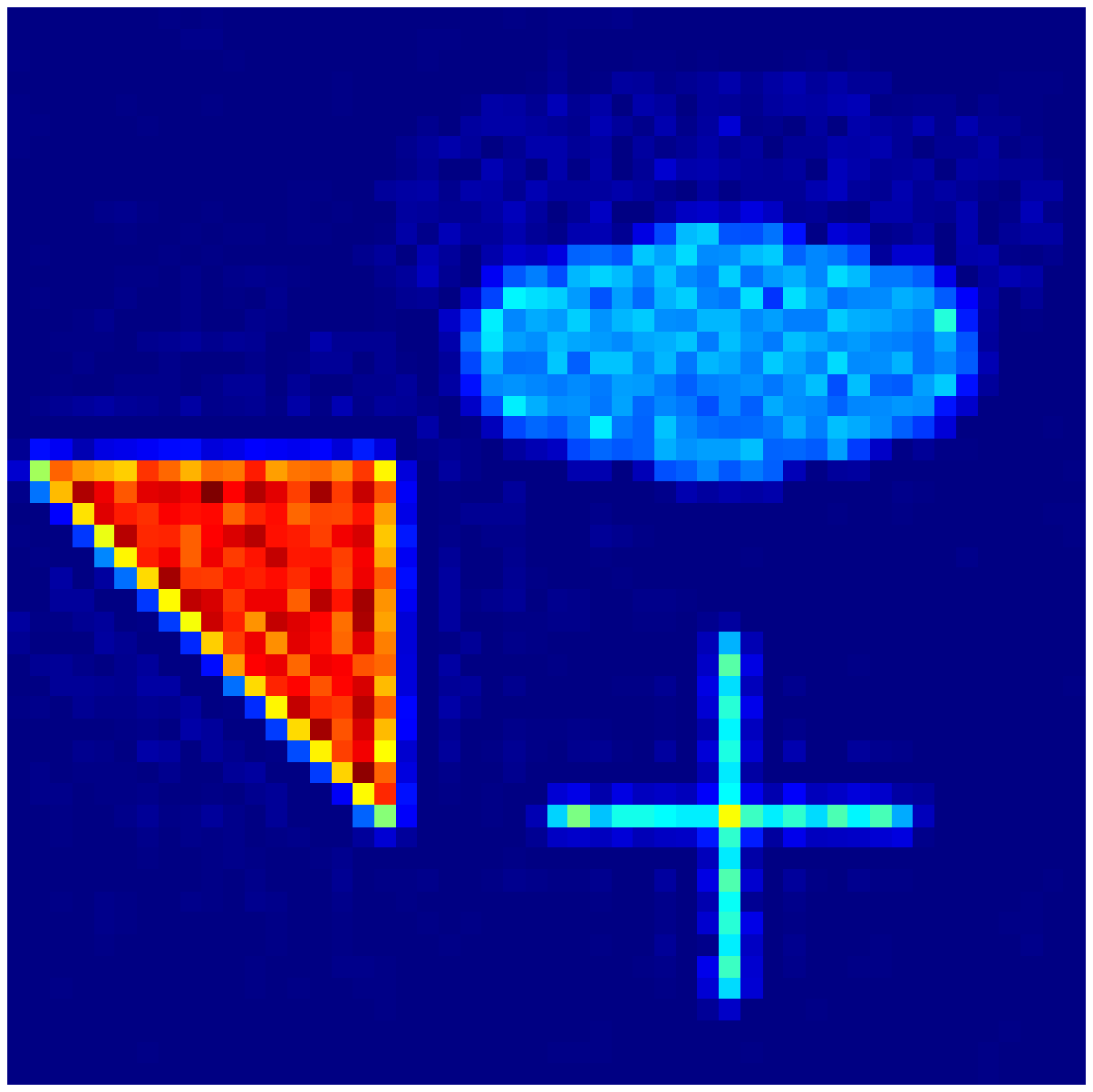} \hspace{-0.5cm} &  \includegraphics[width=1.5cm]{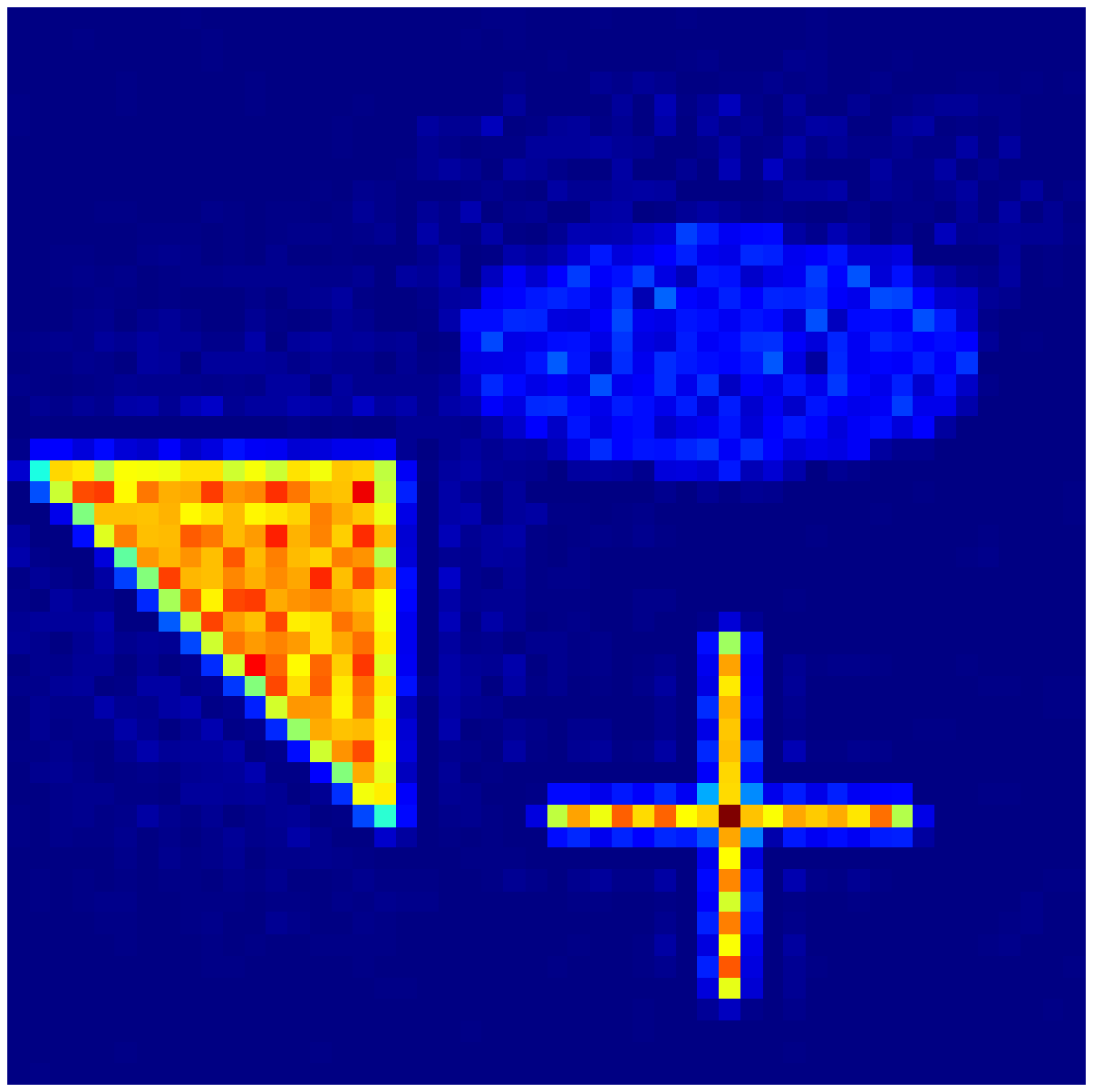} \hspace{-0.5cm} & 
\includegraphics[width=1.5cm]{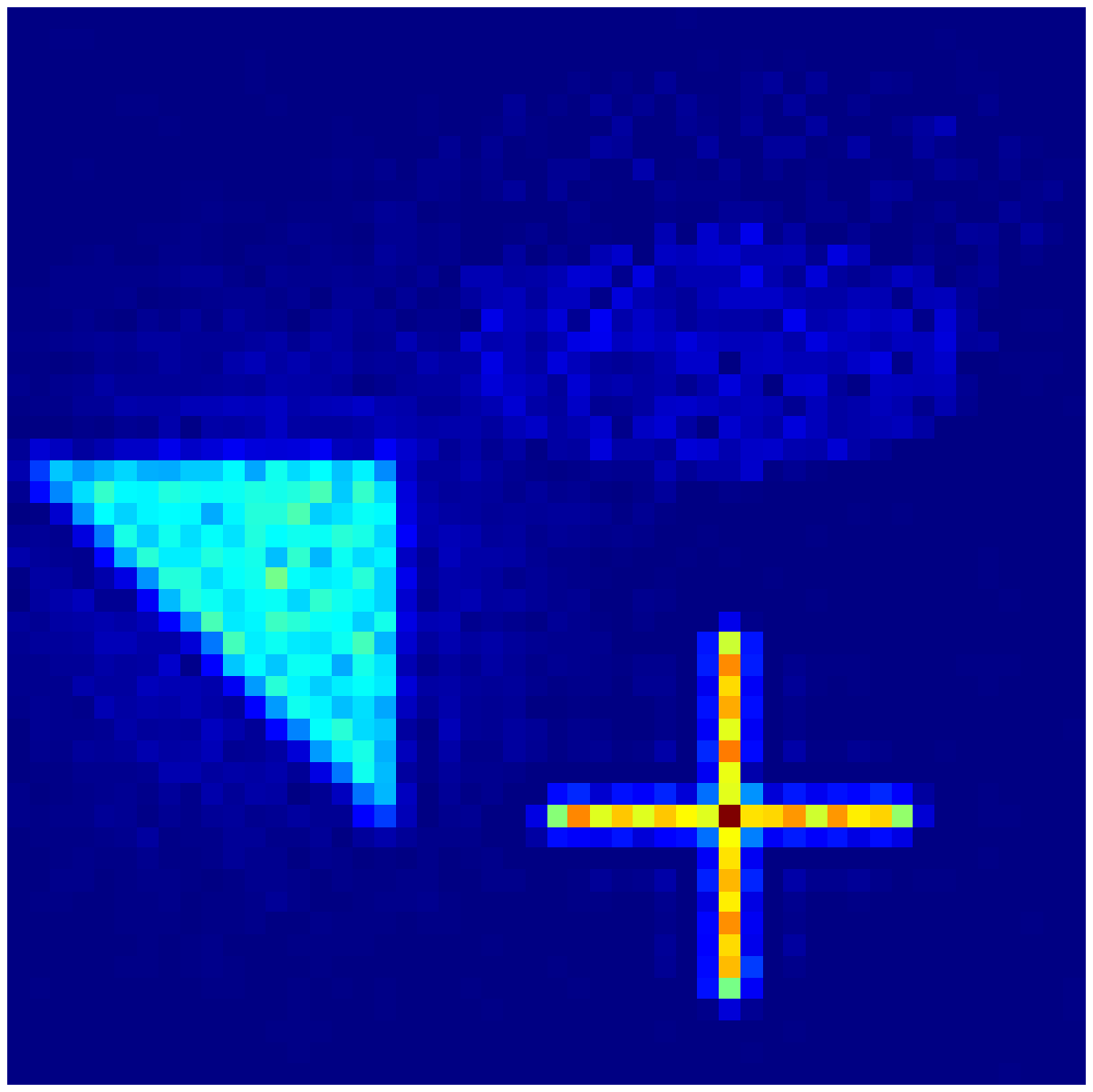} \hspace{-0.5cm} & 
\includegraphics[width=1.5cm]{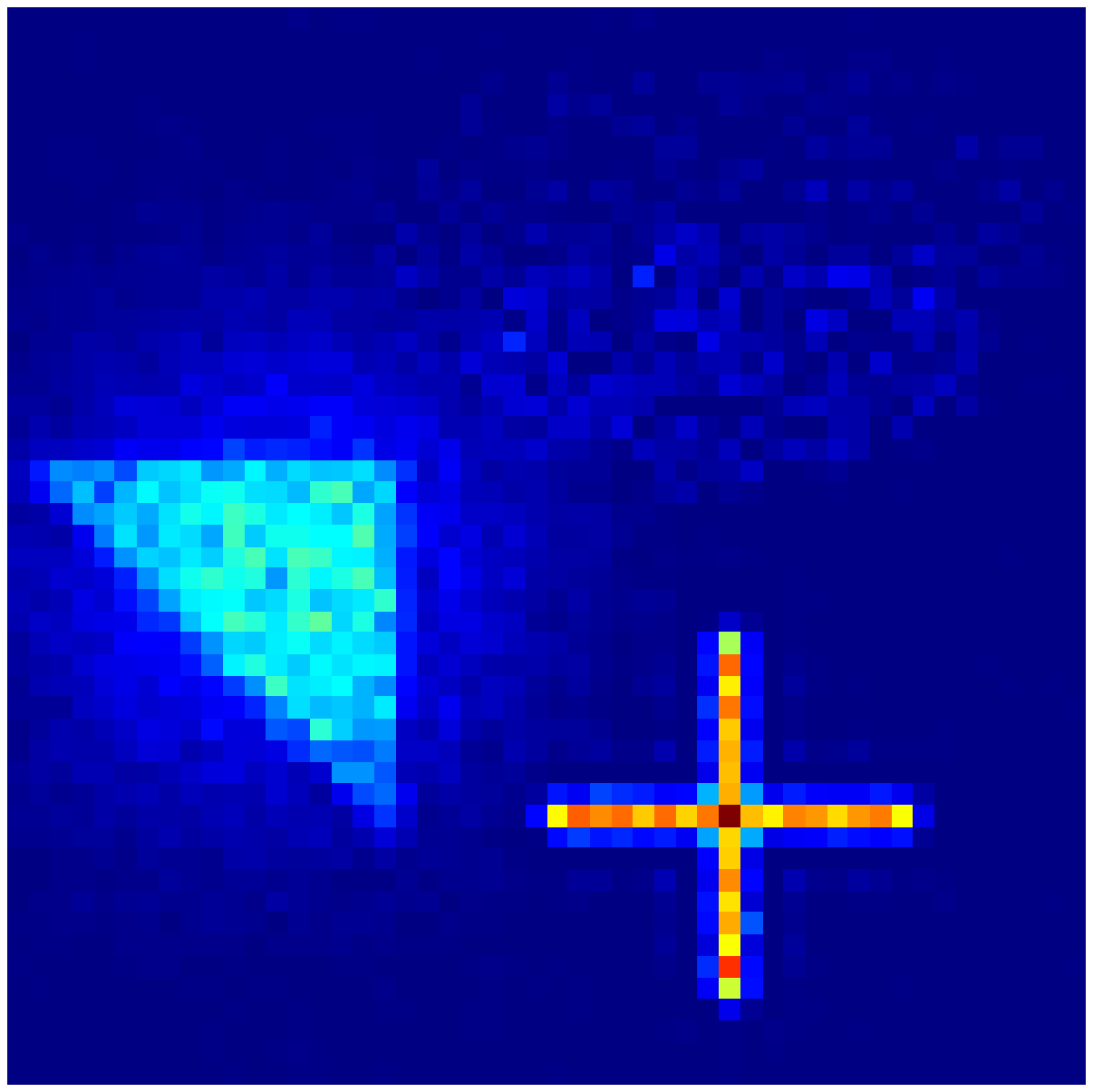} \hspace{-0.5cm} & \includegraphics[width=1.5cm]{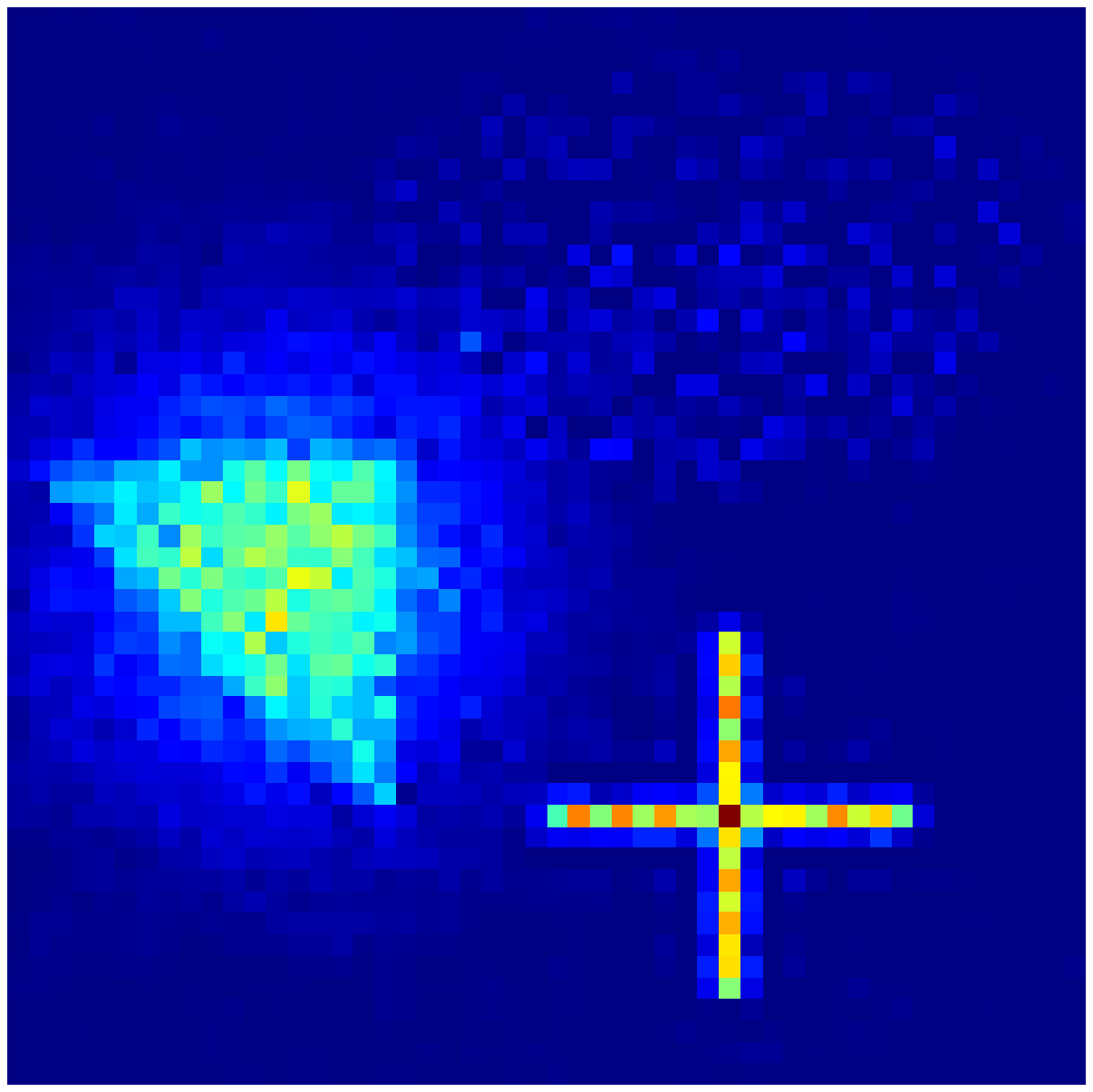}  
\end{tabular}
\caption{Reconstructions with $\ell_1$ and with a combination of $\ell_1$ and $\ell_{2,1}$ regularization for the dynamic image deblurring problem.  These methods are based on FGMRES, and the first time point has been omitted for visualization ease.}\label{fig:DB_recon2}
\end{center}
\end{figure}

One can also observe that for this example the methods based on FLSQR (Figure \ref{fig:errors_deblur_dynamic_LSQR}) perform better than the methods based on FGMRES (Figure \ref{fig:errors_deblur_dynamic_GMRES}), at the computational cost of an extra matrix-vector-product with the adjoint $\bfA\t$ per iteration. Nevertheless, we provide results for both solvers to highlight the benefit of using group sparsity (and a combination of sparsity and group sparsity) with respect to other regularizers, and we reiterate that this can be done using different flexible Krylov methods seamlessly.
Reconstructed images corresponding to $\ell_1$ regularization (top row) and a combination of $\ell_1$ and $\ell_{2,1}$ regularization (bottom row) are presented for methods based on LSQR in Figure \ref{fig:DB_recon1} and for methods based on GMRES in Figure \ref{fig:DB_recon2}. In particular, we would like to draw attention to the cross in the bottom-right of the images at the later time points (images on the far right). One can see that the reconstructions using $\ell_1$ regularization (top row) are more blurred that the ones using a combination of $\ell_1$ and $\ell_{2,1}$ (bottom row), which appear much crisper and closer to the true solution in Figure \ref{fig:DB_problem} (top row) both for methods based on FLSQR and on FGMRES.

\subsection{Anomaly detection in atmospheric inverse modeling}
\label{sec:OCO2}
In this experiment, we are interested in the inverse problem of efficiently and accurately detecting anomalies from observed data in atmospheric inverse modeling. This is relevant, for example, in large-scale anomalous emissions of greenhouse gasses and air pollution detection. A major challenge for anomaly detection is that inverse models using standard priors are not able to capture both anomalies (e.g., natural gas blowouts) and smooth regions (e.g., broad-scale emissions) simultaneously. A new solution decomposition framework was described in \cite{chung2022hybrid} where the desired parameters are represented as a sum of different stochastic components, and different priors are used for each component.  In particular, a Gaussian prior was used for the smooth component of the reconstruction and a sparsity prior was used for capturing anomalies.  However, anomalies are typically sparsely distributed in space but consistent over time, and the previous framework cannot capture such phenomena.  Thus, we have extended the solution decomposition framework for group sparsity regularization and present the results here.

In this example we consider a realistic linear atmospheric inverse model, 
where the aim is to estimate CO$_2$ fluxes across North America at $3$-hourly temporal resolution over $41$ days (approximately $6$ weeks from late June through July 2015) and at $1^\circ\times1^\circ$ latitude-longitude spatial resolution.  This setup corresponds to $3,222$ unknowns per $3$-hour time interval.  In particular, $\bfx_{\rm true} \in \bbR^{328 \cdot 3222}$ corresponds to the CO$_2$ fluxes from NOAA's CarbonTracker product (version 2019b). 
Synthetic satellite observations contained in $\bfb \in \bbR^{19,156}$ are generated as in \eqref{eq:linear_pbm} and mimic those from NASA's OCO-2 satellite, where $\bfA$ simulates an atmospheric transport model and $\bfe$ is added Gaussian noise to represent measurement errors. The components of the noise $\bfe$ are considered to be uncorrelated, so the covariance matrix $\bfR$ is $\sigma^2 \bfI$, where $\sigma=1.1267$ has been chosen so that $\bfe$ has a noise level of $\sigma \|\bfe\|/\|\bfA \bfx_{\rm true}\| = 1$. We remark that although these noise levels seem high for classic inverse problems settings, they are realistic in real data inverse modeling studies using OCO-2 data \cite{miller2018characterizing,miller2020impact}. We refer the interested reader to \cite{MillerSaibaba2020,liu2021} for additional detail on the specifics of the problem setup.

Although a decomposition of $\bfx_{\rm true} = \bfs + \bfxi$ where $\bfs,\bfxi \in \bbR^{328 \cdot 3222}$ is not available, we observe that, similar to actual atmospheric models, the true fluxes contain a combination of large, sparsely distributed values which correspond to anomalies (e.g., fires, anthropogenic emissions, or anomalies in biospheric fluxes) and smooth, broad regions of surface fluxes with small-scale variability. 
For defining the prior for $\bfxi$, we follow similar approaches \cite{yadav2016statistical, MillerSaibaba2020} and consider prior covariance matrix, $\bfQ = \lambda^{-2} \bfQ_t \kron \bfQ_s$ where $\bfQ_t$ represents the temporal covariance and $\bfQ_s$ represents the spatial covariance in the fluxes.  These covariance matrices are defined by kernel functions \begin{align}
    k_t(d_t;\theta_t) & = \left\{ \begin{array}{ll}
         1-\frac{3}{2}\left(\frac{d_t}{\theta_t}\right)+\frac{1}{2}\left(\frac{d_t}{\theta_t}\right)^3  & \mbox{ if } d_t\leq \theta_t,\\
         0 & \mbox{ if } d_t > \theta_t,
    \end{array} \right. \\
    k_s(d_s;\theta_s) & = \left\{ \begin{array}{ll}
         1-\frac{3}{2}\left(\frac{d_s}{\theta_s}\right)+\frac{1}{2}\left(\frac{d_s}{\theta_s}\right)^3  & \mbox{ if } d_s\leq \theta_s,\\
         0 & \mbox{ if } d_s > \theta_s,
    \end{array} \right. 
\end{align} where $d_t$ is day difference between two unknowns, $d_s$ is spherical distance between two unknowns, and $\theta_t,\theta_g$ are kernel parameters. In this setting, we set $\theta_t=9.854$ and $\theta_s = 555.42$, as in \cite{MillerSaibaba2020}. 
For the group sparsity regularizer for $\bfs$, we define $3222$ groups with each group corresponding to a spatial location.  Thus, the groups are defined to include 41 days of 3-hourly time intervals. 

\begin{figure}[bt]
    \centering
    \includegraphics[height = 2.6 cm, width = \textwidth]{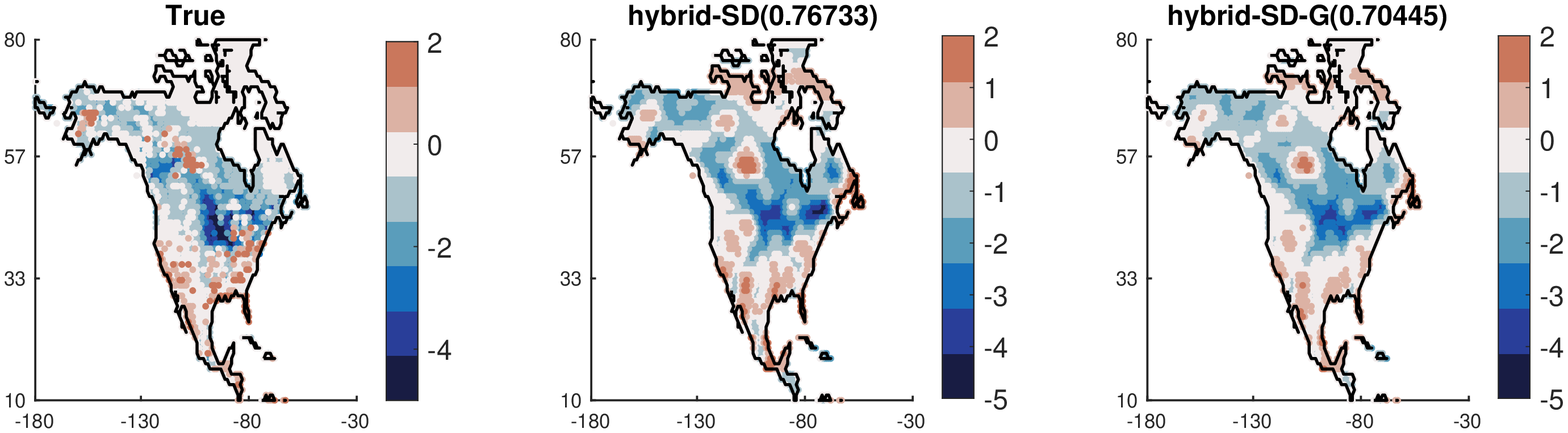}
    \caption{Anomaly detection OCO-2 example. We provide the averaged true fluxes for reference along with the averaged computed reconstruction of the fluxes (in $\mu$ mol $m^{-2} s^{-1}$) for hybrid-SD and hybrid-SD-G. Relative reconstruction error norms for the averaged spatio-temporal fluxes are provided in the titles, and all results correspond to using the DP selected regularization parameters.}
    \label{fig:ex3_recon}
\end{figure}

\begin{figure}[bt]
    \centering
    \includegraphics[height = 6.6 cm, width = .75\textwidth]{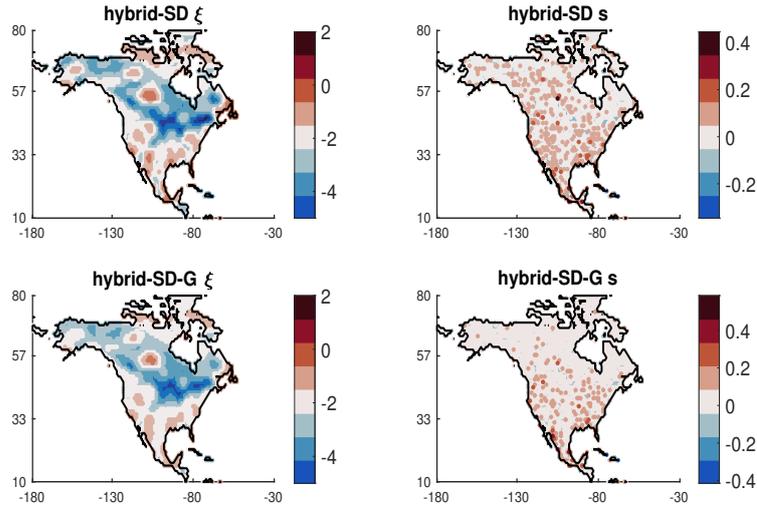}
    \caption{OCO-2 example.  Reconstructed $\bfxi$ and $\bfs$ images for the atmospheric inverse problem using hybrid-SD (top row) and with group sparsity hybrid-SD-G (bottom row).}
    \label{fig:ex3_split}
\end{figure}

We compute spatio-temporal reconstructions using hybrid-SD-G and compare the results to hybrid-SD (referred to as \texttt{sdHybr} in \cite{chung2022hybrid}), where the regularization parameters $\lambda$ and $\alpha$ are selected using the DP as defined in \eqref{dp}. The temporal-averaged images of the flux reconstructions are presented in Figure \ref{fig:ex3_recon}, along with the temporal average of the true image. We observe that both hybrid-SD-G and hybrid-SD average reconstructions are able to capture both sources and sinks present in the true average image, with the hybrid-SD-G reconstruction having a slightly smaller relative reconstruction error for the average image.  The main benefit of the solution decomposition framework is the ability to obtain two components of the solution. In Figure \ref{fig:ex3_split}, we provide the reconstructions of the individual components $\bfxi$ and $\bfs$ that form the solution for hybrid-SD (top row) and hybrid-SD-G (bottom row). It is clear that group sparsity regularization provides a smoother background and is able to distinguish persistent anomalies (as opposed to spurious false positive anomalies) better than standard sparsity regularization. 

The results of this case study demonstrate that hybrid-SD-G can yield accurate results for complex, spatio-temporal atmospheric inverse modeling that are inherently different from the ones obtained using hybrid-SD. Comparisons of hybrid-SD to existing reconstruction methods can be found in \cite{chung2022hybrid}.

\section{Conclusions and future work}\label{sec:conclusions}

This paper presents a suite of numerical algorithms based on flexible Krylov methods for solving linear inverse problems with group sparsity regularization. The idea of using group structure to achieve better sparse recovery has received much attention, but automated and efficient approaches for practical implementation are still lacking.  This work addressed that gap by proposing new iterative methods that are efficient, since they only require matrix-vector and possibly matrix-transpose-vector products, and automatic, by selecting regularization parameters on-the-fly. 
Moreover, these approaches exploit flexible preconditioning techniques to avoid inner-outer schemes by building a single solution subspace from which to compute solutions. 

Additional regularization can be added to the projected problem to avoid semi-convergence or to guarantee convergence of the proposed method (at a higher computational cost). This scheme allows for the regularization parameter to be chosen adaptively using suitable parameter choice criteria. In particular, the results in this paper are shown for hybrid-FGMRES, hybrid-FLSQR, IRW-FLSQR, and hybrid-SD, all with group sparsity regularization and using the DP to find appropriate regularization parameters.

Moreover, the scheme presented in this paper is highly general, since very different regularization terms can be modeled in the group sparsity framework. This is highlighted in the numerical examples, which show the performance of the method in different applications and for different group sparsity modalities. In particular, an example of overlapping group sparsity is shown using the natural tree-structure of wavelet decompositions and non-overlapping group sparsity is tested in solutions that have spatio-temporal components where the solution is temporally persistent and sparse in space.

Future work includes extensions to other applications and other constraints (e.g., nonnegativity).  For example, group sparsity regularization has been used for nonnegative matrix factorization \cite{kim2012group} and for group-based dictionaries in neuroimaging using fMRI \cite{zhao2018functional}. 



\backmatter



\bmhead{Acknowledgments}
The authors gratefully acknowledge support from the Cambridge Mathematics of Information in Healthcare Hub (CMIH), University of Cambridge.  The authors would also like to thank Fred Vickers-Hastings and Pol del Aguila Pla for early discussions on this work.

This work was partially supported by the National Science Foundation program under grants DMS-2026841, DMS-2245192, and DMS-2208294.  Any opinions, findings, and conclusions or recommendations expressed in this material are those of the author(s) and do not necessarily reflect the views of the National Science Foundation.


\section*{Statements and declarations}


\bmhead{Conflict of interest} The authors have no competing interests to declare that are relevant to the content of this article.


\bibliography{references}

\end{document}